\documentclass[a4paper,11pt]{article}
\usepackage[utf8]{inputenc}
\usepackage{float}
\usepackage{graphicx}
\usepackage{amsthm}
\usepackage{latexsym}
\usepackage{amsmath}
\usepackage{amssymb}
\usepackage{amsfonts}
\usepackage{multirow}
\usepackage{amsthm}

\usepackage{enumerate}
\usepackage{graphics}
\usepackage{fullpage}
\usepackage{color}
\usepackage{comment}

\usepackage{bussproofs}

\usepackage{booktabs, tabularx}

\usepackage{xspace}
\usepackage{xcolor,colortbl}

\usepackage{mathrsfs}

\usepackage{tikz}
\usetikzlibrary{trees,arrows,automata,shapes,decorations}

\newtheorem{thm}{Theorem}[section]
\newtheorem{defi}[thm]{Definition}
\newtheorem{lem}[thm]{Lemma}
\newtheorem{cor}[thm]{Corollary}
\newtheorem{rem}[thm]{Remark}
\newtheorem{exam}[thm]{Example}
\newtheorem{prop}[thm]{Proposition}

\newcommand{\eps}{\varepsilon}

\newcommand{\ot}{\leftarrow}

\newcommand{\RN}[1]{%
  \textup{\uppercase\expandafter{\romannumeral#1}}%
}

\newcommand{\bM}{{\mathbf M}}
\newcommand{\mb}{\mathbf}

\newcommand{\myComment}[1]{}

\title{Knowability as continuity:\\
a topological account of informational dependence}
\author{Alexandru Baltag \& Johan van Benthem}
\date{}
\begin{document}
\maketitle

\begin{abstract}
We study knowable informational dependence between empirical questions, modeled as continuous functional dependence between variables in a topological
setting. We also investigate epistemic  independence in topological terms and show that it is compatible with functional (but non-continuous) dependence. We then proceed to study a stronger notion of knowability based on uniformly continuous dependence. On the technical logical side, we determine the complete logics of languages that combine general functional dependence, continuous dependence, and uniformly continuous dependence.
\end{abstract}
\section{Introduction and preview}

\paragraph{Dependence in the world} Dependence occurs when objects or phenomena constrain each other's behavior. These constraints  range from complete functional determination to milder forms of correlation.  At the opposite extreme lies independence:  complete freedom of behavior, regardless of what the other object or phenomenon does. A common  mathematical model for all this uses a family of  variables $X$ for relevant items that can take values in their ranges, where correlations arise when not all a priori possible value combinations can occur. The laws of physics provide many examples: Newton's law of gravity constrains the values that can arise for the masses of two objects, their distances, and their accelerations toward each other. Dependence is at the heart of the physical universe and our theories about it, but just as well our social world and the many dependencies in human behaviors.  Logicians have long been interested in modeling this pervasive notion and studying its basic laws \cite{Vaa07}, \cite{BaBe21}. In this line, the present paper offers a new dependence logic in a richer topological setting that applies to empirical inquiry.

\vspace{-2ex}

\paragraph{Dependence, knowledge and information flow} In a complementary sense to its physical manifestation, dependence can also be seen as an information-carrying relation. If $y$ depends functionally on $X$, then, given the information which values  the variables in $X$ have right now, we automatically know the value of $y$. And weaker forms of correlation license weaker forms of knowledge transfer. This dual perspective was
proposed in Situation Theory, where realism about constraints between situations in the world \cite{BaPe} led to logical theories of information flow for human agents \cite{BaSe}. But dependence also underlies (Dynamic) Epistemic Logic, the basic theory of semantic information \cite{BaRe16}, \cite{HandEL}, \cite{Be11}. If some players get a card each from a  set whose composition is common knowledge, then the card dealt to the first player constrains the card combinations available for the other players. Or epistemically, if we know the card of the first player, we also know more about which cards are held by which other player.
In this paper, we will mainly pursue this informational-epistemic perspective on dependence.

\vspace{-2ex}

\paragraph{Variables as {\it wh}-questions} As used in empirical sciences, `variables' differ in an essential way
from the First-Order Logic notion with the same name: while in FOL variables are simple placeholders, with no meaning of their own, in empirical science they are specific ``what" \emph{questions} (e.g. what is the position of a particle, or its velocity, momentum or acceleration, etc?). The values of a variable are the possible exact answers to the question (e.g. the exact position, etc). A variable $y$ will take different values in different possible states of the world, so in this sense variables are \emph{functions} from a state space $S$ to a range of values $\mathbb{D}_y$. On the other hand, the corresponding \emph{propositional question} can be understood as a \emph{partition} of the state space: each cell of the partition represents (the set of worlds satisfying) a specific answer to the question. We obtain a correspondence between variables and (propositional) questions: any variable $y: S\to \mathbb{D}_y$ induces a partition on $S$, whose cells are given by the sets $y^{-1}(s):=\{w\in S: y(w)=y(s)\}$, for all $s\in S$.

\vspace{-2ex}

\paragraph{Functional dependence and propositional determination}
Unlike in FOL, in empirical sciences  the values of different variables are not necessarily  independent: answering some questions may already be enough to answer other questions. In this sense, \emph{dependence is all about logical relations between questions}.
The standard logical setting for this is provided by the notion of \emph{functional dependence} between variables/questions, expressed as a \emph{global} statement $=(X;y)$ in \cite{Vaa07}, and denoted in this paper by $D(X; y)$: essentially, this captures the existence of a function $F$ that maps the values of the variables in $X$ at all states into the value of $y$ at the same states. There is also a \emph{local}, state-dependent version $D_Xy$, introduced in \cite{BaBe21}, which can be thought of as semantic \emph{determination}: fixing the $X$-values to the current ones (realized at the current state $s$) also fixes the value for $y$. In the epistemic reading, a dependence statement $D_X y$ can be seen as a \emph{conditional knowledge} claim: an observer can know the current value of $y$ if given the current values of the variables in $X$. As a special case, when $y$ is a Boolean variable associated to some proposition $\varphi$, we obtain a natural `determination modality' $D_X\varphi$, which says that fixing the $X$-values to the current ones (at $s$) `determines' the truth of $\varphi$ (i.e. fixes its value to `true').

The complete logic of these notions was investigated in \cite{BaBe21}. Our concern in this paper is lifting an important idealization in the preceding framework, which is necessary to make dependence logics function in a setting of empirical investigation.

\vspace{-2ex}

\paragraph{From idealized sharp values to empirical measurement} Standard scenarios in epistemic logic are discrete, and it then makes sense to work with \emph{sharp values} of variables, such as truth values for proposition letters. However, in the empirical sciences, and in daily life, values are usually determined by imprecise observations or \emph{measurements}. Even if not in error, these deliver only intervals or other sets of values containing the real value.
Of course, one can perform more measurements and combine outcomes, or turn to a more precise measuring device, but sharp individual values for variables will still not be delivered. Thus, in empirical inquiry, the   view of dependence as conditional knowledge of sharp values becomes unrealistic.

Here we arrive at the main question of this paper. How should we model dependence and its logic in an empirical setting where only measurements with limited precision are available? We need a richer style of modeling for \emph{outcomes of observations} or measurements, and we need a way of describing how empirical knowledge can arise though a  process of \emph{approximation}.

\vspace{-2ex}
\paragraph{Topological models for empirical inquiry and knowledge} As it happens, a suitable  richer style of modeling exists. \emph{Topological models} have been around since the 1930s for modeling intuitionistic and modal logics \cite{BeBe07}, where the surplus of topological spaces over plain sets  elucidates epistemic notions such as knowledge or verifiability. In our present setting, open sets in a topological space represent  possible outcomes of experiments or measurements, while propositions of interest are modeled by arbitrary subsets. In this way, topological notions acquire epistemic meaning. For instance, a property $P$ denoting an open set is `verifiable' in that, for any point $s$ satisfying $P$, there exists an experiment that can tell us that the point has that property. Dually, properties denoting closed sets are falsifiable. In this setting, even technical  separation properties for topological spaces  make epistemic sense: e.g., the $T_0$-property says that any  point can be identified uniquely by experiments, since  any other point can be separated by some measurement.  The  area of Epistemic Topology exploits all these features, witness important studies like \cite{KK96}, \cite{DaPaMo96}, \cite{BalGierSmets}. Our analysis of dependence follows up on this, by looking at further structure in dependence models with topologies over the value ranges for the variables.
\vspace{-2ex}
\paragraph{Empirical variables as topological maps}
We are thus lead to a natural modeling of empirical variables $y$ as functions from a state space $S$ onto a \emph{topological space} $(\mathbb{D}_y, \tau_y)$: in this way, in addition to the exact value $y(s)$ of variable $y$ at state $s$ (an idealization which typically is not accessible by measurements), we also have \emph{observable approximations} of this value, given by open neighborhoods $U\in \tau_y$ with $y(s)\in U$.
Finally, topological spaces allow us to treat \emph{sets of variables} $X$ as single new entities taking tuples as values, by using the product topology on the Cartesian product $\Pi_{x\in X} \mathbb{D}_x$ value spaces for the separate $x \in X$. This uniform treatment of single variables and finite sets of these greatly simplifies epistemic analysis.

\vspace{-2ex}

\paragraph{Propositionalization: empirical questions as topologies} If we move from empirical variables to the corresponding propositional questions, this requires refining the partitional view of questions into a \emph{topology} on the state space $S$; any open neighborhood of the current state $s$ in the question's topology will represent a partial, ``inexact" answer to the given question: these are the ``feasible" answers, that one can determine by observations (measurements).\footnote{The underlying partition of $S$ into cells (representing the `exact' answers) can still be recovered from this topological representation: the cell of state $s$ consists of all states that are ``topologically indistinguishable" from $s$, i.e. that belong to exactly the same open sets as $s$.} The above-mentioned correspondence between variables and questions-as-partitions extends to their topological versions: the topology $\tau_y$ on the range of values $\mathbb{D}_y$ of variable $y$ naturally induces a topology on the state domain $S$, having open sets of the form $\mathsf{y}^{-1}[U]$, for any open $U\in\tau_y$.

\vspace{-2ex}
\paragraph{Continuity as information-carrying dependence} We have presented dependence $D_Xy$ as delivering conditional knowledge of $y$ given $X$. But note that this required being given the exact value(s) of $X$. In an empirical setting such precision is seldom  achievable: all we can observe are measurable approximations. Thus, a useful and truly `epistemic' dependence should allow us to \emph{approximate} the value of $y$ arbitrarily closely given suitable approximations for the values of $X$. In its global form, this requires that the dependence function $F$ mapping $X$-values to $y$-values be \emph{globally continuous} with respect to the underlying topologies. This is the definition of our \emph{known epistemic dependence} $K(X;y)$: the observer \emph{knows} that she can get to know the value of $y$ as accurately as needed if given accurate enough approximations of $X$. What makes this dependence ``epistemic" is the fact that it requires only observable approximations, and what makes it ``known" is its globality. As in the case of exact dependence, there is also a \emph{local} version $K_Xy$ of epistemic dependence, which only requires that the dependence function $F$ is \emph{continuous on some open neighborhood} $U$ of the exact value(s) $X(s)$ (of $X$ at the current state $s$). We read $K_Xy$ as \emph{knowable} epistemic dependence: the observer \emph{can come to know} (after doing an accurate enough measurement of $X$) that she can get to know the value of $y$ as accurately as needed if given accurate enough approximations of $X$.\footnote{There exists an even ``more local" version $k_Xy$ of epistemic dependence, requiring only point-continuity at $X(s)$. But this is in a sense ``too local", too sensitive to the actual value of $X$, at the risk of being ``unknowable": $k_Xy$ may hold without ever being known to the observer, no matter how accurate her observations. In epistemological terms, point-continuity is not an `introspective' fact.}


\paragraph{Special case: conditional knowability of a proposition} As in the case of exact dependence, we can apply epistemic dependence $K_Xy$ to the case when $y$ is a Boolean variable associated to a proposition $\varphi$, obtaining a natural notion of \emph{knowability of a proposition given (accurate enough approximations of) $X$}: $K_X\varphi$ holds at $s$ if there is if there is some $X$-open neighborhood of $s$ (in the topology induced by $X$ on the state space) all of whose members satisfy $\varphi$. This means that one can come to know that $\varphi$ was the case by measuring $X$ with enough accuracy. The modality $K_X\varphi$ connects back to the standard topological semantics for modal logic, more precisely to the so-called \emph{interior semantics}: note that $s$ satisfies $K_X$ iff it is in the interior of the set of states satisfying $\varphi$ (where interior is taken wrt the topology induced by $X$ on the state space). As such, this notion fits well with the recent epistemic interpretations of the interior operator in terms of ``knowability".

\vspace{-2ex}

\paragraph{Further notions: full independence versus topological independence} In mathematical practice, \emph{independence} of variables often makes for greater simplicity than  dependence. The probability of two independent events is just their product, no  analysis of correlations is required. Independence, too, has a natural interpretation in the above models. Two sets of variables $X$ and $Y$ are fully \emph{independent} if, given the states that occur in our model, knowing the exact values of the variables in $X$ tells us nothing new about the values of $Y$ (i.e. it still leaves them free to take on any combination of values that the model allows them anywhere). Independence is  zero-correlation of $Y$-values with $X$-values, functional dependence is total correlation.\footnote{There are also intermediate forms of correlation, as in the above epistemic models, but we will not study these as such in this paper. However, various forms of value correlation can be expressed in the logical languages for dependence models to be introduced later.} But for empirical variables, this notion of independence is epistemically irrelevant (since typically we cannot know the exact values of $X$). The richer topological setting offers a more subtle concept: two sets of variable $X$ and $Y$ are \emph{topologically independent} if no matter how accurate measurements of $X$ we are given, we can never learn anything new about $Y$. Interestingly enough, topological independence is compatible with global functional dependence! So we can have a total correlation $D_XY$ (which in principle would give us complete knowledge of the exact values of $Y$ if we only were given the exact values of $X$), while at the same time having topological independence (so that no approximate observation of $X$ can tell us anything about $Y$). In our paper, we characterize such situations mathematically in terms of the dependence function being ``everywhere surjective" \cite{Berai18}, thus significantly expanding the current interfaces of Epistemic Topology with mathematics. More philosophically, the extended framework obtained in this way allows us to reason about in-principle learnability and unlearnability on a par.

\vspace{-2ex}
\paragraph{Special topologies} While our discussion may have suggested that we abandon standard (relational) epistemic logic, the topological setting contains in fact standard epistemic models as the special case of \emph{discrete topologies} on the range of values (where all subsets of the space are open). In this case, the topologies induced on the state space are partitional, and we can thus interpret the corresponding propositional questions as ``agents" (modeled using the accessibility relations associated to these partitions). $D_XY$ and $K_XY$ collapse to the same notion in this case (capturing `epistemic superiority' of group $X$ over $Y$ \cite{BalSme20}), and similarly $D_X \varphi$ and $K_X\varphi$ collapse to the standard $\mathsf{S5}$ notion of \emph{distributed knowledge}. Another special class of topologies,  to be used extensively in our later technical proofs, are \emph{Alexandrov spaces} where each point has a smallest oopen neighborhood. These spaces match the usual models for  modal  $\mathsf{S4}$ and support a reduction to relational models where the accessibility relation is the standard topological \emph{specialization order} $s \leq t$. In what follows, such relational models will  allow us to use standard modal techniques.


\paragraph{`Knowing how': from plain continuity to uniform continuity} Even a known epistemic dependence $K(X;y)$ leaves something to be desired. The inverse map associated with a continuous function $F$, from open neighborhoods of values $F(s)$ to open neighborhoods of $s$, is \emph{state dependent}. So, the agent knows \emph{that} she can approximate $y$ as accurately as desired given an accurate enough approximation of $X$, but she may not know in advance \emph{how} accurate needs the measurement of $X$ be (in order to obtain the desired accuracy for $y$).
To have in advance such an explicit `epistemic know-how' (as to which approximation works) she needs the dependence function $F$ to be \emph{uniformly continuous}. Formalizing this stronger notion of \emph{conditional knowability how} (to find $y$ from $X$) requires an even richer setting. We chose to use for this the well-known framework of \emph{metric spaces} \cite{Wil70}, mostly for reasons of accessibility to a general philosophical audience; but uniform continuity can also be defined in a more qualitative way in the much more general framework of \emph{uniform spaces} \cite{Isbell}, and the even more general setting of quasi-uniform spaces, and our notions and results do carry through to these qualitative generalizations.


\bigskip

The \emph{content of this paper} can be viewed  in two ways. Technically, it presents a natural next step for current Logics of Dependence, moving from set-based models to topological ones, and leading to several new completeness, representation, and decidability theorems proved by extending modal techniques to this setting. More philosophically, however,  what we offer in analyzing and defining notions of information-based dependence is an invitation to Epistemic Topology, an area of formal philosophy where many new concepts become available once we go beyond its already established interfaces with mathematics and computer science. In line with this, we see both our definitions and our theorems as main contributions of this paper.

\vspace{-2ex}

\paragraph{Plan of this paper} Section 2 introduces the basic notions of empirical wh-questions, in two equivalent versions: \emph{empirical variables} as topological maps, and \emph{propositional questions} as topologies. Section 3 introduces and investigates various notions of dependence and the corresponding modalities: exact functional dependence (in both its global and local versions) and the determination modality,  knowledge and knowability of a proposition given a set of variables, and most importantly the key notion of \emph{continuous (epistemic) dependence} (in both its global and local versions). Section 4 contains the formal syntax and semantics of our decidable \emph{logic of continuous dependence} $\mathsf{LCD}$, a sound and complete proof system to reason about epistemic dependence, as well as some concrete examples.
Section 5 discusses the notions of full independence and topological independence, and their interplay with functional dependence in a topological setting.
Section 6 extends the setting to uniformly continuous dependence on metric models, and presents the syntax and semantics of the decidable \emph{logic of uniform dependence} $\mathsf{LUD}$ (which in fact also comprises all the other forms of dependence and modalities encountered in this paper), as well as a sound and complete axiomatization of this logic. Section 7 lists some further directions of immediate epistemic relevance, including richer dynamic logics of approximation, a notion of computable dependence understood as Scott continuity, the use of `point-free topologies', etc. A summary of results and general conclusions is found in Section 8. All the proofs of completeness and decidability are found in the Appendix.

\section{Basic framework: empirical questions}

In this section, we present a topological setting for imprecise observations.  Our first task is the introduction of \emph{empirical questions}, for which only imprecise partial answers may be observable. These will come in \emph{two variants}: `\emph{wh}-questions' (whose answers are non-propositional objects or `values': numbers, names or some other form of data), and \emph{propositional} questions (whose answers are propositions). The first type (wh-questions) will be modeled as `\emph{empirical variables}', i.e. quantities that take values in some appropriate value space. Formally, empirical variables are \emph{maps} $X: S\to (\mathbb{D}, \tau)$ from the state space $S$ to some topological space $(\mathbb{D}, \tau)$ (where the open neighborhoods of the actual value of the variable in a state representing the measurable approximations of that exact value). The second type (propositional questions) will be represented as \emph{topologies} $\tau\subseteq {\mathcal P}(S)$ on a state space $S$, with the (inexact) propositional answers being given by the open neighborhoods of the actual state. This generalizes the traditional `partitional' view of propositional questions (in which the possible exact answers form a partition of the state space) to the case of empirical questions, for which typically only inexact, approximative answers are epistemically available: such answers can be more or less precise, and so they are only \emph{partial} answers, which can thus overlap, or even refine each other.

As we shall see, the two settings (wh-questions as empirical variables, and propositional questions as topologies) are in some sense equivalent: we can always `propositionalize' any wh-question, and conversely we can associate a quantitative variable to every propositional question. Moreover, these correspondences are structure-preserving, and so \emph{the two settings obey the same logic} $\mathsf{LCD}$. However, the correspondence is \emph{not} one-to-one: many different wh-questions correspond to the same propositional one. This reflects the intuitive fact that \emph{the first setting (of wh-questions as variables) is} \emph{richer}, giving us more \emph{information} (in the form of the `actual values' of the variable). For this reason, we will choose this first setting as the preferred formal semantics for our logic $\mathsf{LCD}$ (as a logic of dependence between empirical variables).

\subsection{Preliminaries on topology}

While in this paper we assume familiarity with basic notions of General Topology, we briefly review some of them here, mainly in order to fix the notation.

\medskip

\par\noindent\textbf{Topologies, neighborhoods and local bases} \, Recall that a topological space $(\mathbb{D},\tau)$ consists of a set $\mathbb{D}$ (the domain, consisting of `points' or `values' $d\in \mathbb{D}$), together with a family $\tau\subseteq \mathcal{P}(\mathbb{D})$ of `open' sets, assumed to include $\emptyset$ and $\mathbb{D}$ itself, and to be closed under finite intersections and arbitrary unions. The complements of open sets are called 'closed'.
For a point $d\in \mathbb{D}$, we put
$$\tau (d)\,\, :=\,\, \{ U\in \tau: d\in U\}$$
for the \emph{family of open neighborhoods} of $d$. A family of opens $\mathcal{B}\subseteq \tau$ is a \emph{basis} (or `base') for the topology $\tau$ if every open set $U\in\tau$ is a union of open sets from $\mathcal{B}$. In this case, we also say that $\tau$ is \emph{the topology generated by the basis} $\mathcal{B}$.

\medskip
\par\noindent\textbf{Interior and closure} \, The well-known topological \emph{interior} and \emph{closure} operators are defined as usual, as maps $Int, CL: \mathcal{P}(\mathbb{D})\to \mathcal{P}(\mathbb{D})$, defined by putting for all sets of values $D\subseteq \mathbb{D}$:
$$Int(D)=\bigcup\{O\subseteq D: O\in \tau\}, \,\, Cl(D)=\bigcap\{F\supseteq D: (\mathbb{D}-F)\in \tau\}$$
In words, $Int(D)$ is the largest open subset of $D$, while $Cl(D)$ is the smallest closed superset of $D$.

\medskip

\par\noindent\textbf{Specialization preorders, topological indistinguishability and separation axioms}
\, Given a space $(\mathbb{D},\tau)$, we denote by $\leq\,\subseteq \mathbb{D}\times \mathbb{D}$ the \emph{specialization preorder} for its topology, defined by putting, for all $d,c\in \mathbb{D}$:
$$d\leq c  \,\,\, \mbox{ iff } \,\,\, \tau(d)\subseteq \tau(c) \,\,\, \mbox{ iff } \,\,\, \forall U\in \tau (d\in U\Rightarrow c\in U).$$
We also denote by $\simeq$ the \emph{topological indistinguishability} relation, defined as:
$$d\simeq c \,\,\, \mbox{ iff } \,\,\, d\leq c\leq d \,\,\, \mbox{ iff } \,\,\, \tau(d)=\tau(c) \,\,\,  \mbox{ iff } \,\,\, \forall U\in \tau (d\in U\Leftrightarrow c\in U).$$
For any point $d\in \mathbb{D}$, we will denote its equivalence class modulo the indistinguishability relation $\simeq$ by
$$[d]_\tau \,  :=\, \{d'\in  \mathbb{D}: d\simeq d'\}.$$
The space $(\mathbb{D},\tau)$ is \emph{$T_0$-separated} if topologically indistinguishable points are the same, that is, $\tau(d)=\tau(d')$ implies $d=d'$ (for all $d,d'\in \mathbb{D}$); equivalently, iff every equivalence class $[d]_\tau$ is just a singleton $\{d\}$; and again equivalently: iff the specialization preorder $\leq$ is a partial order (i.e., also anti-symmetric).

The space $(\mathbb{D},\tau)$ is \emph{$T_1$-separated} if the specialization preorder is the identity, i.e. $\tau(d)\leq\tau(d')$ implies $d=d'$.

\medskip\par\noindent\textbf{Continuity at a point} A function $F:(\mathbb{D},  \tau_{\mathbb{D}})\to (\mathbb{E}, \tau_{\mathbb{E}})$ between topological spaces is \emph{continuous at a point} $x\in \mathbb{D}$ if we have
$$\forall U\in \tau_{\mathbb{D}} (F(x))\, \exists V\in \tau_{\mathbb{E}} (x):\, F(V)\subseteq U,$$
or equivalently: for every open neighborhood $U\in \tau(F(x))$ of $F(x)$, its preimage is an open neighborhood $F^{-1}(U)\in \tau(x)$ of $x$.

\medskip\par\noindent\textbf{Global and local continuity} $F:(\mathbb{D},  \tau_{\mathbb{D}})\to (\mathbb{E}, \tau_{\mathbb{E}})$ is \emph{(globally) continuous} if it is continuous at all points $x\in \mathbb{D}$. $F$ is \emph{locally continuous at (around) a point $x\in \mathbb{D}$} if it is continuous at all points in some open neighborhood $O\in\tau(x)$ of $x$.

\medskip

\par\noindent\textbf{Subspaces and quotients}
Every subset $D\subseteq \mathbb{D}$ of (the domain of) a topological space $(\mathbb{D},\tau)$ determines a \emph{subspace} $(D, \tau_D)$, obtained by endowing $D$ with its \emph{subspace topology}
$$\tau_D \, : =\, \{U\cap D: U\in \tau\}$$
Every equivalence relation $\sim\subseteq \mathbb{D}\times \mathbb{D}$ on (the domain of) a topological space $(\mathbb{D},\tau)$ determines a \emph{quotient space} $( \mathbb{D}/\sim, \tau/\sim)$, having as set of values/points the set
$$\mathbb{D}/\sim \, :=\, [d]_\sim: d\in \mathbb{D}$$
of \emph{equivalence classes} $[d_\sim: =\{d'\in \mathbb{D}: d\sim d; \}$ modulo $\sim$, endowed with the \emph{quotient topology}
$\tau/\sim$, given by putting for any set $\mathcal{A}\subseteq \mathbb{D}/\sim$ of equivalence classes:
$$\mathcal{A} \in \tau/\sim \, \mbox{ iff } \, \bigcup \mathcal{A}=\{d\in  \mathbb{D}:[d]_\sim\in \mathcal{A}\}\in\tau.$$
The quotient topology is the largest (``finest") topology on $\mathbb{D}/\sim$ that makes continuous the \emph{canonical quotient map} $\bullet/\sim$, that maps every point $d\in \mathbb{D}$ to its equivalence class $[d]_\sim$.

\medskip

\par\noindent\textbf{Compactness and local compactness}
A topological space $(\mathbb{D},\tau)$ is \emph{compact} if every open cover of $\mathbb{D}$ has a finite subcover; i.e., for every collection $\mathcal{C}\subseteq \tau$ of open subsets s.t. $\bigcup \mathcal{C}= \mathbb{D}$, there exists a \emph{finite} subcollection $\mathcal{F}\subseteq \mathcal{C}$ s.t.  $\bigcup \mathcal{F}= \mathbb{D}$. A subset $K\subseteq \mathbb{D}$ of a topological space $(\mathbb{D},\tau)$ is \emph{compact} if the subspace $(K,\tau_K)$ determined by it is compact. A space $(\mathbb{D},\tau)$ is \emph{locally compact} if every point has a compact neighborhood; i.e. for every $d\in \mathbb{D}$ there exists some open set $U\in \tau$ and some compact set $K\subseteq \mathbb{D}$ s.t. $d\in U\subseteq K$.

\medskip\par\noindent\textbf{Metric spaces, pseudo-metric spaces and ultra-(pseudo-)metric spaces}
A \emph{pseudo-metric space} $(\mathbb{D},d)$ consists of a set $\mathbb{D}$ of points, together with a `pseudo-metric' (or `pseudo-distance') $d: \mathbb{D}\to [0,\infty)$, satisfying three conditions: $d(x,x)=0$; $d(x,y)=d(y,x)$; $d(x,z)\leq d(x,y)+d(y,z)$. A \emph{metric space} is a pseudo-metric space $(\mathbb{D},d)$ satisfying an additional condition, namely: $d(x,y)=0$ implies $x=y$. In this case $d$ is called a \emph{metric} (or `distance'). Every pseudo-metric space is a topological space, with a basis given by \emph{open balls}
$$\mathcal{B}(x,\varepsilon):=\{y\in \mathbb{D}: d(x,y)<\eps\},\,\,\, \, \mbox{ with $x\in \mathbb{D}$ and $\varepsilon\in (0,\infty)$.}$$
Note that the interior $Int(A)$ of any set $A\subseteq \mathbb{D}$ is the union of all the open balls included in $A$:
$$Int(A)=\bigcup \{ \mathcal{B}(x,\varepsilon)\subseteq A: \varepsilon>0\} =
\{x\in \mathbb{D}: \mathcal{B}(x,\varepsilon) \subseteq A \mbox{ for some } \eps>0 \};$$
while dually the closure $Cl(A)$ is the intersection of all the \emph{closed balls} that include $A$:
$$Cl(A)=\bigcap \{ \overline{\mathcal{B}(x,\varepsilon)}\supseteq A: \varepsilon>0\}, \mbox{ where }  \overline{\mathcal{B}(x,\varepsilon)}:= \{y\in \mathbb{D}: d(x,y)\leq\eps\}.$$
Metric spaces are $T_1$-separated, hence also $T_0$-separated.\footnote{But pseudo-metric spaces are not necessarily so: in fact, a pseudo-metric space is $T_0$-separated only iff it is a metric space.}

A familiar concrete example of metric spaces are the \emph{Euclidean spaces} $\mathbb{R}^n$ of integer dimensions $n>0$, consisting of $n$-tuples $(x_1, \ldots, x_n)$ of real numbers $ x_1, \ldots, x_n\in \mathbb{R}$, endowed with the Euclidean distance $d(x,y):=\sqrt{\sum_i (x_i -y_i)^2}$. Euclidean spaces have the important property that they are locally compact (but not compact).

Another example are the \emph{ultra-(pseudo-)metric spaces}: a (pseudo-)metric  $d: \mathbb{D}\to [0,\infty)$ is said to be an ultra-(pseudo-)metric if it satisfies a stronger version of the above triangle inequality, namely $d(x,z)\leq max \{d(x,y),d(y,z)\}$.\footnote{A concrete example of ultra-metric spaces are the so-called $p$-adic fields $\mathbb{F}_p$, for any prime number $p$: in a certain mathematical sense, these are ultra-metric analogues of the field of real numbers $\mathbb{R}$.}

\medskip

\par\noindent\textbf{Uniform continuity}
A function $F:(\mathbb{D},  d_\mathbb{D})\to (\mathbb{E}, d_\mathbb{E})$ between two pseudo-metric spaces is \emph{uniformly continuous} if we have
$$\forall \varepsilon>0 \, \exists \delta>0\, \forall x\in \mathbb{D}\, \forall y\in \mathbb{D}\, \left(\, d_\mathbb{D}(x,y)< \delta \Rightarrow d_\mathbb{E}(F(x), F(y)) < \varepsilon\,\right).$$
In contrast, global continuity corresponds in a metric space to swapping the quantifiers $\exists \delta$ and $\forall x$ in the above statement to the weaker form

$$\forall  \varepsilon>0 \,  \forall x\in \mathbb{D} \, \exists \delta>0\,  \forall y\in \mathbb{D}\,
\left(\, d_\mathbb{D}(x,y) < \delta \Rightarrow d_\mathbb{E}(F(x), F(y)) < \varepsilon\,\right).$$
A function $F:(\mathbb{D},  d)\to (\mathbb{E}, d)$ is \emph{locally uniformly continuous at (around) a point} $x\in \mathbb{D}$ if its restriction $F|O$ to (the subspace determined by) some open neighborhood $O=\mathcal{B}(x, \delta)$ is uniformly continuous.

These various forms of continuity may collapse to the same notion in especially favorable situations, of which we give here two examples.

\begin{prop} (Heine-Cantor Theorem)
If $F:(\mathbb{D},  d_\mathbb{D})\to (\mathbb{E}, d_\mathbb{E})$ is a map from a compact metric space to another metric space, then $F$ is (globally) continuous iff it is (globally) uniformly continuous.

If $F:(\mathbb{D},  d_\mathbb{D})\to (\mathbb{E}, d_\mathbb{E})$ is a map from a locally compact metric space to another metric space, then $F$ is locally continuous (on a neighborhood) around a point $x\in \mathbb{D}$ iff it is locally uniformly continuous (on a neighborhood) around $x$.
\end{prop}

For our second example, we need a few more definitions.

\medskip

\par\noindent\textbf{Lipschitz, locally Lipschitz and pseudo-locally Lipschitz functions}
A function $F:(\mathbb{D},  d_\mathbb{D})\to (\mathbb{E}, d_\mathbb{E})$ between two pseudo-metric spaces is \emph{Lipschitz} if there exists a constant $K>0$ (called `Lipschitz constant') s.t. for all $x,y\in \mathcal{D}$ we have that
$$d_\mathbb{D}(x,y)\leq K d_\mathbb{E}(F(x), F(y)).$$
$F$ is \emph{locally Lipschitz} if for every $x\in\mathbb{D}$ there exists some open neighborhood $O=\mathcal{B}(x, \delta)\in \tau_\mathbb{D}(x)$ of $x$, such that the restriction $F|O$ to (the subspace determined by) $O$ is Lipschitz.

Finally, $F$ is \emph{pseudo-locally Lipschitz}\footnote{Unlike the other concepts in this Preliminaries' section, the notion of pseudo-local Lipschitz functions seems to be novel.} if for every $x\in\mathbb{D}$ there exists some open neighborhood $O\in \tau_\mathbb{D}(x)$ of $x$, and some open neighborhood $U\in \tau_{\mathbb{E}}(F(x))$ of $F(x)$, such that the restriction $F|(O\cap F^{-1}(U))$ to (the subspace determined by) the intersection of $O$ with the preimage $F^{-1}(U)$ is Lipschitz.

\medskip

\begin{prop}
A Lipschitz function $F:(\mathbb{D},  d_\mathbb{D})\to (\mathbb{E}, d_\mathbb{E})$ between pseudo-metric spaces is uniformly continuous.

For every point $x\in \mathbb{D}$, a locally Lipschitz function $F:(\mathbb{D},  d_\mathbb{D})\to (\mathbb{E}, d_\mathbb{E})$ is locally uniformly continuous around $x$.

A pseudo-locally Lipschitz function $F:(\mathbb{D},  d_\mathbb{D})\to (\mathbb{E}, d_\mathbb{E})$ is continuous at some point $x\in \mathbb{D}$ iff it is locally uniformly continuous around $x$ (and so also iff it is locally continuous around $x$).
\end{prop}

\begin{proof} The proofs of all the above results are standard, with the exception of the last item. To justify it, assume $F$ is pseudo-locally Lipschitz, and let $O\in \tau_{\mathbb{D}}(x)$ and $U\in \tau_{\mathbb{E}}(F(x))$ be s.t. $F$ is Lipschitz on $O\cap F^{-1}(U)$. Assume now that $F$ is continuous at $x$ in $\mathbb{D}$, hence we have $F^{-1}(U)\in\tau_{\mathbb D}$, and thus $O\cap F^{-1}(U)\in \tau _{\mathbb D}(x)$ is an open neighborhood of $x$. Since $F$ is Lipschitz on the open neighorhood $O\cap F^{-1}(U)$, it is also uniformly continuous on it (by the first item of our Proposition), and therefore it is locally uniformly continuous around $x$.\end{proof}

\medskip

In conclusion, \emph{the notions of local continuity and local uniform continuity are equivalent on locally compact spaces}; while \emph{all three notions of localized continuity (-continuity at a point, local continuity around the point, and local unform continuity around the point) coincide for pseudo-locally Lipschitz functions}! As we'll show, these facts have interesting epistemological and logical consequences, with a special relevance to our completeness and expressivity results.

\subsection{Empirical variables as inexact wh-questions}

In this section, we fix a state space $S$. A \emph{proposition} is just a subset $P\subseteq S$ of the state space $S$. To stress this reading, we sometimes write $s\models P$ for $s\in P$, and read it as ``$P$ is true at $s$ ".

\medskip

As traditionally understood, a \emph{wh-question} is one that starts with ``what", ``who", ``which", ``where", ``when" etc. Examples are ``What is the speed of light?", ``Where on Earth I am?" (i.e. what is my position on the globe), ``When did you first arrive to Amsterdam?". The exact answer to a wh-question is not necessarily a proposition, but can be some other type of object, e.g. a number, a pair or triplet of coordinates, a date, a name etc. Here, we will think of each such object as the ``value" assigned to a variable in a given state of the world. Note that, unlike the variables of predicate logic (which are just abstract, interchangeable placeholders with no intrinsic meaning, that can take any arbitrary values independent of each other), ours are ``variables" in the sense used in empirical sciences and databases: they denote specific quantities or qualitative data (e.g. time, position, velocity, money, persons, shapes etc.), each with its own range of possible values, that may be subject to restrictions and inter-variable correlations. The generic form of a wh-question can thus be reduced to ``What is the value of variable $X$?", where ``the value" refers to the current value of $X$ in the actual world.

When calling our variables ``empirical", what we mean is that \emph{their exact value might typically be inaccessible to observation}. Instead, one can usually observe only \emph{approximations} of that value
 (e.g. a time interval, an area on the surface of the Earth, a range of speeds, a set of possible suspects, etc). So, when thinking of an empirical variable as a wh-question, its \emph{feasible answers} (the ones one could obtain via some kind of empirical process e.g. an observation or experiment) \emph{are typically ``imprecise" or partial answers}, in the sense that may not completely match the exact answer.

We now proceed to formalize this key notion within a topological setting.

\smallskip


\begin{defi} (\textbf{Empirical variables as topological maps})
Given a state space $S$, an \emph{empirical variable}  is a surjective map $X:S \to (\mathbb{D}_X, \tau_X)$ from states $s\in S$ to `values' $X(s)\in \mathbb{D}_X$ in a $T_0$-separated topological space $(\mathbb{D}_X, \tau_X)$.
\end{defi}

Intuitively,  the open neighborhoods $U\in \tau_X (X(s))$ represent \emph{observable approximations} of the value of $X$ at state $s$, obtainable e.g. as a result of a \emph{measurement} of the value of $X$. Both the surjectivity condition, saying that $\mathbb{D}_X=X(S)$,
and the requirement of $T_0$-separation are innocuous: we can always enforce them by taking first a
subspace (restricting to $X[S]$ both the codomain $\mathbb{D}_X$ and
its topology), and then a quotient (identifying all
indistinguishable points in $\mathbb{D}_X$). The point of these requirements is that they simplify several
definitions and statements of  results below.

\medskip

\par\noindent\textbf{Example: Euclidean variables} A well-known example is that of \emph{Euclidean variables} $X: S\to \mathbb{R}^n$, for some $n\in \mathbb{N}$: here, $\mathbb{D}_X=\mathbb{R}^n$ is the $n$-dimensional Euclidean space for some $n\in \mathbb{N}$ (consisting of all real-number vectors $x=(x_1, \ldots, x_n)$), and  $\tau_X$ is the standard Euclidean topology, generated by the family of all $n$-dimensional open balls of the form
$$\mathcal{B}(x, r)=\{y\in \mathbb{R}^n : d(x,y)<r\}$$
where $x\in R^n$, $r>0$, and $d(x,y)=\sqrt{\sum_i (x_i-y_i)^2}$ is the standard Euclidean distance. Any open ball $B(x, r)$, with $X(s)\in B(x,r)$, can be considered as a measurement of the value of $X(s)$ with margin of error $r$.

\medskip
\par\noindent\textbf{Example: a particle in space} For a more concrete example, consider a \emph{point-like particle in three-dimensional space}. The \emph{state space} is $S=\mathbb{R}^3$.  Consider the question ``What is the $x$-coordinate of the particle?". As a wh-question, this is an empirical variable $X: \mathbb{R}^3 \to (\mathbb{R}, \tau)$, with $X(a,b,c):=a$ and the natural topology $\tau$ is given by the family of rational open intervals $(a,b)$, where $a,b\in \mathbb{Q}$ with $a<b$.\footnote{The restriction to intervals with rational endpoints corresponds to the fact that actual measurements always produce rational estimates.} Suppose the particle is actually in the point $(1, 0,3)$: so this triplet $s(1,0,3)$ =is the \emph{actual state}. Then the exact answer to the wh-question is the actual exact value $1$ of the variable $X$ at state $s$, while the feasible approximations at $s$ are rational intervals $(a,b)$ with $a<1<b$. The similar wh-questions regarding the other two coordinates $y$ and $z$ are similarly represented by empirical variables $Y,Z: \mathbb{R}^3 \to (\mathbb{R}, \tau)$, with the same topology $\tau$ as before.

\medskip

All topological notions introduced in the Introduction relativize to the range space $(\mathbb{D}_X, \tau_X)$ of some empirical variable $X$, simply by labeling them with the subscript $X$: so e.g.
$\leq_X\,\subseteq \mathbb{D}_X \times \mathbb{D}_X$ is the specialization preorder for the topology $\tau_X$.
Note that we don't need a special symbol for the corresponding indistinguishability relation on $\mathbb{D}_X$: since the space is assumed to be $T_0$, indistinguishability  on $\mathbb{D}_X$ simply boils down to \emph{equality} $=$.

\medskip\par\noindent\textbf{Special case: exact questions} An \emph{exact wh-question} on a state space $S$ is just a variable $X:S \to (\mathbb{D}_X, \tau_X)$ mapping states into values living in a \emph{discrete} space, i.e. s.t. $\tau_X=\mathcal{P}(\mathbb{D}_X)$. Intuitively, this means that at every state $s\in S$ one can observe the exact value $\{X(s)\}$ of the variable $X$. The topology is then irrelevant: the complete answer to the question can be determined, so no partial approximations are really needed.

\subsection{Abstraction: propositional questions}

A special case of empirical questions are the ones whose answers are propositions: we call these ``propositional questions".  Traditionally, a propositional question is taken to be a partition of the state space $S$: the partition cells are the possible (complete) answers to the question, while the unique cell that contains a given state $s$ is the true answer at state $s$.
But, once again, in empirical sciences the complete answer may be impossible to determine by observations; typically, only partial answers are empirically given. Hence, the partition into complete answers must be replaced by a family of propositions, called `feasible' or `observable' answers.

\smallskip

\begin{defi}(\textbf{Empirical questions of the propositional kind})
A \emph{(propositional) empirical question} is just a topology $\tau\subseteq \mathcal{P}(S)$ on the set of states $S$. If the actual state is $s$, then the open neighborhoods $O\in \tau$ with $s\in O$ represent the feasible or `observable' answers. So $\tau(s)$ is the family of all partial answers to the question that are (true and) observable at state $s$.
\end{defi}

The closure properties of $\tau(s)$ make obvious sense in this interpretation. Closure under binary intersections reflects the fact that the conjunction of two (true) feasible answers is itself (true and) feasible: indeed, $P\wedge Q$ is empirically established iff both $P$ and $Q$ are. Closure under arbitrary unions means that any arbitrary disjunction of (true) feasible answers is itself a (true and) feasible answer to the question: indeed, $\bigvee_i P_i$ is empirically established iff at least one of the $P_i$'s is.\footnote{Note that the corresponding closure condition for arbitrary conjunctions does \emph{not} apply to feasible answers: indeed, while for establishing an infinite disjunction is sufficient (and necessary) to establish only one of the disjuncts, establishing an infinite conjunction $\bigwedge_i P_i$ would necessitate to first establish \emph{all} of the conjuncts. This would require waiting an infinite amount of time, hence it is unfeasible in an empirical context.} Finally, the tautological 'answer' (corresponding to the trivially true proposition given by the whole state space $S$) is surely feasible, albeit uninformative.

\medskip

\par\noindent\textbf{Complete answers are not necessarily feasible}. The \emph{complete 'answer'} to the question $\tau$ at state $s$ is the equivalence class $[s]_\tau:=\{w\in S: s\simeq w\}$
of $s$ wrt the topological indistinguishability relation $\simeq$ for $\tau$. The complete answer is the propositional analogue of the `exact' value of an empirical variable. The fact that the complete answer might not be feasible is reflected in the fact that \emph{this equivalence class is not necessarily an open set}. The indistinguishability relation $\simeq$ induces a partition on $S$, whose cells are the equivalence classes modulo $\simeq$, corresponding to all the possible complete answers to the question.

\medskip

\par\noindent\textbf{Example continued: the particle in space} Recall our example of a point-like particle situated in the point $(1, 0,3)\in S=\mathbb{R}^3$. Instead of answering the wh-question ``what is the $x$-coordinate?" by simply specifying the exact value $X=1$ of the corresponding empirical variable, we can interpret the question as a \textit{propositional} one, whose complete answer is $\{(1,y,z): y, z\in \mathbb{R}\}$, which in English corresponds to the proposition ``The $x$-coordinate of the particle is $1$". Similarly, if only an approximation $(a,b)$ of $X$ (with $a,b\in \mathbb{Q}$ and $a<1<b$) is feasible or measurable, then the same information can be packaged in a \emph{partial} answer $\{(x,y,z)\in \mathbb{R}^3: a<x<b\}$ to the corresponding propositional question: in English, this is the proposition ``The $x$-coordinate is between $a$ and $b$".

\medskip

As seen in this example, \emph{every wh-question can be ``propositionalized"}, i.e. converted into a propositional question, by \emph{abstracting away the actual values}. We proceed now to formalize this process of abstraction from values.

\smallskip

\par\noindent\textbf{From variables to propositional questions: the (weak) $X$-topology on $S$}.
To any empirical variable $X:S\to (\mathbb{D}_X, \tau_X)$, we can associate a propositional empirical question $\tau_X^S$, called
the \emph{$X$-topology on $S$}: this is the so-called `weak topology' induced by $X$ on $S$, defined as \emph{the coarsest topology on $S$ that makes $X$ continuous}. More explicitly, $\tau_X^S$ is
given by
$$\tau_X^S \,\, :=\,\, \{X^{-1}(U): U\in \tau_X\}.$$
We can interpret the $X$-topology $\tau_X^S$ as the result of `propositionalizing' the wh-question $X$.
An \emph{$X$-neighbourhood} of a state $s\in S$ is just a neighbourhood of $s$ in the $X$-topology on $S$, i.e. a set $N$ with $s\in X^{-1}(U)\subseteq N$ for some $U\in\tau_X$. As usual, for a state $s\in S$, we denote by $\tau_X^S(s)$ the family of open $X$-neighborhoods of $s$. Similarly, the \emph{$X$-interior} $Int_X(P)$ of any set $P\subseteq S$ is the interior of $P$ in the $X$-topology, and the same for the \emph{$X$-closure} $Cl_X(P)$. Finally, note that the $X$-topology $\tau_X^S$ is \emph{not} necessarily $T_0$-separated (unlike the topology $\tau_X$ on the value range of $X$, which was assumed to be $T_0$): indeed, the value of a variable $X$ can well be the same in two different states.\footnote{However, it would be natural to further assume that any two states that agree on the values of \emph{all} the relevant variables are \emph{the same} state: that would allow us to identify our `states' with tuples of values (or rather assignments of values to all the relevant variables), obtaining a `concrete' representation of the state space, in the style typically used in Physics and other empirical sciences. This additional assumption will be embodied in the so-called ``concrete models", defined in the next section as a special case of our more general topo-models.}

\medskip

\par\noindent\textbf{$X$-relations on states} The relations $\leq_X$ and $=$ on the range space $\mathbb{D}_X$ naturally induce corresponding relations on \emph{states} in $S$, defined by putting for all $s,w\in S$:
$$s\leq_X^S w \,\,\, \mbox{ iff } \,\,\, X(s)\leq_X X(w) \,\,\,  \mbox{ iff } \,\,\,\forall U\in \tau_X (X(s)\in U\Rightarrow X(w)\in U),$$
$$s =_X^S w \,\,\, \mbox{ iff } \,\,\, X(s) = X(w) \,\,\,  \mbox{ iff } \,\,\, \forall U\in \tau_X (X(s)\in U\Leftrightarrow X(w)\in U).$$
The first relation $\leq_X^S$ coincides with the \emph{specialization preorder for the weak ($X$-)topology} $\tau_X^S$ on states; while the second relation $=_X^S$ corresponds to the (indistinguishability for the $X$-topology, which by $T_0$-separation is the same as the) \emph{equality of $X$-values} on the two states.

\medskip

\par\noindent\textbf{Special case: partitions as exact propositional questions} Note that the propositional counterpart of an \emph{exact} wh-question $X:S \to (\mathbb{D}_X, \mathcal{P}(\mathbb{D}_X))$ on a state space $S$ is a \emph{partition} of $S$, whose partition cells are sets of the form $X^{-1}(d)$, with $d\in \mathbb{D}_X$. Conversely, every partition of $S$ can be viewed as a propositionalized exact question. So, \emph{from a purely propositional perspective, exact questions are just partitions of the state space}.

\subsection{Equivalence between variables and (propositional) questions}

As we saw, we can go from wh-questions (modeled as empirical variables) to propositional questions, by abstracting away from the actual values: replacing the variable $X:S\to (\mathbb{D}_X, \tau_X)$ with the (weak) $X$-topology $\tau_X^S\subseteq \mathcal{P}(S)$. We show now that we can also go back (from questions to empirical variables), essentially by taking the complete answers as actual values, with the topology induced by the feasible answers.

\medskip

\par\noindent\textbf{From propositional questions to variables: the quotient topology}.
Given any propositional empirical question (topology) $\tau\subseteq \mathcal{P}(S)$, we can convert it into a wh-question by taking as our associated empirical variable the \emph{canonical quotient map} $\bullet/\simeq$ from $S$ to the \emph{quotient space} $S/\simeq$ modulo the topological indistinguishability relation for $\tau$; more precisely, we take
$$\mathbb{D}_{X}:=S/{\simeq} =\{[s]_\tau:s\in S\}$$ to be the set of all equivalence classes $[s]_\tau:=\{w\in S: s\simeq w\}$ wrt $\simeq$; $X: S\to S/{\simeq}$ is the quotient map $\bullet/\simeq$, given by
$$X(s):=[s]_\tau;$$
and we endow $\mathbb{D}_X$ with the \emph{quotient topology} $\tau_X$, given by
$$\tau_X:= \{U\subseteq \mathbb{D}_X: X^{-1}(U)\in \tau\} =\{X(O): O\in \tau\}.\footnote{
This last equality is due to the fact that our quotient is modulo the indistinguishability relation for $\tau$.}
$$
It is easy to see that these two operations are inverse to each other. First, if we start with an empirical variable  $X: S\to (\mathbb{D}_X, \tau_X)$, take its associated $X$-topology $\tau_X^S$ on $S$, then take the canonical quotient map $\bullet/=_X^S$ (modulo the indistinguishablity relation for $\tau_X^S$, which as we saw is exactly the relation $=_X^S$ of equality of $X$-values), then the result is `isomorphic' to the original variable $X$: more precise, the topological range space $(\mathbb{D}_X, \tau_X)$ (with the quotient topology) is homeomorphic to the original space $(\mathbb{D}_X, \tau_X)$ via the canonical homeomorphism $\alpha$ that maps every value $X(s)$ to the equivalence class $[s]$ modulo $=_X^S$; and moreover, the quotient variable  $\bullet/=_X^S$ coincides with the functional composition $\alpha\circ X$ of the homeomorphism $\alpha$ and the original variable $X$. Vice-versa, if we start with a propositional question/topology $\tau$ on $S$, take its canonical quotient map $\bullet/\simeq$ (modulo its topological indistinguishability relation $\simeq$), and take the associated weak topology $\tau_{\bullet/\simeq}^S$ on $S$,  we obtain exactly the original topology $\tau$.

\medskip

This equivalence between empirical variables on $S$ (wh-questions) and topologies on $S$ (propositional questions) extends to all the notions defined in this paper, and will form \emph{the basis of two equivalent semantics for our logic} (topo-dependence models and standard topo-models).
However, note that this equivalence is \emph{not} a bijection. While the second correspondence (going from propositional questions to variables, via our quotient construction) is injective, the first correspondence (going from variables to propositional questions, via `propositionalization') is not injective: many different variables can give rise to the same weak topology on the state space. This fact captures the intuition that the setting of empirical variables is in some sense `richer' than the one of propositional questions.

\subsection{Joint questions}

In an inquiry, we are typically interested in answering multiple questions simultaneously: e.g. we want to find the spatial $x$, $y$ and $x$ coordinates of a point, as well as its temporal coordinate $t$ (with respect to some system of reference). We can combine such a set of questions into a single \emph{joint question}, e.g. ``what are the space-time coordinates of the point?". A propositional answer to such a joint question will be a \emph{conjunction} of propositional answers to each of the underlying simple questions, so a  joint propositional question can be seen a the coarsest topology that refines all underlying topologies; while a value/object answer to a joint wh-question will be a \emph{tuple of values} (each answering one of the underlying wh-questions), so a joint wh-question can be modeled using the product topology.

We first start with the propositional case.

\medskip

\par\noindent\textbf{Sets of propositional questions as joint questions: the join topology} Given a set $\tau=\{\tau_i: i\in I\}$ of topologies on a state space $S$ (interpreted as empirical questions), we can regard it as a single question, given by the \emph{supremum} or ``join" $\bigvee_i \tau_i$ (in the lattice of topologies on $S$ with inclusion) of all the topologies $\tau_i$. This is the topology generated by the union $\bigcup_i \tau_i$ of the underlying topologies: the coarsest topology on $S$ that includes all of them. The indistinguishability relation $\simeq_\tau$ for this join topology coincides with the intersection $\bigcap_i \simeq_i$ of the indistinguishability relations $\simeq_i$ of each topology $\tau_i$. The complete answer to the join propositional question $\tau$ at state $s$ is thus the intersection of the complete answers at $s$ to all the propositional questions $\tau_i$; while a partial answer to the join question, i.e. an open set $O\in \tau$, is an intersection $\bigcup_i O_i$ of partial answers $O_i\in \tau_i$ to all the underlying questions $\tau_i$.

\vspace{1ex}

\par\noindent\textbf{Some special cases} When $\tau=\emptyset$ is the empty family of topologies, $\bigvee \emptyset:=\{\emptyset, S\}$ is the \emph{trivial (`indiscrete') topology} on $S$. When $\tau=\{\tau_1\}$ is a \emph{singleton}, the joint question corresponding to the set $\tau$ is the same as the question/topology $\tau_1$.

\vspace{1ex}

\smallskip

We now move on to the more interesting case of joint empirical variables.

\medskip

\par\noindent\textbf{Sets of variables as joint wh-questions: the product topology} A \emph{finite set} $X=\{X_i: i\in I\}$ of empirical variables $X_i:S\to (\mathbb{D}_{i}, \tau_{i})$ can itself be regarded as a \emph{single variable} (hence, our use of the same notation $X$ for both variables and sets of variables, a common practice in, e.g., Statistics when dealing with random variables), essentially given by taking the product map $\Pi_{i\in I} X_i: S\to \Pi_{i\in I} (\mathbb{D}_i, \tau_i)$ into the topological product space (suitably restricted in its range to make this map surjective). More precisely, given such a finite set  $X=\{X_i: i\in I\}$ of empirical variables, we can associate to it a single variable, also denoted by $X$, as follows:

\begin{itemize}
\item the set of values
$\mathbb{D}_X:= \{(x_i(s))_{i\in I}: s\in S\}\subseteq
\Pi_{i\in I} \mathbb{D}_{i}$;

\item the
topology $\tau_X$ is the restriction to $\mathbb{D}_X$ of the \emph{product topology} on $\Pi_{i\in I} \mathbb{D}_{i}$, generated by the restrictions to $\mathbb{D}_X$ of all products $\Pi_{i\in I} U_{i}$ of open sets (in their own topologies $\tau_i$);

\item finally, the map $X: S\to \mathbb{D}_X$ is given by $X(s)=(x_i(s))_{i\in I}$.
\end{itemize}

Essentially, we can interpret the products $\Pi_{i\in I} U_{i}$ as the possible results of a \emph{simultaneous measurement} of all these variables. Since there might be correlations between the variables, the actual range of results $\mathbb{D}_X$ as defined above is not necessarily the whole Cartesian product, but a subset of it. Putting it in our earlier interrogative terms, the variable $X$ represents the \emph{joint question} obtained by simultaneously asking all the wh-questions in the set $X$, and $\Pi_{i\in I} U_{i}$ is a \emph{joint approximate answer} to all these questions.

\vspace{1ex}

\par\noindent\textbf{Some special cases} When $X=\emptyset$ is empty,   $\mathbb{D}_{\emptyset}=\{\lambda\}$ with $\lambda=(\,)$ the empty string, $\tau_{\emptyset}=\{\emptyset, \mathbb{D}_{\emptyset}\}$ the \emph{trivial topology} on $\mathbb{D}_{\emptyset}$ (which in this case coincides with the discrete topology!), and the map $X:S\to \mathbb{D}_{\emptyset}$ given by $X(s)=\lambda$ for all $s\in S$.
Note that, since the weak $\emptyset$-topology $\tau_\emptyset^S=\{\emptyset, S\}$ is the trivial topology on $S$, its indistinguishability relation $=_\emptyset^S$ is the \emph{universal relation} on $S$ (relating every two states). When $X=\{x\}$ is a \emph{singleton}, the single variable corresponding to the set $X$ is the same as $x$ itself.

\vspace{1ex}

Given the preceding, from now on we will use the notation $X$ for \emph{both single variables and for finite sets of variables}, identifying the set $X$ with the associated single variable $X$ (and dually identifying a single variable $x$ with the associated set $\{x\}$).

\bigskip

\par\noindent\textbf{Sanity check: equivalence between the two notions of joint questions}
It is easy to see that above-mentioned equivalence between propositional empirical questions and empirical variables extends to the notion of joint question/variable. Indeed,  the weak topology induced on $S$ by the product variable $X=\Pi_{i\in I} X_i: S\to \Pi_{i\in I} (\mathbb{D}_i, \tau_i)$ is exactly the join supremum) of all the topologies $\tau_i^S$ (with $i\in I$).

\section{Information-carrying dependence as topological continuity}\label{EpistemicDep}

In this section, we present a conceptual analysis of notions of epistemic dependence that make sense in a topological setting for imprecise observations. The outcome is a series of definitions, for which we prove some simple clarificatory characterizations. Deeper technical results on the resulting topological dependence logics are deferred to the next section.

We start by reviewing the older notion of \emph{exact (non-epistemic) dependence} $D_XY$ as it occurs in the Logic of Functional Dependence $\mathsf{LFD}$ in \cite{BaBe21}, then we move to exploring various types of known/knowable dependence.

For this purpose, we fix once again a set of epistemic states or `worlds' $S$: intuitively,  these are all the states that are consistent with some (implicit) agent's background information.

\subsection{Background: exact dependence in $\mathsf{LFD}$}\label{ExactDep}

The form of dependence studied in Dependence Logic ($DL$) \cite{Vaa07}, denoted in $DL$ literature by $=(X;Y)$ and here by $D(X;Y)$, is both \emph{exact} and \emph{global}: it requires the existence of a global functional dependence, that maps all exact values of the variable (or set of variables) $X$ into the exact values of the (set of) variable(s) $Y$. The Logic of Functional Dependence ($\mathsf{LFD}$) introduced in \cite{BaBe21} goes beyond this, by introducing a \emph{local}, but still \emph{exact}, version of dependence $D_XY$: the current exact value of $X$ (in the actual world) uniquely determines the current exact value of $Y$ (in the same actual world).\footnote{This local version of dependence introduced in $\mathsf{LFD}$ is more fundamental, in a sense, than the global one of $DL$: it turns out that local $D_XY$ (combined with local determination modality $D_X\varphi$) can define global $D(X;Y)$, but not the other way around.} In addition, in $\mathsf{LFD}$ one considers propositional dependence operators $D_XP$ (``exact determination"): the (current) exact value of $X$ determines $P$ to be true, These modalities that can be seen as a special case of local dependence $D_XY$, in which the determined variable $Y$ is propositional (denoting a Boolean function $Y_P: S\to \{0,1\}$ from states to truth values, seen as the characteristic function of a set of states $P\subseteq S$) and its uniquely determined value is $1$ (`true').

In this section we review the basics of $\mathsf{LFD}$. Note that these definitions are set-theoretical in nature: they do not use the topological structure of our space of values, but only the (set of) values themselves. This expresses the fact that this form of dependence is exact and `informational' (i.e. the exact value of $X$ carries full information about the exact value of $Y$), rather than approximate and `epistemic' (having to do with an observer's ability to infer as good estimates as needed for the value of $Y$ from sufficiently precise approximations of $X$'s value). The epistemic versions will come in subsequent sections.

We start with the special case of propositional dependence operators:

\medskip

\par\noindent\textbf{Exact determination of a proposition by a variable} Given a state space $S$, a (finite set of) empirical variable(s) $X:S\to  (\mathbb{D}_X, \tau_X)$ and a proposition $P\subseteq S$, we say that $P$ is \emph{determined by $X$ at state $s$}, and write $s\models D_XP$, if the (exact) value of $X$ (at $s$) carries the information that $P$ is true (at $s$):
$$s\models D_X P \,\, \, \mbox{ iff } \, \,\, X^{-1}(X(s))\subseteq P  \, \,\, \mbox{ iff } \, \,\, \forall w\in S ( s=_X^S w \Rightarrow s\models P),$$
where $=_X$ is the equality of $X$-values (defined as above: $s=_X w$ iff $X(s)=X(w)$). Note that $D_X$ is just a standard relational modality, having $=_X$ as its accessibility relation. In particular, since $s =_\emptyset t$ holds for all $s, t\in S$, $D_\emptyset P$ is the \emph{universal modality} (quantifying over all states):
$$s\models D_\emptyset P \,\, \, \mbox{ iff } \, \,\, w\models P \mbox{ for all $w\in S$}.$$
This motivates an abbreviation $A(P)$, saying that $P$ is true in \emph{all} states:
$$A (P)\,\, :=\,\, D_\emptyset P.$$


\par\noindent\textbf{Abstraction: exact determination by a propositional question} By abstracting away from the specific values of the variable, we can `lift' this notion to the level of propositional questions.  Given a state space $S$, a proposition $P\subseteq S$ and a (empirical) propositional question (i.e. topology) $\tau$ on $S$, we say  $P$ is \emph{determined by (the answer of) $\tau$ at state $s$}, and write  $s\models D_\tau P$, if the (complete) answer to $\tau$ (at $s$) carries the information that $P$ is true (at $s$):
$$s\models D_\tau P \,\, \mbox{ iff } \,\, [s]_\tau\subseteq P,$$
where recall that $[s]_\tau=\{w\in S: s\simeq w\}$is the equivalence class of $s$ wrt the indistinguishability relation $\simeq$ for the topology $\tau$. One can easily recognize $D_\tau$ as a standard relational modality having topological indistinguishability as its underlying accessibility relation:
$$s\models D_\tau P \,\, \mbox{ iff } \,\, [s]_\tau\subseteq P \, \, \mbox{ iff } \,\, \forall w\in S ( s\simeq w \Rightarrow s\models P) .$$
Note again that the topology itself plays no role, but only the associated equivalence relation $\simeq$, or equivalently the corresponding partition of the state space into equivalence classes: in effect, this is really about the information carried by the underlying ``exact" (partitional) question.

\medskip

\par\noindent\textbf{Exact dependence between variables}
Given two (finite sets of) empirical variables $X:S\to  (\mathbb{D}_X, \tau_X)$ and $Y:S\to  (\mathbb{D}_Y, \tau_Y)$ over the same state space $S$, we say that \emph{$X$ exactly determines $Y$ at state $s$}, and write $s\models D_XY$, if the value of $X$ at $s$ uniquely determines the value of $Y$:
$$s\models D_X Y \,\, \, \mbox{ iff } \,\,\,  X^{-1}(X(s)) \subseteq  Y^{-1}(Y(s)) \, \,\, \mbox{ iff } \, \,\, \forall w\in S ( s=_X^S w \Rightarrow s=_Y^S w).$$
Since $s =_\emptyset t$ holds for all $s, t\in S$, the statement $D_\emptyset Y$ holds iff $Y$ is a \emph{constant}: its value is the same in all states. We can thus introduce an abbreviation $C(Y)$ (saying that ``$Y$ is a constant"):
$$C(Y)\,\, :=\,\, D_\emptyset Y$$

\smallskip

\par\noindent\textbf{Abstraction: exact dependence between propositional questions} By abstracting again from the specific values, this notion is lifted to the level of propositional questions: given topologies $\tau$ and $\tau'$ on the same state space $S$, we say that \emph{$\tau$ exactly determines $\tau'$ at state $s$}, and write $s\models D_\tau \tau'$, if the exact answer to $\tau$ at $s$ uniquely determines the exact answer to $\tau'$:
$$s\models D_\tau \tau' \,\,\, \mbox{ iff } \,\, \, [s]_\tau \subseteq  [s]_{\tau'},$$
where $[s]_\tau$ and $[s]_{\tau'}$ are equivalence classes wrt indistinguishability relations for $\tau$ and $\tau'$.

Once again, the topologies themselves play no role, but only the corresponding partitions of the state space into equivalence classes. In effect, this notion is purely about the dependence between partitional questions: the (true) answer to the first question gives us full information about the (true) answer to the second question.

\medskip

\par\noindent\textbf{The logic of functional dependence} The operators $D_X P$ and $D_X Y$ are at the heart of the simple Logic of Functional Dependence ($\mathsf{LFD}$)\footnote{This logic was introduced, in a purely set-theoretical setting, devoid of topological features, in \cite{BaBe21}.}, with the following syntax:

\vspace{-3mm}

$$
\begin{array}{c cc cc ccc cc cc cc cc}
\varphi :: =  & P \vec{x}   &|& \neg \varphi &|& \varphi\wedge\varphi    &|& D_X\varphi    &|&
D_XY
\end{array}
$$
where the letters $X$ range over finite sets of variables (coming from a fixed set $Var$), the letters $P$ are predicate symbols (coming from a finite set of such symbols) each having a specific arity $n$, and $\vec{x}=(x_1 \ldots x_n)$ is any $n$-tuple of variables in $Var$.

\medskip

Given the above discussion, it should be clear that this logic has \emph{two equivalent interpretations}: one in terms of \emph{empirical variables} (regarded purely \emph{set-theoretically}, as simple assignments of values to states, with \emph{no} topological structure), and another one in terms of \emph{partitions} (or equivalence relations, hence the name ``relational semantics").

In particular, variable-based
\emph{dependence models} come with: a state space $S$; a way to assign to each syntactic variable $x$ some map $\mathbf{x}:S\to  \mathbb{D}_x$ that associates to each state $s\in S$ some value $\mathbf{x}(s)$ (in some range of values $\mathbb{D}_x$); and a way to interpret each $n$-ary predicate symbol $P$ as an $n$-ary relation $I(P)$ between values. The value assignment can be naturally extended to (finite) sets of variables $X\subseteq Var$, by putting $\mathbf{X}(s):= (\mathbf{x}(s))_{x\in X}$ for all $s\in S$, and then equality of values $=_X$ can be defined as above, by putting: $s=_X w$ iff $\mathbf{x}(s)=\mathbf{X}(w)$. The variable-based semantics of the atoms $P\vec{x}$ is given as usual by putting
$$s\models P x_1\ldots x_n \,\, \mbox{ iff } \,\, (\mathbf{x}_1(s), \ldots, \mathbf{x}(s))\in I(P).$$
The equivalent \emph{relational semantics} only assumes given a family of equivalence relations (or partitions) $=_x$ on $S$, one for each basic variable $x$, and simply treats the atoms $Px_1\ldots x_n$ using a valuation (subject to additional constraints, see the Appendix for details). In both versions, the semantics for Boolean connectives uses the classical Tarskian clauses, while the semantics of exact dependence and determination operators is given by putting:
$$s\models D_X\varphi  \,\, \mbox{ iff } \,\, \forall w\in S (s=_X w \Rightarrow w\models \varphi)$$
$$s\models D_XY \,\, \mbox{ iff } \,\, \forall w\in S (s=_X w \Rightarrow s_Y w)$$
The operators $A\varphi$ (for universal modality) and $C(X)$ (for ``$X$ is a constant) can be defined as abbreviations in this logic, as mentioned above: $A\varphi := D_\emptyset\varphi$, and $C(X):= D_\emptyset X$.

\medskip

\par\noindent\textbf{Global dependence} As already mentioned, the standard notion of (functional) dependence $=(X;Y)$ in Dependence Logic $DL$ is \emph{global}: using our notation $D(X;Y)$, this is usually defined by putting
$$S\models D(X;Y) \,\, \, \mbox{ iff } \,\,\, \exists F: \mathbb{D}_X \to \mathbb{D}_Y \mbox{ s.t. } F\circ \mathbf{X}=\mathbf{Y} \mbox{ holds on } \mathbb{D}_X,$$
which can equivalently be stated in terms of equality of values:
$$S\models D(X;Y) \,\, \, \mbox{ iff } \,\,\, \forall s,w\in S (s=_X w \Rightarrow s=_Y w).$$
For propositional questions, given by topologies $\tau, \tau'$ on $S$, the corresponding notion is simply the global version of local dependence $D_\tau \tau'$:
$$S\models D(\tau; \tau')  \,\, \mbox{ iff } \,\, [s]_\tau \subseteq [s]_{\tau'} \mbox{ holds for all $s\in S$}$$
(where recall that $[s]_\tau$ and $[s]_{\tau'}$ are the equivalence classes of $s$ wrt the indistinguishability relations for the topologies $\tau$ and $\tau'$). Once again, note that this is a purely set-theoretic notion, in which the topologies do not play any role in themselves, but only vis the induced equivalence relations: global functional dependence is also in effect a relation between partitional questions.
Also, it is again easy to see that these two notions of global dependence fit together, via the correspondence given by the weak topology; i.e., we have that
$$S\models D(X;Y) \,\, \mbox{ iff } \,\, X^{-1}(X(s)) \subseteq Y^{-1}(Y(s)) \mbox{ holds for all $s\in S$}.$$
In $\mathsf{LFD}$, global dependence is not a primitive notion, but can be defined via the abbreviation:
$$D(X; Y) \,\, :=\,\, A D_XY.$$
Indeed, it is easy to see that the semantics of the formula $AD_XY$ matches the above semantic clause for $D(X;Y)$.

\medskip

\par\noindent\textbf{A sound and complete proof system} Here is a version of the axiomatic proof system  $\mathbf{LFD}$ for the Logic of Functional Dependence, whose completeness is the main technical result in \cite{BaBe21}:

\begin{table}[h!]
\begin{center}
{\small
\begin{tabularx}{\textwidth}{>{\hsize=0.8\hsize}X>{\hsize=1\hsize}X>{\hsize=0.8\hsize}X}

\toprule

\textbf{(I)} & \textbf{Axioms and rules of Propositional Logic}
 \vspace{1mm} \ \\
 \textbf{(II)} &\textbf{$S5$ Axioms for Determination}: \vspace{1mm}  \ \\
($D$-Necessitation)& From $\varphi$, infer $D_X\varphi$   \vspace{1mm}\ \\
($D$-Distribution) & $D_X (\varphi \to \psi) \to (D_X\varphi \to D_X\psi)$  \vspace{1mm}\ \\
(Factivity: axiom $T$) & $D_X\varphi \to \varphi$   \ \\
(Axiom $4$) & $D_X\varphi \to D_X D_X\varphi$ \ \\
(Axiom $5$) & $\neg D_X\varphi \to D_X \neg D_X\varphi$ \ \\
\vspace{1mm} \ \\
\textbf{(III)} & \textbf{Axioms for exact dependence}: \vspace{1mm} \ \\
(Inclusion) & $D_X Y$, provided that $Y\subseteq X$ \vspace{1mm}\ \\
(Additivity) & $\left(D_X Y\wedge D_X Z\right) \to D_X (Y\cup Z)$  \vspace{1mm}\ \\
(Transitivity) & $\left(D_X Y\wedge D_Y Z\right)\to D_X Z$ \vspace{1mm}\ \\
(Determined Dependence) & $D_XY \to D_X D_X Y$ \vspace{1mm}\ \\
(Transfer) &  $D_X Y \to \left(D_Y\varphi \to D_X\varphi \right)$  \vspace{1mm} \vspace{1mm}\ \\
(Determined Atoms) & $Px_1\ldots x_n \to D_{\{x_1, \ldots, x_n\}}Px_1\ldots x_n $  \vspace{1mm} \ \\
\bottomrule
\end{tabularx}
}
\end{center}
\vspace{-0.6cm}
\caption{The proof system $\mathbf{LFD}$.}\label{tb1}
\end{table}

\vspace{1.5ex}

\subsection{Propositional knowledge and conditional knowability of a proposition}

Next, we take a look  at some epistemic attitudes in an empirical setting. We start with propositional `knowledge', then proceed to analyze when  a proposition $P$ is `knowable' given (a good enough approximation of the value of) some empirical variable $X$.

\medskip

\par\noindent\textbf{Propositional knowledge as universal quantification over epistemic possibilities} We say that a proposition $P\subseteq S$ is \emph{known}, and write $K(P)$, if $P$ holds in all epistemically possible states, i.e., if in fact we have that $P = S$. So, in our setting, the knowledge operator $K$ is simply the above-defined universal modality, quantifying over all epistemic states:
$$K(P)=A(P).$$

\medskip

\par\noindent\textbf{Conditional knowability, given (approximate values of) a variable}
Given a proposition $P\subseteq S$ and an empirical variable $X:S\to  (\mathbb{D}_X, \tau_X)$, we say that \emph{$P$ is knowable given $X$ at state $s$} (or knowable conditional on $X$), and write $s\models K_X P$, if the truth of $P$ at $s$ is determined by \emph{some} (good enough) approximation of the value $X(s)$, i.e. if we have
$$ \exists O\in \tau_X (X(s)):\, X^{-1}(O)\subseteq P.$$
Intuitively, this means that
one \emph{can come to know that $P$ is true}, \emph{after observing a sufficiently accurate estimate of $X(s)$}.

Once again, there is a natural analogue of this notion for propositional questions, obtained by abstracting away from the actual values:

\medskip

\par\noindent\textbf{Abstraction: conditional knowability as interior operator} Given a proposition $P\subseteq S$ and a propositional question (topology) $\tau\subseteq {\mathcal P}(S)$, we say that \emph{$P$ is knowable given question (topology) $\tau$ at state $s$}, and write $s\models K_\tau P$, if one can come to know that $P$ is true at $s$ after learning some (true feasible) answer to question $\tau$; i.e. if \emph{there exists some feasible answer $U\in \tau$ such that $s\in U\subseteq P$}. It is obvious that we have:
$$s\models K_\tau P  \,\, \mbox{ iff } \,\, s\in Int_\tau (P).$$

\par\noindent\textbf{Equivalence between the two notions} One can easily see that conditional knowability given an empirical variable $X$ is equivalent to conditional knowability given its associated propositional question $\tau_X^S$:
\begin{prop}
The following are equivalent, for any variable $X$ and proposition $P\subseteq S$:
\begin{enumerate}
\item $s\models K_X P$
\item $s\in Int_{X} (P)$.
\end{enumerate}
\end{prop}

Our notion of conditional knowability (of a proposition given a variable or question) can thus be seen as a variant of the well-known interior semantics for modal logic.

\medskip

\par\noindent\textbf{Special case: Knowledge as unconditional knowability} When we take the empty set of variables $X=\emptyset$, we obtain the above notion of (unconditional) \emph{knowledge} $K(P)$:
$$K(P):= K_{\emptyset} P$$

\subsection{Knowing the value of a variable}

In a next step, we look at the various ways in which one can be said to \emph{know the value} of an empirical variable.

\medskip

\par\noindent\textbf{Exact knowledge} For a start, we  say that a variable $X$ is \emph{exactly known}, and write $K(X)$, if there exists a value $d\in \mathbb{D}_X$ such that the proposition $K(X=d)$ is known, i.e. if $X(t)=d$ holds for \emph{all} states $t\in S$; equivalently, iff  $\mathbb{D}_X$ is a singleton $\mathbb{D}_X=\{d\}$. So, $K(X)$ holds iff $X$ is a \emph{constant} map:
$$K(X)=C(X).$$

\medskip

\par\noindent\textbf{Approximate observations} The act of observing a given variable $X$ at state $s$ within some given approximation $U\in \tau_X(X(s))$ can be modeled as \emph{restriction of the state space}: we move from the original space $S$ to the subspace $X^{-1}(U)=\{s\in S: X(s)\in U\}$. All the variables $X:S\to (\mathbb{D}_X, \tau_X)$ are thus also restricted to $U$, so we obtain $X: S'\to (\mathbb{D}'_X, \tau'_X)$ where
$S':= X^{-1}(U)$, $\mathbb{D}'_X:= U$, and $\tau'_X:=\tau_X|U =\{O\cap U: O\in \tau_X\}$ is the subspace topology on $U$ .

\medskip

\par\noindent\textbf{Approximate knowledge} For empirical variables, we also have a more inexact form of knowledge, namely \emph{approximate knowledge}. Given any open set $U\in \tau_X$, we say that the \emph{value of $X$ is known with approximation $U$} if we have $K(X\in U)$, that is, $X(t)\in U$ for \emph{all} states $t\in S$.

\medskip

\par\noindent\textbf{Arbitrarily accurate knowledge} We say that the \emph{value of $X$ is known with arbitrary accuracy at state $s$}, and write $s\models k(X)$ if: for every open neighborhood $U\in \tau_X (X(s))$ of $X(s)$, the value of $X$ is known with approximation $U$; i.e., we have that
$$\forall U\in \tau_X (X(s))\forall t\in S ( s\in U \Rightarrow t\in U).$$

\begin{prop}
\begin{enumerate}
\item $s\models k (X)$ iff the singleton $\{X(s)\}$ is \emph{dense} w.r.t. the topology $\tau_X$ (i.e., $Cl(\{X(s)\})=\mathbb{D}_X$).
 \item  \emph{Knowledge of arbitrarily accurate knowledge of $X$ is the same as exact knowledge of $X$}: i.e. we have $K k(X)\Leftrightarrow K(X)$.
\item If $(\mathbb{D}_X,\tau_X)$ is $T_1$-separated, then \emph{arbitrarily accurate knowledge of $X$ is the same as exact knowledge of $X$}: i.e.  we have $k(X)\Leftrightarrow K(X)$.
\end{enumerate}
\end{prop}

\subsection{Known dependence versus knowable dependence}

The local version of (exact) dependence is \emph{not necessarily known}, or even knowable: $D_XY$ does not imply $K D_XY$, nor it implies $K_X D_XY$.
Things change if we consider the notion of \emph{global} functional dependence $D(X;Y)$ introduced above: using the universal modality $A$, we can express the fact that $D(X;Y)$ holds globally via the identity
$$D(X; Y) =A D_XY.$$

\par\noindent\textbf{Global dependence is known dependence}
This follows directly from the interpretation of the universal modality $A$ as ``knowledge" $K$, which together with the above identity gives us
$$D(X;Y)= K D_XY$$
So $D(X;Y)$ captures a situation in which the existence of an exact dependence of $Y$ on $X$ is \emph{known} at $s$.

But now we can also consider a \emph{virtual} version of this:

\medskip

\par\noindent\textbf{Knowable dependence} A knowable dependence of $Y$ on $X$ holds at $s$ iff the existence of an (exact) dependence $D_XY$
\emph{can} be known if one is given some sufficiently accurate estimate of the value $X(s)$. This is captured by the expression
$$K_X D_XY.$$
This holds at a state $s$ whenever $D_XY$ holds on some open neighborhood $O\in \tau_X^S(s)$.

\medskip

However, an important observation is that (for inexact variables) both known dependence $K D_XY$ and knowable dependence $K_X D_XY$ may still be epistemically useless, as far as knowability of $Y$ from  $X$ is concerned. They do \emph{not} automatically guarantee that \emph{any} estimate of the value of $Y$ (no matter how vague) is ever known after observing the value of $X$ with \emph{any} accuracy (no matter how precise)!

What we need instead is a notion of \emph{epistemic} dependence, i.e. one that ensures \emph{knowledge transfer}: we should be able to \emph{infer} the value of $Y$ with any desired accuracy from a sufficiently accurate measurement of the value of $X$. This forms the topic of the next section.

\subsection{Epistemic dependence as continuous correlation}

We now introduce the key notion of our paper: the concept of \emph{epistemic} dependence between (inexact) empirical variables. This is not the same as known or knowable dependence: what we want is rather an \emph{information-carrying dependence}: given two (sets of) empirical variables $X$ and $Y$, when can one be said to \emph{know} how to find the value of $Y$ (with any desired accuracy) given (a sufficiently accurate estimate of) the value of $X$? The answer will be given below by the notion of \emph{known epistemic dependence}, to be written $K(X;Y)$, which requires global continuity of the dependence map.

A related, but weaker notion is that of \emph{knowable epistemic dependence} $K_XY$: this is the case in which a known epistemic dependence between $X$ and $Y$ \emph{can} be acquired after learning some (sufficiently accurate) estimate of $X$. Topologically, this is a more local property, obtained by requiring continuity only over a neighborhood of the current value of $X$. We proceed now to formally define these key concepts.

\medskip

\par\noindent\textbf{Known epistemic dependence} We introduce now a global concept $K(X;Y)$, obtained by adding the \emph{continuity} requirement to the definition of (global) functional dependence:
$$S\models K(X;Y) \,\, \mbox{ iff } \,\,  \mbox{ there is a continuous map $F:\mathbb{D}_X \to \mathbb{D}_Y$ s.t. $F\circ  X=Y$}.$$
To understand why continuity makes this global notion of dependence ``epistemic" (in contrast to global exact dependence $D(X;Y)$), it is useful to provide  a number of other useful equivalent formulations of this notion:

\begin{prop}\label{DepChar}
The following are equivalent for empirical variables $X$ and $Y$:
\begin{enumerate}
\item $S\models K(X; Y)$
\item $\forall s\in S \forall U\in \tau_Y (Y(s)) \, \exists V\in \tau_X (X(s)):\, X^{-1}(V)\subseteq Y^{-1}(U)$
\item $\forall s\in S \forall U \in \tau_Y^S(s) \exists V\in \tau_X^S(s) \, V\subseteq U$
 \item for all states $s\in S$ and sets $P\subseteq S$, if $s\in Int_Y (P)$ then $s\in Int_{X}(P)$
\item $Y: (S, \tau_X^S)\to (\mathbb{D}_Y, \tau_Y)$ is continuous (in the $X$-topology on $S$)
\item the identity function $id:(S, \tau_X^S)\to (S, \tau_Y^S)$ is continuous, if we endow its domain with the $X$-topology and its codomain with the $Y$-topology.
\end{enumerate}
\end{prop}

Essentially, this second and third clause say that, for \emph{every} state $s\in S$, one \emph{can know the value of $Y(s)$ with any desired accuracy}, \emph{after observing a sufficiently accurate estimate of the value of $X(s)$}. Since this fact holds globally (for all states), it is \emph{known} to the observing subject: so the meaning of $K(X;Y)$ is that the subject \emph{knows} how to determine the value of $Y$ with any desired degree of accuracy, if she was given accurate enough approximations of the value of $X$.

\smallskip

By relaxing now this knowledge condition to \emph{knowability}, we get a more ``local" version $K_XY$ of the same concept, obtained by restricting it to an open neighborhood of the actual state:

\medskip

\par\noindent\textbf{Epistemic dependence} We write $s\models K_XY$, and say that there \emph{exists a (knowable) epistemic dependence between $X$ and $Y$ at state $s$}, if a known dependence $K(X;Y)$ can be achieved after observing some accurate enough approximation $O\in \tau_X(s)$ of the value $X(s)$ (of variable $X$ at a state $s$). This happens when $K(X;Y)$ holds \emph{on the subspace} $O$ (even if it does not hold on the whole original space $S$):
$$s\models K_XY \,\, \mbox{ iff } \,\,   \mbox{ there exists $O\in \tau_X^S(s)$ s.t. $O\models K(X;Y)$}.$$
By unfolding this definition, we can easily see that we have: $s\models K_XY$ iff
there exists $O\in \tau_X^S(s)$, and a continuous map $F:X(O) \to \mathbb{D}_Y$, s.t. $F\circ X= Y$ holds on $O$.

The next result gives us other useful characterizations of $K_XY$.

\begin{prop}\label{LocDepChar}
The following are equivalent for empirical variables $X$ and $Y$:
\begin{enumerate}
\item $s\models K_X Y$
\item there exists an open neighborhood $O$ of $X(s)$, and a continuous map $F:O \to \mathbb{D}_Y$, s.t. $F\circ X= Y$ holds on $X^{-1}(O)$
\item $\exists O\in \tau_X^S(s) \forall t\in O \forall U \in \tau_Y^S(t) \exists V\in \tau_X^S(t)) \, V\subseteq U$
 \item there is some open $X$-neighborhood $O$ of $s$, s.t. for all states $t\in O$ and sets $P\subseteq O$, if $t\in Int_Y (P)$ then $t\in Int_{X}(P)$
\item $Y$ is continuous in the $X$-topology on some $X$-open neighborhood of $s$
\item the restriction of the identity function $id:(O, \tau_X^O)\to (S, \tau_X^S)$ to some open $X$-neighborhood $O\in\tau_X^S(s)$ is continuous, if we endow its domain with (the topology $\tau_X^O$ induced on $O$ by) the $X$-topology and its codomain with the $Y$-topology.

\end{enumerate}
\end{prop}

So, unlike the local version of exact dependence $D_XY$, \emph{local epistemic dependence $K_XY$ is inherently knowable}: i.e. we always the validity
$$K_XY \, \Rightarrow K_X K_X Y.$$
As in the case of exact dependence, it turns out that \emph{the local version of epistemic dependence is more fundamental than the global one}: $K_XY$ can define $K(X;Y)$ via the equivalence
$$K(X;Y)  = KK_XY =AK_XY,$$
but not the other way around.

\subsection{A side issue: continuity at a point}\label{point-cont}

There is in fact an even ``more local" form of dependence $k_XY$, that requires only \emph{continuity at a point}. We call this notion ``conditional knowability of $Y$ given $X$ ". As we'll see, this form of dependence is in general epistemically \emph{opaque}. As such, it will not play any role in this paper, but it may still be useful to define it, in order to better understand our key notion of epistemic dependence $K_XY$: the two are closely related, but nevertheless subtly different in important ways.

\medskip

\par\noindent\textbf{Conditional knowability of variables as point-continuous dependence} We write $s\models k_XY$, and say that \emph{$Y$ is conditionally knowable given $X$ at state $s$} if we have
$$\forall U\in \tau_Y (Y(s)) \, \exists V\in \tau_X (X(s)):\, X^{-1}(V)\subseteq Y^{-1}(U),$$
i.e. if one can approximate the value(s) of $Y$ with any desired accuracy if given a sufficiently accurate estimate of the value(s) of $X$.

We can characterize conditional knowability $k_XY$ in a similar way to the way we did for $K_XY$ and $K(X;Y)$:

\begin{prop}\label{LocDepChar}
The following are equivalent for empirical variables $X$ and $Y$:
\begin{enumerate}
\item $s\models k_X Y$
\item $ \forall U \in \tau_Y^S(s) \, \exists V\in \tau_X^S(s) \, V\subseteq U$
 \item for all sets $P\subseteq O$, if $s\in Int_Y (P)$ then $s\in Int_{X}(P)$
\item $Y$ is continuous at point $s$ in the $X$-topology
\item the identity function $id:S\to S$ is continuous at point $s$, if we endow its domain with the $X$-topology and its codomain with the $Y$-topology
\end{enumerate}
\end{prop}

However, unlike $K_XY$ and $K(X;Y)$, conditional knowability $k_XY$ is \emph{itself not necessarily knowable} (based on observing $X$):
$$k_XY \not\Rightarrow K_X k_XY$$
In other words, even if $Y$ is actually knowable given $X$ at $s$, this very fact may forever remain unknown:
no matter how precise
her measurements of $X(s)$, the agent may never be in a position to know this fact. The epistemic opacity of this form of dependence is due to its extreme sensitivity to small
deviations or errors: if the value of $X$ happened to be even slightly different than the exact current value $X(s)$, then $k_XY$ would no longer
hold. As a result, any imprecision, however small, is enough to cast doubt upon conditional knowability.

\medskip

\par\noindent\textbf{Counterexample: Thomae's function} An extreme form of unknowable conditional dependence is given by \emph{Thomae's function}. Let $X$ and $Y$ be two single variables, $S=\mathbb{D}_X=\mathbb{D}_Y=\mathbb{R}$ (the set of real numbers), with $\tau_X=\tau_Y$ being the standard topology on $\mathbb{R}$, and we set $X=id$ to be the \emph{identity} map, and $Y$ be given by Thomae's function:
$$Y(s) := \frac{1}{q} \, \mbox{ if $s=\frac{p}{q}$ with $p\in \mathbb{Z}, q\in \mathbb{N}$ s.t. $gcd(p,q)=1$}, \,\,\,\,\, \mbox{ and } \,\,\,\,\, Y(s)=0 \, \mbox{ if $s\not\in \mathbb{Q}$ is irrational}.$$
This function is continuous at all irrational points, and discontinuous at all rational points.\footnote{Note that $Y(0)=Y(0/1)=1$, since by convention $q=1$ is the only natural number s.t. $gcd(0,q)=1$.} As such, \emph{$Y$ is never continuous on any open interval}, and as a consequence, we have in this case that
$$k_XY \Rightarrow \neg K_X k_XY$$
is valid on the space $S$. This is the extreme opposite of inherent knowability: in this space, \emph{whenever $k_XY$ is true, then it is unknowable}!

\medskip

In fact, conditional knowability $k_XY$ is only knowable (based on observing $X$) when we have a full-fledged epistemic dependence $K_XY$, as witnessed by the following equivalence:
$$K_XY \, \Leftrightarrow \, K_X k_XY$$
In other words,\emph{ the `knowable' version of conditional knowability $k_XY$ is exactly the epistemic dependence $K_XY$}. For this reason, we
consider $k_XY$ to be a somewhat irrelevant notion from an epistemological point of view, and we will not explore it in any depth in this paper.\footnote{Nevertheless, we will later specify sufficient conditions in which $k_XY$ \emph{is} inherently knowable (and thus equivalent to $K_XY$): as we'll see, such a `paradisiacal' epistemic situation is guaranteed in a special kind of metric dependence models (``pseudo-locally Lipschitz models").}

\section{The logic of continuous dependence}\label{LCD}

In this section we develop a formal language for an epistemic Logic of Continuous Dependence $\mathsf{LCD}$, that can express both (exact) functional dependence  and (continuous) epistemic dependence.
The language we will use includes the syntax of the Logic of Functional Dependence $\mathsf{LFD}$ introduced in Section \ref{ExactDep}, but the semantics is topological, and the syntax is enriched with topological \emph{interior modalities} $K_X\varphi$ (for `knowability' of a proposition $\varphi$ given an empirical variable $X$) and \emph{continuous-dependence atoms} $K_XY$ (for `epistemic dependence' between variables $X$ and $Y$).
We provide a sound and complete proof system for reasoning about these forms of dependence, and we illustrate our setting with some examples.

\subsection{Syntax and semantics}

We assume given a (finite or infinite) set of \emph{basic variables} $Var$, and a \emph{relational vocabulary} $(Pred, ar)$ consisting of a set $Pred$ of \emph{predicate symbols}, together with an \emph{arity map} $ar: Pred\to N$, associating to each symbol $P\in Pred$ a natural number $ar(P)\in N$. We denote by $X\subseteq Var$ finite sets of variables in $V$, and denote by $\vec{x}=(x_1\ldots x_n)$ finite tuples of variables.

\begin{defi}[\textbf{Syntax of $\mathsf{LCD}$}] The logic of continuous dependence has a language defined recursively by the following clauses:

\vspace{-3mm}

$$
\begin{array}{c cc cc ccc cc cc cc cc}
\varphi :: =  & P \vec{x}   &|& \neg \varphi &|& \varphi\wedge\varphi    &|& D_X\varphi      &|&  K_X \varphi &|&
D_XY &|& K_X Y
\end{array}
$$
where $n$ is the arity of $P$ and $\vec{x}=(x_1, \ldots, x_n)$ is any $n$-tuple of variables in $Var$.
\end{defi}

\par\noindent\textbf{Abbreviations} \, Knowledge $K\varphi$, knowable dependence $K_X Y$, known dependence $K(X;Y)$ and exact knowledge $K(Y)$ are defined in this language as abbreviations:
$$K\varphi \,\,\, :=\,\,\, K_{\emptyset} \varphi, \,\,\,\,\,\,\,\,\,\,\,\,
K(X; Y)  \,\,\, :=\,\,\, K K_X Y, \,\,\,\,\,\,\,\,\,\,\,\,
K(Y)  \,\,\, =\,\,\, K_{\emptyset} Y$$
It is easy to see that these abbreviations fit the semantic definitions of the corresponding operators given in previous sections, e.g. we have
$$s\models K K_XY \,\,\, \mbox{ iff } \,\,\, S\models K(X;Y).$$

\begin{defi}[\textbf{Topo-dependence models}] A \emph{typed topological model} is a multi-typed structure $M=(\mathbb{D}_x, \tau_x, I)_{x\in Var}$, indexed by variables $x\in Var$ (each thought as having its own distinct type), where: for each $x\in Var$, $\mathbb{D}_x$ is a set, giving the range of values of variable $x$; $\tau_x$ is a topology on $\mathbb{D}_x$; and $I$ is an interpretation function, mapping each predicate symbol $P$ of arity $n$ into an $n$-ary relation on the union $\bigcup_{x\in Var} \mathbb{D}_x$.

A \emph{topological dependence model} (`topo-dependence model' for short) is a structure $\bM=(M, S, \mathbf{x})_{x\in Var}$, consisting of: a typed topological model $M=(\mathbb{D}_x, \tau_x, I)_{x\in Var}$; a set $S$ of abstract \emph{states}; and, for each variable $x\in Var$, a corresponding empirical variable, i.e. a surjective map $\mathbf{x}:S\to \mathbb{D}_x$.

For every set $X\subseteq Var$ of variables,  $\mathbf{X}:S\to \mathbb{D}_X=\{(x(s))_{x\in X}: s\in S\}$ stands for the single joint empirical variable associated to the set $\{\mathbf{x}: x\in X\}$, as defined in the preceding section: $\mathbf{X}(s):=(\mathbf{x}(s))_{x\in X}$ for all $s\in S$.
We also use all the associated notation introduced in the previous section: in particular, for every $X\subseteq Var$, we consider the $\mathbf{X}$-topology on $S$, given by $\tau_{\mathbf{X}}^S:= \{\mathbf{X}^{-1}(U): U\in \tau_{\mathbf{X}}\}$, and the $\mathbf{X}$-interior $Int_{\mathbf{X}}:\mathcal{P}(S)\to \mathcal{P}(S)$ (defined as the interior operator in the $\mathbf{X}$-topology).
\end{defi}

\begin{defi} [\textbf{Semantics of $\mathsf{LCD}$}] \emph{Truth of a formula} $\varphi$ in a topo-dependence model  $\bM=(M,S)$ at a state $s \in S$ (written $\bM, s \models \varphi$, where we drop the index $\bM$ when the model is understood) is defined as in $\mathsf{LFD}$ for the atoms, as well as for all the operators of $\mathsf{LFD}$, and in addition by the following clauses:
$$s\models K_X \varphi \,\,\, \mbox{ iff } \,\,\, s\in Int_{\mathbf{X}} \|\varphi \|$$
$$s\models K_X Y \,\,\, \mbox{ iff } \,\,\,  \exists O\in \tau_{\mathbf{X}}^S(s) \, \exists
F:\mathbf{X}(O)\to \mathbb{D}_Y \mbox{ s.t. } F \mbox{ continuous and } F\circ \mathbf{X}= \mathbf{Y} \mbox{ holds on } O.$$
\end{defi}


The paradigmatic examples in topology are metric spaces, which give us our paradigmatic type of topo-dependence models:

\smallskip

\par\noindent\textbf{Special case: metric models}
A \emph{metric dependence model} (`metric model', for short) is a topo-dependence model $\bM$, in which each of the topologies $\tau_x$ is given by a metric $d_x$ on $\mathbb{D}_x$. The underlying typed metric model $M=(\mathbb{D}_x, d_x, I)_{x\in Var}$ comes with designated metrics $d_x$ instead of topologies, but each of them induces of course a topology $\tau_x$ having the family of open $d_x$-balls as its basis.  For joint empirical variables $\mathbf{X}$ given by finite sets of variables $X\subseteq Var$, the induced (product) topology $\tau_{\mathbf{X}}$ on $\mathbb{D}_X$ is also metric, and is easily seen to be generated by open balls
wrt the (the restriction to $\mathbb{D}_X$ of) the so-called \emph{Chebyshev distance} (also known as the \emph{supremum metric}, or $L_\infty$ metric):\footnote{
This metric is uniformly equivalent with the more standard Euclidean metric on the product space.}
$$d_X (\vec{u},\vec{v}) \, :=\, sup \{d(u_x, v_x): x\in X\}, \mbox{ for all } \vec{v}=(v_x)_{x\in X}, \vec{u}=(u_x)_{x \in X}\in \Pi_{x\in X} \mathbb{D}_{x}.$$
Note that, for (finite) \emph{non-empty} sets of variables $X\neq\emptyset$, this amounts to taking the maximum of all distances
$$d_X (\vec{u},\vec{v}) \, :=\, max \{d(u_x, v_x): x\in X\},$$
while for $X=\emptyset$, recall that $\mathbb{D}_\emptyset=\{\lambda\}$ where $\lambda=()$ is the empty string, so that we have that
$$d_\emptyset (\lambda, \lambda) \, :=\, sup \, \emptyset =0  \, \mbox{ (so that $d_\emptyset$ is still a metric!)}$$
The corresponding (weak) $\mathbf{X}$-topologies $\tau_{\mathbf{X}}^S$ induced on the state set $S$ are not necessarily metric, but they are \emph{pseudo-metric}, being generated by open balls of the form
$$\mathcal{B}_X(s, \varepsilon)\,  :=\, \{w\in S: d_{\mathbf{X}}^S (s,w)<\varepsilon\},$$
with $s\in S$, $\varepsilon>0$ and where the pseudo-metric $d_{\mathbf{X}}^S$ on $S$ is given by setting:
$$d_{\mathbf{X}}^S (s,w)\, :=\,  d_{\mathbf{X}}(\mathbf{X}(s), \mathbf{X}(w)).$$
Note that as a result, $d_{\mathbf{X}}^S$ coincides with the supremum pseudo-metric induced by single-variable pseudo-distances:
$$d_{\mathbf{X}}^S (s,w) \, = \,  sup \{ d_x^S(s,w): x\in X\}$$
With these notations, it is easy to see that in a metric model our topological semantics amounts to the following:
$$s\models D_X \varphi \, \, \mbox{ iff } \,\, \forall w\in S \,(d_X^S(s,w)=0\Rightarrow w\models \varphi)$$
$$s\models D_XY \, \, \mbox{ iff } \,\, \forall w\in S \,(d_X^S(s,w)=0\Rightarrow d_Y^S(s,w)=0)$$
$$s\models K_X\varphi  \, \, \mbox{ iff } \,\, \exists \delta>0 \,\forall w\in S\, (d_X^S(s,w)<\delta \Rightarrow w\models \varphi)$$
$$s\models K_XY  \, \mbox{ iff } \, \exists \delta_0>0 \, \forall t\in \mathcal{B}_X(s, \delta_0) \, \forall \varepsilon>0 \, \exists \delta>0 \, \forall w\in \mathcal{B}_X(s,\delta_0)
\, ( d_X^S(w,t)<\delta \, \Rightarrow \, d_Y^S(w,t)< \varepsilon)$$

\medskip

In many applications, it is  useful to think of the \emph{states} employed so far as being `concrete', i.e., represented by means of tuples of values for each of the fundamental variables; in other words, \emph{assignments} of values to variables:

\smallskip

\par\noindent\textbf{Special case: concrete models}
A \emph{concrete topo-dependence model} is a structure $\bM=(M, S)$, consisting of: a typed topological model $M=(\mathbb{D}_x, \tau_x, I)_{x\in Var}$; and a set $S\subseteq \Pi_{x\in Var} \mathbb{D}_x$ of \emph{`concrete'  states}, i.e. type-respecting assignments of values to variables. This structure is subject to the additional requirement that $\mathbb{D}_x=\{s(x): s\in S\}$ for every $x\in Var$.\footnote{Once again, this condition is innocuous (we can always restrict the codomain $\mathbb{D}_x$ to the actual range of values taken by $x$), and is just meant to enforce the surjectivity of the corresponding empirical variables.} A concrete topo-dependence model is indeed a special case of a topo-dependence model: we can  associate to each variable $x\in Var$ an empirical variable $\mathbf{x}: S\to \mathbb{D}_x$, given by
$\mathbf{x}(s)  := s(x)$.

\smallskip

\par\noindent\textbf{Example: Euclidean models} An example of topo-dependence models that are both metric and concrete are \emph{Euclidean models}. For a finite set $Var$, a Euclidean model is simply given by a subset $S\subseteq \mathbb{R}^{Var}$ of the Euclidean space of dimension $|Var|$, consisting of assignments $s:Var\to \mathbb{R}$ of real values $s(x)\in \mathbb{R}$ to variables $x\in Var$. The metric on each copy of $\mathbb{R}$ is the standard Euclidean distance $d(s,w)=\sqrt{\sum_i (s(x)-y(s))^2}$, and so the topology $\tau_x$ is the subspace topology induced on $\mathbb{D}_x=\{s(x): s\in S\}$ by the standard Euclidean topology on $\mathbb{R}$.

\subsection{The proof system $\mathbf{LCD}$}

We are now prepared to present our axiomatic proof system $\mathbf{LCD}$ for the Logic of Continuous Dependence. The axioms and rules are given in  Table \ref{tb2}. The system includes, as group (I) in the Table, the axioms and rules of the system $\mathsf{LFD}$ for the Logic of Functional Dependence, as already listed in Table \ref{tb1}. The additional axioms and rules are divided into three further groups: (II) Axioms for propositional knowability, (III) Axioms for knowable (epistemic) dependence, and (IV) Connecting Axioms (that connect the epistemic/inexact notions with the corresponding exact notions).

Note the analogy of the $\mathbf{LCD}$ axioms in groups (II) and (III) of Table 2 with the similar $\mathbf{LFD}$ axioms in groups (II) and (II) of Table 1 (which are of course also included in  $\mathbf{LCD}$, within group (I) of Table 2). This is to be expected, since epistemic dependence is in a sense the ``feasible" generalization of exact dependence to inexact variables; mathematically speaking, $K_X$ is the topological analogue of $D_X$. But note first that \emph{not all} the $\mathsf{LFD}$ have analogues for epistemic dependence: axiom (5) for $D_X$ has no analogue for $K_X$, and the same applies to the $\mathsf{LFD}$ axiom of Determined Atoms; and second, that even in the case of the analogue axioms which \emph{do} hold for inexact dependence,
their \emph{meaning} is different in the topological setting. We discuss this shift in some detail.

\begin{table}[h!]
\begin{center}
{\small
\begin{tabularx}{\textwidth}{>{\hsize=0.8\hsize}X>{\hsize=1\hsize}X>{\hsize=0.8\hsize}X}

\toprule

\textbf{(I)} & \textbf{Axioms and rules of $\mathbf{LFD}$}
 \vspace{1mm} \ \\

 \textbf{(II)} &\textbf{$S4$ Axioms for Knowability}: \vspace{1mm}  \ \\
($K$-Necessitation)& From $\varphi$, infer $K_X\varphi$   \vspace{1mm}\ \\
($K$-Distribution) & $K_X (\varphi \to \psi) \to (K_X\varphi \to K_X\psi)$  \vspace{1mm}\ \\
(Veracity) & $K_X\varphi \to \varphi$   \ \\
(Positive Introspection) & $K_X\varphi \to K_X K_X\varphi$
\vspace{1mm} \vspace{1mm}\ \\
\textbf{(III)} & \textbf{Axioms for Knowable Dependence}: \vspace{1mm} \ \\
($K$-Inclusion) & $K_X Y$, provided that $Y\subseteq X$ \vspace{1mm}\ \\
($K$-Additivity) & $\left(K_X Y\wedge K_X Z\right) \to K_X (Y\cup Z)$  \vspace{1mm}\ \\
($K$-Transitivity) & $\left(K_X Y\wedge K_Y Z\right)\to K_X Z$ \vspace{1mm}\ \\
(Knowability of Epistemic Dependence) & $K_X Y \to K_X K_X Y$ \vspace{1mm}\ \\
(Knowability Transfer) &  $K_X Y \to \left(K_Y\varphi \to K_X\varphi \right)$
\vspace{1mm} \vspace{1mm}\ \\
\textbf{(IV)} & \textbf{Connecting Axioms}: \vspace{1mm} \ \\
(Knowable Determination) & $K_X\varphi \to D_X\varphi$ \vspace{1mm} \ \\
(Knowable Dependence) & $K_XY \to D_XY$  \vspace{1mm} \ \\
(Knowledge of Necessity) & $A\varphi \to K\varphi$  \vspace{1mm} \ \\
(Knowledge of Constants) & $C(Y) \to K(Y)$ \vspace{1mm} \ \\
\vspace{1mm} \vspace{1mm}\ \\
\bottomrule
\end{tabularx}
}
\end{center}
\vspace{-0.6cm}
\caption{The proof system $\mathbf{LCD}$, with notations
$A\varphi:=D_\emptyset\varphi$, $C(Y):=D_\emptyset Y$, $K\varphi:=K_\emptyset\varphi$, $K(Y):=K_\emptyset Y$.}\label{tb2}
\end{table}

\vspace{1.5ex}

\par\noindent\textbf{What the $\mathsf{LCD}$ axioms mean}  Mathematically, the $S4$ axioms for  $K_X$ in Group {\bf II} capture the main properties of the interior operator, as given by the Frechet axioms of topology. In addition, Veracity asserts that knowable facts are true.\footnote{This axiom is listed here only for expository reasons, since is in fact derivable from the axiom of Knowable Dependence, together with the $\mathbf{LFD}$ axiom of Factivity. Similarly, one can easily see that the old $\mathbf{LFD}$ rule of $D$-Necessitation is now derivable in $\mathsf{LCD}$ from $K$-Necessitation together with Knowable Dependence.}

The Axioms of Knowable Dependence in Group {\bf III} are similar to the corresponding $\mathsf{LFD}$ axioms for exact dependence, and they are familiar from Database Theory, where they are known as the Armstrong axioms of dependence. Inclusion expresses ``epistemic superiority": (the approximate answers to) a larger set of questions carry all the information (and more) that is carried by (the approximate answers to) a subset of those questions. Additivity says that, if (the inexact answers to) two sets of questions are knowable then (the inexact answers to) all the questions in their union are knowable. Transitivity captures a version of Dretske's Xerox Principle \cite{Dretske}: if $X$ carries all the information about $Y$, and $Y$ carries all the information about $Z$, them $X$ carries all the information about $Z$.

Mathematically, these Axioms of Knowable Dependence capture the main properties of the category of topological spaces, with product space as a categorical product:
the Inclusion axiom holds because of the \emph{continuity of the projection maps} $\pi^X_Y: \mathbb{D}_X\to \mathbb{D}_Y$ for $Y\subseteq X$ (where recall that $\mathbb{D}_X:=X(S)$ is a subspace of $\Pi_{x\in X} \mathbb{D}_x$ with the product topology); the Additivity axiom captures the \emph{universality property of the product}, saying that a pair of continuous functions $F:U_X \to \mathbb{D}_Y$ and $G:V_X \to \mathbb{D}_Z$ on two opens $U_X, V_X\subseteq \mathbb{D}_X$ gives rise to a continuous function $(F,G): U_X\cap V_X \to \mathbb{D}_Y\times \mathbb{D}_Z$ into the product space, given by $(F,G)(x):= (F(x), G(x))$; and Transitivity captures the \emph{closure of continuous functions under composition}.
Next, the axiom of Knowability of Epistemic Dependence  expresses the fact that $K_X Y$ is `topologically local' (-it holds at a point only if it holds on a whole open neighborhood of that point), and that it is stronger, in general,  than simple functional dependence. Epistemically, this says that epistemic dependence $K_X Y$ is indeed a knowable dependence: whenever it holds, one can come to know it (after observing some good enough approximation of $X(s)$).
Going further, the Knowability Transfer axiom captures the continuity of dependence: if $F$ is continuous, then the inverse $F^{-1}(U)$ of any open subset $U$ of its domain is open. Epistemically: if $Y$ is knowable given $X$, then any proposition that is knowable given $Y$ is also knowable given $X$.

Finally, the Connection Axioms in Group {\bf IV} spell out the relationship between epistemic (inexact) dependence and knowability and their exact analogues in  $\mathsf{LFD}$. Knowable Determination says
that knowability implies determination: if a statement is knowable based on some approximation of the value of $X$, then the truth of that statement is determined by the (exact) value of $X$. Knowable Dependence is the analogue claim for empirical variables: if $Y$ is epistemically dependent on $X$, then the exact value of $Y$ could be determined of one was given the exact value of $X$. Mathematically, this just follows from the fact that a continuous dependence \emph{is} a (functional) dependence. The last two Connection Axioms tell us that the converses of these statements hold (only) in the special case when $X=\emptyset$ is the empty set of variables: mathematically, this is because $\tau_\emptyset$ is by definition the \emph{discrete topology} (so that $Int_{\emptyset}$ is the \emph{global quantifier} on $S$, while $K_\emptyset Y$ simply means that $Y$ is a constant). Epistemically, Knowledge of Necessity says that, if a statement is necessary (true in all epistemically possible worlds), then it is known. Similarly, Knowledge of Constants says that, if $X$ a constant (having the same value in all epistemically possible worlds), then its value is known.

\medskip

\begin{thm}\label{Completeness} The proof system  $\mathbf{LCD}$ in Table \ref{tb1} is a sound and complete axiomatization of the logic of continuous dependence on any of the following classes of models: (a) arbitrary topo-dependence models; (b) metric dependence models; (c) concrete topo-dependence models; (d) pseudo-locally Lipschitz metric dependence models. Moreover, the logic $\mathsf{LCD}$ is decidable. \end{thm}

\begin{proof} Proofs for these assertions  are found in the Appendix to this paper. \end{proof}

\begin{rem}[A first-order perspective] Why do we succeed in axiomatizing the  logic of topological dependence? It is well-known that many modal logics, including the basic dependence logic $\mathsf{LFD}$, admit  a faithful effective \emph{translation} into the language of first-order logic. Thus, in principle, these logics are completely axiomatizable, and other properties follow as well. For our topological semantics, an obvious translation would use a \emph{three-sorted} first-order language, with sorts for states, values of variables, and open sets. However, the standard requirement on topological spaces that opens be closed under arbitrary unions is not first-order. But then, in our motivation for topological models, we thought of open sets as outcomes of possible measurements. Usually not all open sets qualify for this. For instance, in the reals, open intervals are just a \emph{base} for the topology: closed under finite intersections, but not under arbitrary unions. If we think of our second sort as an open base for the topology, requirements remain first-order expressible. Moreover, it is easy to add explicit first-order descriptions for the behavior of  variables viewed as names of functions, and as a result, at least in principle, we are still in a  first-order setting that predicts axiomatizability in general for the systems of reasoning investigated in this paper. Of course, to find  explicit complete logics, we still have to do real work. \end{rem}

Note that the point-continuity atom $k_XY$ introduced in Section \ref{point-cont} (capturing the epistemically opaque notion of ``conditional knowability of $Y$ given $X$") is \emph{not} a part of the syntax of $\mathsf{LCD}$. And indeed, we can use part (d) of Theorem \ref{Completeness} to show that:

\begin{prop}\label{expressk}
Conditional knowability $k_XY$ is \emph{not} expressible in the language of $\mathsf{LCD}$. As a consequence, the language $\mathsf{LKk}$ based only on the operators $K_X\varphi$ and conditional knowability atoms $k_X Y$ (i.e. the language obtained by replacing in $\mathsf{LCD}$ the local continuity atoms $K_XY$ by point-continuity atoms $k_XY$) is more expressive than $\mathsf{LCD}$.\footnote{The fact that $\mathsf{LKk}$ is at least as expressive as $\mathsf{LCD}$ follows already from the fact that $K_XY$ is definable in  $\mathsf{LKk}$ (as $K_X k_XY$).}
\end{prop}

\begin{proof} If $k_XY$ were definable in $\mathsf{LUD}$, then the formula $k_XY\wedge \neg K_XY$ would be expressible in $\mathsf{LCD}$. This formula is obviously satisfiable: just take the Counterexample constructed in Section \ref{point-cont} (involving Thomae's function). By the soundness of the proof system $\mathbf{LCD}$, this formula must then be consistent with the axioms and rules of this proof system. But by part (d) of Theorem \ref{Completeness}, the formula must be satisfiable in a pseudo-locally Lipschitz models: this is a contradiction (as pseudo-locally Lipschitz models validate the implication $k_XY\Rightarrow K_XY$, thus contradicting our formula).\end{proof}

This result suggests an obvious question:

\medskip

\noindent\textbf{Open problem} \, \emph{What is the complete logic of the language $\mathsf{LKk}$ based only on the operators $K_X\varphi$ and conditional knowability atoms $k_X Y$? Is that logic decidable?}

\subsection{A concrete example}

After all these abstract notions, some simple concrete illustration may be helpful.

\begin{exam}[Car speed detection]
An aircraft police officer determines a car's velocity $v$, by measuring the time $t$ it takes  the car to pass between two highway markings, say 1 km apart. An abstract topo-dependence model for this situation has a set $W$ of possible worlds or `states', as well as two variables $v,t: W\to (0,\infty)$, whose common range $\mathbb{D}_v=\mathbb{D}_t=(0,\infty)$ comes the standard Euclidean topology $\tau_v=\tau_t$.\footnote{We assume that it is known that the car is already in motion, so we exclude the case $v=t=0$. We could also consider a concrete model for the same situation: the set of possible states consists then of possible assignments for the basic variables $W=\{(x,y)\in (0,\infty)\times (0,\infty): x\cdot y=1\}$, with $v(x,y)=x$, $t(x,y)=y$.}  A basis of observable properties consists of the open intervals $(\alpha, \beta)\subseteq (0,\infty)$, with rational endpoints. The variables stand in a (continuous) dependence relation, for all states $s\in W$:
$$s(v)=\frac{1}{s(t)}$$
\end{exam}

Assuming that the legal speed limit is $120$ km per hour, there is a unary predicate $S$ (`speeding') denoting $S=(120,\infty)\subseteq \mathbb{D}_v$. In terms of ontic  dependence, the exact value of the speed carries the information on whether the car is speeding or not:
$$Sv\Rightarrow D_v Sv, \,\,\,\, \,\,\, \neg Sv\Rightarrow D_v \neg Sv.$$
Since the speed functionally depends on the time, we also have the (global) dependence: $$W\models D(t; v).$$ Putting these together (and using the Transfer Axiom), we see that the models also validates: $$Sv\Rightarrow D_t Sv, \,\,\,\, \,\,\, \neg Sv\Rightarrow D_t \neg Sv.$$ In other words: the \emph{exact value} of the time (needed to pass between the two markings) \emph{carries the information} on whether the car is speeding or not.

By itself, this does \emph{not} make the police officer \emph{know} that the car is speeding. But as the speeding property $S$ is an \emph{open} set,
the officer can \emph{come to know} that the car is speeding (if this is indeed the case):
$$Sv\Rightarrow K_v Sv.$$
So the fact that the car is speeding is \emph{knowable}: in principle, the officer can learn this by performing an  accurate enough  measurement of the speed $v$. However, the velocity $v$ will typically not be directly available to him for measurement. Luckily though, the  functional dependence between speed and time is \emph{continuous}, so we have the (global, and therefore known) epistemic dependence:
$$W\models K_t v$$
Once again, by putting the last two statements together, and using this time the Knowability Transfer axiom, we obtain: $$Sv\Rightarrow K_t Sv,$$ In other words: by doing an \emph{accurate enough measurement of the time} (needed to pass between the two markings), the policeman can \emph{come to know} that the car is speeding.

In contrast,  the police may never get to know that the car is \emph{not} speeding, since the complement of  $(0,\infty)-S=(0,120]$ is not an \emph{open} set: the implication
$$\neg Sv\Rightarrow K_v \neg Sv$$
is not valid in the above model, and neither is the time-dependent knowability version
$$\neg Sv\Rightarrow K_t \neg Sv$$
Of course, this epistemic difficulty only arises when the car is \emph{exactly} at the speed limit $v=120$: in this case, the car is not speeding, but the policemen will never be sure of this, no matter the accuracy of his measurements.
On the other hand, $x=120$ is the \emph{only} counterexample in this case: since the interior $Int((0,120])=(0,120)$, every state $s$ with $v(s)\not=120$ will satisfy the implication $\neg Sv\Rightarrow K_v \neg Sv$ (and thus also the implication $\neg Sv\Rightarrow K_t \neg Sv$). Hence, if the speed is \emph{not exactly} 120,  the policeman can eventually come to know whether the car is speeding or not (assuming  there is no absolute limit to the accuracy of his measurements or perception).

\vspace{1ex}

There are also further statistical aspects to this practical scenario, e.g. with aggregating results of repeated measurements, but these are beyond the scope of this paper.

\section{Epistemic Independence}

\noindent\textbf{Known dependence versus epistemic dependence}\, Recall the exact dependence operators $D_X\varphi$ and $D_XY$ of the logic $\mathsf{LFD}$, as well as the associated `epistemic' abbreviations
$$K\varphi \, :=\, D_{\emptyset} \varphi, \,\,\,\,\,\,\,\,\,\, K(Y)\, :=\, D_{\emptyset} Y$$
(both equivalent to the corresponding notions defined in $\mathsf{LCD}$), and the `global' dependence
$$D(X;Y)\, := \, K D_XY,$$
which captures a form of \emph{known dependence}: it is known that the exact value of $X$ determines uniquely the exact value of $Y$. Let us now draw some comparisons.

It is easy to see that our notion of epistemic dependence implies the existence of a knowable exact dependence, as shown by the following validities:

\vspace{-1ex}

$$K_X \varphi \, \, \, \Rightarrow \,\,\, K_X D_X \varphi  \, \, \, \Rightarrow \,\,\, D_X\varphi,$$

\vspace{-2ex}

$$K_X Y  \, \, \, \Rightarrow \,\,\, K_X D_X Y  \,\, \, \Rightarrow \,\,\, D_X Y,$$

\vspace{-2ex}

$$K(X; Y) \Rightarrow D(X; Y).$$

But the \emph{converse implications fail}: epistemic dependence $K_X Y$ is stronger than knowable dependence $K_X D_X Y$, and known global epistemic dependence $K(X;Y)$ is stronger than known global dependence $D(X;Y)$. The first notion requires as a surplus that good enough approximations of $X$ are enough to give any desired estimate of $Y$. Topologically, the distinction shows in the existence of non-continuous dependence maps.

In fact, even more extreme cases are possible: the \emph{exact dependence might be known}, while at the same time it is \emph{known that there is no knowable epistemic dependence}. This situation is topologically characterized by the following result.

\begin{prop}
The following are equivalent in topo-dependence models $\bM$:
\begin{enumerate}
\item $D(X;Y) \wedge K \neg K_XY$ holds in $\bM$ (at any/all states)
\item there exists some map $F:\mathbb{D}_X\to \mathbb{D}_Y$ s.t. $F\circ X=Y$ and $F$ has a dense set of discontinuities (i.e., $Cl\{d\in \mathbb{D}_X: F \mbox{ discontinuous at } d\}=\mathbb{D}_X$).
\end{enumerate}
\end{prop}

There is nevertheless one type of variables for which known/knowable exact dependencies are equivalent to their epistemic counterparts. A variable $X:S\to (\mathbb{D}_X, \tau_X)$ is called \textit{exact} if its value topology is discrete (i.e. $\tau_X=\mathcal{P}(\mathbb{D}_X)$).

\begin{prop} For \emph{exact} variables $X$, the following implications are valid:

\vspace{-1ex}

$$K_X \varphi \Leftrightarrow D_X \varphi$$
$$K_X Y \Leftrightarrow D_X Y$$
\end{prop}

\noindent\textbf{Complete ignorance as topological independence} Going further with the above discrepancies,  we now move to even more extreme situations where, while $X$ globally determines $Y$, it is known that no observable estimate of $X$ will give \emph{any} information about $Y$!
 The dual opposite of epistemic dependence is the case when \emph{no} observable approximation of the value of $X$ can give \emph{any} information concerning the value of $Y$. We refer to this notion as \emph{epistemic (or `topological') independence}, and denote it by $Ig_XY$ (from `ignorance').

 The notion $Ig_XY$ is not a simple negation $\neg K_XY$, but a much stronger concept,  capturing a type of \emph{zero-knowledge} situation: all $X$ is completely uninformative as far as $Y$ is concerned.

$$s\models Ig_XY \,\, \, \mbox{ iff } \,\,\, Y(U)=\mathbb{D}_Y \mbox{ for all } U\in \tau_X^S(s).$$

\vspace{1ex}

This is to be distinguished from \emph{simple ontic independence} $I_XY$ between the exact values in the sense of the original dependence logic $\mathsf{LFD}$, given by
$$s\models I_XY \,\, \, \mbox{ iff } \,\,\, Y\circ X^{-1}(X(s))=\mathbb{D}_Y.$$
\indent Again we draw some comparisons. For a start, both notions have a \emph{global counterpart}:
$$I(X; Y)\,\, :=\,\, K I_XY, \,\,\,\,\,\,\,\,\, \,\, Ig(X;Y) \,\, :=\,\, K Ig_XY.$$

\noindent The relationships between these notions are given by the following implications:
$$I_XY \, \Rightarrow \, Ig_XY, \,\,\,\,\,\,\,\,\,\,\,\, \,\,\, \, I(X;Y) \, \Rightarrow \, Ig(X;Y)$$
$$I(X;Y) \, \Rightarrow \, I_XY, \,\,\,\,\,\,\,\,\,\,\,\, \,\,\, \, Ig(X;Y) \, \Rightarrow \, Ig_XY$$

\vspace{1ex}

In particular then, epistemic independence is \emph{weaker} than ontic independence! Moreover, as in the case of probabilistic independence, \emph{both global versions are symmetric}, i.e. we have:
$$I(X;Y) \, \Leftrightarrow \, I(Y;X), \,\,\,\,\,\,\,\,\,\,\,\, \,\,\, \, Ig(X;Y) \, \Leftrightarrow \, Ig(Y;X).$$

The most interesting conceptual observation is that \emph{we can have full (known) exact dependence while at the same time having full (known) epistemic independence}.

\begin{prop}
The following are equivalent for topo-dependence models $\bM$:
\begin{enumerate}
\item $D(X;Y) \wedge Ig(X;Y)$ holds in  $\bM$ (at any/all states)
\item there exists some everywhere-surjective map $F:\mathbb{D}_X\to \mathbb{D}_Y$ s.t. $F\circ X=Y$.
\end{enumerate}
\end{prop}

Here, we used the following mathematical notion from Analysis: $F:(\mathbb{D},\tau)\to (\mathbb{D}',\tau')$ is \emph{everywhere-surjective} iff we have $F(U)=\mathbb{D}'$ for all non-empty open sets $U\in \tau\setminus\{\emptyset\}$. Everywhere-surjectivity implies everywhere-discontinuity (hence also the above-mentioned density of the set of discontinuities), but it is a much stronger condition.

\begin{exam}Each everywhere-surjective function $F:\mathbb{D}_X\to \mathbb{D}_Y$ yields an example of known exact dependence with (knowledge of) epistemic independence.
A  simple example is the Dirichlet function: $\mathbb{D}_X$ is the set of real numbers with the standard topology, $\mathbb{D}_Y$ is the set $\{0,1\}$ with the discrete topology (which is the same as the subspace topology induced by the standard topology),  with $F(d)= 0$ iff $d$ is rational, and $F(d)=1$ otherwise.
More interesting examples are the ones in which the domain and codomain of $F$ are the same (i.e., $\mathbb{D}_X = \mathbb{D}_Y$ is the common range of both variables), preferably a nice space. See \cite{Berai18} for many examples, where the common range is $Q$, $R$ or other interesting spaces. \end{exam}

\begin{exam}[Pythagorean drivers] Unlearnability in the sense of, say, the Dirichlet function may be rare, though in our earlier car speeding example, one could think of a Platonic police officer trying to book a Pythagorean driver who dematerializes at any irrational value of $t$.
\end{exam}

Unlearnability is related to the  behavior of chaotic dynamical systems, \cite{Lorenz}, where predictions of future behavior may be impossible. However, we must leave this interface of Epistemic Topology and mathematics to another occasion.

Instead, we end with a technical logical issue. It is known that the joint logic $\mathsf{LFDI}$ of ontic independence and ontic dependence is undecidable, \cite{BaBe21}, but topological indepenence had its own behavior. This suggests the following question about the purely topological component of our logic:

\medskip

\noindent\textbf{Open problem} \, \emph{Is the logic of topological independence and dependence decidable?}

\section{Uniform dependence: a stronger notion of knowability}

Up till now, we identified the approximations of the current value $X(s)$ of variable $X$ at state $s$ with open neighborhoods $U\in \tau_X(X(s))$. In this sense, our global continuous dependence $K(X; Y)$ captures known dependence, both in the sense of knowing that/whether and in the sense of knowing how:  for any desired estimate $U_Y\in \tau_Y$, the agent knows how to determine that $Y(s)\in U_Y$ if given a sufficiently good estimate $V_X\in \tau_X$ of $X(s)$. This is full how-to-knowledge: in fact, the agent knows how to find the weakest appropriate estimate $V_X$ that will determine whether  $Y(s)\in U_Y$ or not.

But the topological notion of estimation is \emph{local}: there is no way to compare the accuracies of estimates of different values $X(s)$ and $X(t)$ in different states. Things change radically if we introduce a \emph{global notion of accuracy or error} $\varepsilon$, as e.g. the one given by real-numbered open intervals $(x-\varepsilon,x+\varepsilon)$ of a given length, or more generally by open balls $B(x, varepsilon)=\{y: d(x,y)<\varepsilon)$ in a metric space. In such a setting, full how-to knowledge of a dependence would require that: for any given desired accuracy $\varepsilon_Y$ in determining the value of $Y$, the agent knows how to find an appropriate accuracy $\delta_X$, such that an estimate of $X$ with accuracy $\delta_X$ (in any state $s$) will yield an estimate of $Y$ with accuracy $\varepsilon_Y$ (in the same state $s$). Essentially, this stronger form of (global) epistemic dependency requires \emph{uniform continuity} of the underlying dependence map.

\subsection{Uniform dependence in metric spaces}

The notion of uniform dependence takes us beyond topology: its semantics requires the move to \emph{metric} dependence models.

\medskip

\par\noindent\textbf{Empirical variables over metric spaces}
As already noticed in Section \ref{LCD}, one may consider as a special case empirical variables of the form $X: S\to (\mathbb{D}_X, d_X)$, whose range of values has a \emph{metric topology} induced by a metric $d_X$ on $(\mathbb{D}_X$. From now we will restrict ourselves to such ``metric variables".
Closeness relations of the form $d(X(s),X(t))< \varepsilon$ will then give us global notions of \emph{accuracy} for the values of $X$: margins of error for measurements (or more generally, for observations) of the value of $X$.

\medskip

\par\noindent\textbf{Sets of variables as joint variables} As we already saw in Section \ref{LCD} (when giving the semantics of $\mathsf{LCD}$ on metric dependence models), a finite \emph{set} $X=\{X_i: i\in I\}$ of empirical variables over metric spaces can itself be regarded as a \emph{single metric variable}. More precisely, given such a finite set  $X=\{X_i: i\in I\}$ of metric variables $X_i:S\to (\mathbb{D}_I, d_i)$, each coming with its own metric space, we can associate to it a single metric variable, also denoted by $X$, whose set of values is
as usual $\mathbb{D}_X:= \{(x_i(s))_{i\in I}: s\in S\}\subseteq
\Pi_{i\in I} \mathbb{D}_{i}$, while the metric is (the restriction to $\mathbb{D}_X$ of) the so-called \emph{Chebyshev distance} (also known as the \emph{supremum metric}):
$$d_X (\vec{u},\vec{v}) \, :=\, sup \{d(u_i, v_i): i\in I\}, \mbox{ for all } \vec{v}=(v_i)_{i\in I}, \vec{u}=(u_i)_{i \in I}\in \Pi_{i\in I} \mathbb{D}_{i}.$$

\par\noindent\textbf{Strong epistemic dependence} We say that there is a \emph{strong (or `uniform') epistemic dependence} between empirical variables $X$ and $Y$, and write $U(X;Y)$, if for every margin of error $\varepsilon_Y>0$ for $Y$-measurements there exists some margin of error $\delta_X>0$ for $X$-measurements, such that every estimate of $X$ with accuracy $\varepsilon_X$ entails an estimate of $Y$ with accuracy $\delta_Y$:
$$\forall \varepsilon>0 \, \exists \delta>0 \, \forall s\in S \, \forall t\in S \, \left(\, d_X(X(s), X(t))< \delta \Rightarrow d_Y(Y(s), Y(t))< \varepsilon\, \right).$$
In contrast, if we apply the definition of \emph{known epistemic dependence} $K(X; Y)$ in the metric topology, and unfold the definition in terms of the underlying accuracies, we obtain:
$$\forall \varepsilon>0 \, \forall s\in S\, \exists \delta>0 \, \forall t\in S \, \left(\, d_X(X(s), X(t))< \delta \Rightarrow d_Y(Y(s), Y(t))< \varepsilon\, \right).$$
So the difference between known epistemic dependence and strong epistemic dependence is the position of the universal state quantifier $\forall s$ (before or after the existential accuracy quantifier $\exists \delta$), similarly to the analogue difference between continuity and uniform continuity. This analogy is not accidental,
witness the following result in the style of  earlier characterizations:

\begin{prop}\label{GlobalUniformCont}
Given two empirical variables $X: S\to (\mathbb{D}_X, d_X)$ and $Y: S\to (\mathbb{D}_Y, d_Y)$, the  following are equivalent:
\begin{enumerate}
\item $U(X;Y)$
\item there exists a uniformly continuous map $F_{X;Y}:(\mathbb{D}_X, d_X)\to (\mathbb{D}_Y, d_Y)$ s.t. $F\circ X= Y$ holds on $S$
\item $Y: (S, d_X^S)\to (\mathbb{D}_Y, d_Y)$ is uniformly continuous (wrt the $X$-topology on $S$ induced by the pseudo-metric $d_X^S$)
\item the identity function $id:(S, d_X^S)\to (S, d_Y^S)$ is uniformly continuous (wrt to the $X$-topology given by $D_X^S$ on its domain and the $Y$-topology given by $d_Y^S$ on its codomain).
\end{enumerate}
\end{prop}

Note that strong epistemic dependence $U(X;Y)$ is a \emph{global} notion.
We can also introduce a topologically local notion $s\models U_X Y$ of \emph{locally strong epistemic dependence}:
$$\exists \delta_0>0 \, \forall \varepsilon>0 \, \exists \delta>0 \, \forall t, w\in \mathcal{B}_X(s,\delta_0) \, \left(\, d_X(X(t), X(w))<\delta \Rightarrow d_Y(Y(t), Y(w))< \varepsilon\, \right).$$
Once again, if we denote by $\tau_X$ the topology induced by the metric $d_X$, this is equivalent to:
$$\exists O\in \tau_X(s) \mbox{ s.t. } O\models U(X; Y) \mbox{ (i.e., it holds \emph{globally} on the open \emph{subspace} $O$)}.$$

\begin{prop}
The following are equivalent:
\begin{enumerate}
\item $s\models U_X Y$
\item there exists an open neighborhood $O$ of $X(s)$, and also a uniformly continuous map $F:O \to \mathbb{D}_Y$ s.t. $F\circ X= Y$ holds on $X^{-1}(O)$
\item $Y: (S, d_X^S)\to (\mathbb{D}_Y, d_Y)$ is locally uniformly continuous around $s$
\item the identity function $id:(S, d_X^S)\to (S, d_Y^S)$ is locally uniformly continuous around $s$.
\end{enumerate}
\end{prop}

Moreover, in analogy to the case of (simple) epistemic dependence, \emph{strong epistemic dependence is inherently known}, while \emph{locally strong epistemic dependence is inherently knowable}:
$$U(X;Y)\Rightarrow K U(X;Y), \,\,\,\,\,\,\,\,\,\,\,\,\,\, \, U_XY \Rightarrow K_X U_XY.$$
However, there are some significant differences with our earlier topological notions.
For both simple and epistemic dependence, the global versions are equivalent to known local ones:

$$D(X;Y) \Leftrightarrow K D_XY, \,\,\,\,\,\,\,\,\,\,\,\,\,
K(X;Y) \Leftrightarrow K K_XY$$
Knowing a dependence (between $X$ and $Y$) \textit{is} the same as knowing that you could know the dependence (if given enough info on $X$).
However, this is \textit{not} the case for strong epistemic dependence in our new sense of uniform continuity:
$$U(X;Y) \,\,\, \mbox{is \emph{not} equivalent with} \,\,\, K U_XY$$
The explanation is that knowing how to find $Y$ from $X$ (with any desired accuracy) is \textit{not} the same as knowing that you could know how to do it, if given enough info on $X$. The second notion is in general much weaker!
We can provide counterexample to the equivalence of $U(X;Y)$ and $KU_XY$ by taking a fresh look at our car speeding scenario:

\begin{exam}[Car speeding revisited] A priori, the police officer has a known epistemic dependence $K(t, v)$, but not a strong one: $U(t,v)$ is false! For any  speed estimate $a\,\pm\, \varepsilon$, the officer can find a  $\delta$ such that after measuring time $t$ with accuracy $\delta$, (s)he will know whether the velocity is in interval $(a-\varepsilon, a+\varepsilon)$.  But, given a desired accuracy $\varepsilon$ for the speed $v$ (with no further information on the interval), the officer cannot know in advance what accuracy $\delta$ is needed in measuring time $t$ to determine velocity within $\varepsilon$.
This is because, given a fixed $\varepsilon>0$, there is \textit{no} $\delta>0$ s.t., for \textit{all} $t\in (0,\infty)$:
$$\frac{1}{t}-\frac{1}{t+\delta}<\varepsilon$$

\noindent However, $v$ is \textit{strongly knowable from $t$} in the sense of a locally strong dependence: in fact, $U_t v$ holds globally, so \textit{this strong dependence is known}:
$$K U_t v$$
Indeed, \textit{after} performing \textit{just one} reasonably accurate measurement $(a,b)$ of time $t$ (with $0<a<b<\infty$), the officer will know, for any desired velocity accuracy $\varepsilon$, how to determine the needed time accuracy $\delta$; i.e. one s.t. the following will hold for \textit{all} $t\in [a,b]$:
$$\frac{1}{t}-\frac{1}{t+\delta}<\varepsilon$$
This finishes our counterexample, but there are more interesting things going on in our car-speeding scenario: note that in this example both $U_tv$ and $K_tv$ hold globally, so the two notions are actually equivalent on this model, i.e. we have
$$U_t v\Leftrightarrow K_tv$$
Of course, the left-to-right implication always holds, but the converse gives us a kind of ``epistemic bootstrapping": when the value space we are working in is conducive to epistemic inquiry, epistemic dependence automatically gives us strong (uniform knowability).
\end{exam}

We now proceed to generalize this last observation.

\medskip

\par\noindent\textbf{Important special case: locally compact spaces} When the underlying spaces are \textit{locally compact} (i.e. every value in $D_X$ has a compact neighborhood), \emph{epistemic dependence implies locally strong dependence}: the formula
$$K_X Y\Rightarrow U_XY$$
holds on locally compact spaces. In contrast, \emph{global (known) epistemic dependence still does not imply strong epistemic dependence}: even on locally compact spaces, we still have that
$$K(X;Y)\not\Rightarrow U(X;Y).$$
Of course, the global implication $K(X;Y)\not\Rightarrow U(X;Y)$ will hold on fully \emph{compact} spaces, but those are less pervasive, being specific to very circumscribed situation: in particular, our Euclidean space is \emph{not} compact, though it is locally compact.

\medskip

\par\noindent\textbf{Interpretation: bootstrapping knowability} Since Euclidean spaces are locally compact, the natural topology of our space is arguably epistemically fertile: fit to `bootstrap knowability'. Every epistemic dependence between Euclidean variables $X$ and $Y$ is a locally strong dependence.
Indeed,  if such a dependence is knowable (i.e., if $K_XY$ holds between Euclidean variables), then it is also knowable that (given enough information about $X$) you could come to know \emph{how} to find $Y$ from $X$ with any desired accuracy. Since Euclidean variables are the ones most often encountered in empirical science, one could say that Nature makes us a free gift: it is enough to gain knowledge of an epistemic dependence to obtain potential `how-to knowability', and thus, if you wish, pass from science to engineering.

\begin{rem}[Agents and environments]  These observations have a more  general epistemological import.  Epistemic notions may `upgrade' automatically when the learning environment is favorable. What this highlights is that epistemic success involves an interplay of two parties: the nature of the knowing agents and of the environments (adversarial or helpful)  they operate in, cf. also \cite{BaPe}. This point  transcends our technical setting: in a sense, all epistemology is about the balance of the structure of  agents and the world they live in.
\end{rem}

The next class of models is a example of an even more favorable type of environment, that facilitates a kind of `epistemic super-bootstrapping':

\medskip

\par\noindent\textbf{A stronger form of epistemic bootstrapping: pseudo-locally Lipschitz models} A metric dependence model is called \emph{pseudo-locally Lipschitz} if \emph{every basic empirical variable $\mathbf{y}$ is pseudo-locally Lipschitz wrt to every other basic variable $\mathbf{x}$}; i.e.: for any two basic variables $x,y\in V$, the map $\mathbf{y}: (S, d_x^S)\to (\mathbb{D}_y, d_y)$ is pseudo-locally Lipschitz. It is easy to see that this is equivalent to requiring that the\emph{ identity map $id:(S, d_X^S)\to (S, d_Y^S)$ is pseudo-locally Lipschitz} for all \emph{non-empty sets} of variables $X, Y\subseteq V$.  As a consequence, such models afford an even stronger form of epistemic bootstrapping:

\begin{quote}
Every pseudo-locally Lipschitz model validates the implication $k_XY\Rightarrow U_XY$, for all $X,Y\subseteq V$. So in such models, \emph{our three `localized' forms of continuous dependence (point-continuous dependence $k_XY$, locally-continuous dependence $K_XY$ and locally uniformly-continuous dependence $U_XY$) are all equivalent}, for all $X,Y\subseteq V$.
\end{quote}

\begin{rem}['Paradisiacal' environments versus our world] Pseudo-locally Lipschitz models describe a kind of ``paradisiacal" epistemic situation, in which the weakest and most fragile form of dependence (conditional knowability $k_XY$, that typically is epistemically opaque) becomes in fact transparent, being bootstrapped all the way to the strongest form of know-how dependence ($U_XY$).
Unfortunately, while our Euclidean space is compact, typical Euclidean variables encountered in empirical science are not usually pseudo-locally Lipschitz with respect to each other. In this sense, \emph{although our world is epistemically fertile, it is not paradisiacal}!
\end{rem}

\subsection{The logic of uniform dependence}

We proceed to extend the earlier language of $\mathsf{LCD}$ to capture the properties of uniform dependence. We start with a very simple
extension: just add uniform dependence atoms $U(X;Y)$. The resulting logical system $\mathsf{LUD}$ is a straightforward extension of $\mathsf{LCD}$.

\begin{defi}[\textbf{Syntax of $\mathsf{LUD}$}] Given a set of predicates $Pred=\{P, Q, \ldots, \}$ and a set of variables $V=\{x, y, \ldots\}$, the Logic of Uniform Dependence ($\mathsf{LUD}$) has a language given by adding uniform dependence atoms $U(X;Y)$ to the language of $\mathsf{LCD}$; in other words, the set of formulas is recursively given by:
\vspace{-4ex}

$$
\begin{array}{c cc cc ccc cc cc cc cc cccccccc}
\varphi :: =  & P\vec{x}   &|& \neg \varphi &|& \varphi\wedge\varphi          &|& D_X \varphi &|& D_XY &|&  K_X\varphi &|& K_X Y&|& U(X;Y)
\end{array}
$$
where $X,Y\subseteq Var$ are finite sets of variables, $P\in Pred$ is a predicate and $\vec{x}\in Var^*$ is a string of variables of the same length as the arity of $P$.
\end{defi}

\begin{defi} [\textbf{Semantics of $\mathsf{LUD}$}] Given a metric dependence model $\mathbf{M}=(M,S, \mathbf{x})_{x\in Var}$ over some typed metric model $M=(\mathbb{D}_x, d_x, I)_{x\in Var}$,  we define the \emph{satisfaction} relation $s\models \varphi$ by using the usual $\mathsf{LCD}$ recursive clauses for atomic formulas and the $\mathsf{LCD}$ operators (using the metric topologies $\tau_X$ on $\mathbb{D}_X$), and interpreting the uniform dependence atoms $U(X;Y)$ in the obvious way:
$$s\models U(X;Y) \,\, \, \mbox{ iff } \, \,\,  \forall \varepsilon>0 \, \exists \delta>0\, \forall s,t\in S\, (d_X^S(s,t)< \delta\Rightarrow d_Y^S(s,t)< \varepsilon),$$
where, as in Section \ref{LCD}, we use
$$d_X^S (s,t) \,\, :=\,\, sup\{d_x(\mathbf{x}(s) , \mathbf{x}(t)): x\in X\}$$
to denote the pseudo-metric induced on the state space by the Chebyshev distance on $\mathbb{D}_X\subseteq \Pi_{x\in X} \mathbb{D}_x$.
\end{defi}

\begin{thm}\label{Completeness2} The proof system $\mathbf{LUD}$ in Table \ref{tb3} is a sound and complete axiomatization of the Logic of Uniform Dependence on metric dependence models (as well as on concrete metric models). The same proof system is also sound and complete with respect to pseudo-locally Lipschitz models. Moreover, the logic $\mathsf{LUD}$ is decidable.
\end{thm}

The completeness proof is deferred to the Appendix. Essentially, it is based on refining the model construction in the completeness proof for $\mathbf{LCD}$ to obtain a pseudo-metric model, whose associated metric dependence model happens to be pseudo-locally Lipschitz.

\begin{table}[h!]
\begin{center}
{\small
\begin{tabularx}{\textwidth}{>{\hsize=0.8\hsize}X>{\hsize=1.5\hsize}X>{\hsize=0.8\hsize}X}
\toprule
\textbf{(I)} & \textbf{All axioms and rules of the system $\mathbf{LCD}$}
 \vspace{2mm} \ \\
 \textbf{(II)} &\textbf{Axioms for uniform dependence} \vspace{2mm}:
 \ \\
($U$-Inclusion) & $U(X; Y)$, provided that $Y\subseteq X$ \vspace{1.5mm}\ \\
($U$-Additivity) & $\left(\, U(X;Y)\wedge U(X; Z)\, \right) \, \to\, U(X; Y\cup Z)$  \vspace{1.5mm}\ \\
($U$-Transitivity) & $\left(\, U(X; Y)\wedge U(Y; Z)\, \right)\to U(X; Z)$ \vspace{1.5mm}\ \\
(Uniform Dependence is Known) & $U(X;Y) \to K U(X;Y)$ \vspace{1.5mm}\ \\
(Uniformity implies Continuity) & $U(X;Y)\to K(X;Y)$ \vspace{1mm}\ \\
(Uniformity of Knowledge) & $K(Y) \to U(X;Y)$  \vspace{1mm}\ \\

\bottomrule
\end{tabularx}
}
\end{center}
\caption{The proof system $\mathbf{LUD}$.}\label{tb3}
\end{table}

\medskip

Note again the analogy of the $\mathbf{LUD}$ axioms in group (II) of Table 3 with some of the $\mathbf{LCD}$ axioms in groups (III) and (IV) of Table 2 (which are of course also included in  $\mathbf{LUD}$, within group (I) of Table 3). And indeed, if we replace $U(X;Y)$ by $K(X;Y)$ in the above axioms, we obtain valid theorems of  $\mathbf{LCD}$. Once again, this is to be expected, since uniform dependence is a strengthening of epistemic dependence, in the same sense in which epistemic dependence was a strengthening of exact dependence. But note that \emph{more properties are lost} in the transition from $K$ to $U$, in comparison to the earlier transition from $D$ to $K$. There is e.g. no analogue of the Transfer axiom for uniform dependence, since $U(X;Y)$ does not come with an accompanying modality (in contrast to $D_XY$ and $K_XY$, which come together with $D_X\varphi$ and respectively $K_X\varphi$). The reason is that Boolean variables have a discrete range of truth values, so there is no sense in `uniformizing' the meaning of $K_X\varphi$.\footnote{Things would be different if we moved to a multi-valued logic with a continuum set of values: one could then talk about ``uniform knowability" of a proposition $U_X\varphi$.}

Another difference is that, unlike the case of exact dependence and epistemic dependence (where we took the local versions $D_XY$ and $K_XY$ as primitive notions in our syntax), in $\mathsf{LUD}$ we chose to take the \emph{global} version $U(X;Y)$ as basic, and indeed we did not even include the local version $U_XY$ in our syntax. One of the reasons is that, as we already mentioned, $U_XY$ cannot define $U(X;Y)$ in any obvious way (unlike the case of exact dependence and epistemic dependence); and since we are in the end mainly interested in (global) uniformity, we must take $U(X;Y)$ as primitive.

Nevertheless, this limitation suggests extending our language of $\mathsf{LUD}$ with local atoms $U_XY$, obtaining an extended language  $\mathsf{LUD^+}$. The following result shows that this move does increase expressivity:

\begin{prop}\label{expressU}
Local uniform dependence $U_XY$ is \emph{not} expressible in the language of $\mathsf{LUD}$. So the language $\mathsf{LUD^+}$ of both local and global uniform dependence (obtained by adding local atoms $U_XY$ to $\mathsf{LUD}$) is more expressive than the language of $\mathsf{LUD}$.
\end{prop}

\begin{proof} The proof is completely similar to the corresponding inexpressivity result for $k_XY$ (Proposition \ref{expressk}), involving instead the formula $K_XY\wedge \neg U_XY$. This formula is obviously satisfiable on metric dependence models (just take any example in which $X$ is mapped into $Y$ by some function $F$ that is locally continuous around some point, but \emph{not} locally uniformly continuous around it). By the soundness of the proof system $\mathbf{LUD}$, this formula must then be consistent with the axioms and rules of this proof system. By the last part of Theorem \ref{Completeness2} (completeness of $\mathbf{LUD}$ wrt pseudo-locally Lipschitz models), if $U_XY$ were expressible in $\mathsf{LUD}$, then the formula $K_XY\wedge \neg U_XY$ would be satisfiable in a pseudo-locally Lipschitz models. But this is a contradiction, as pseudo-locally Lipschitz models validate the implication $K_XY\Rightarrow U_XY$, thus contradicting our formula.\end{proof}

\begin{rem} In themselves, the inexpressivity results claimed in Propositions \ref{expressk} and \ref{expressU} are not surprising: there was no reason to expect that the missing notions ($k_XY$ and respectively $U_XY$) were expressible in languages ($\mathsf{LCD}$ and respectively $\mathsf{LUD}$) that did not explicitly include them. But such results on expressive limitations are often hard to prove, involving intricate counterexamples. So it comes as a pleasant surprise that in these two cases, the expressivity-limitation results simply follow from the special additional features of the model constructed in the completeness proof!\end{rem}

Furthermore, Proposition \ref{expressU} leads us to yet another obvious question:

\medskip

\noindent\textbf{Open problem} \, \emph{What is the complete logic of the language $\mathsf{LUD^+}$ of both local and global uniform dependence? Is that logic decidable?}

\section{Directions for future work}

Given the analysis presented here, many natural follow-up topics arise, of which we briefly list a few.

\medskip

\par\noindent  \textbf{Further  system issues}  The systems introduced in this paper raise many further technical questions. One immediate issue concerns \emph{definability}. As we have noted, $\mathsf{LCD}$ performs a double epistemization of existing modal dependence logics. First we added topological modalities for knowability $K_X\varphi$ based on measuring values of variables $X$, and after that, we also introduced continuous dependence $K_XY$. Is the latter step necessary, or more precisely, is $K_XY$ definable in the logic $\mathsf{LFD}$ extended only with the topological interior modalities? We believe that the answer is negative, but we have not pursued the definability theory of our topological languages. Other natural open problems concern \emph{axiomatization}.

Another immediate open problem is providing a complete axiomatization of the topological logic $\mathsf{LKk}$ of conditional knowability based on the operators $K_X\varphi$ and the point-continuity atoms for conditional knowability $k_X Y$. As we mentioned, this language is strictly more expressive than $\mathsf{LCD}$. But our methods for proving completeness and decidability do not seem to work for this extended logic. The reason is that, as noted earlier, $k_XY$ is an `unstable' property: it is very sensitive to
any small changes of values. This sensitivity gives rise to technical difficulties. Is this logic decidable, or at least axiomatizable?

Other open problems that we noted were axiomatizing the logic of topological independence, as well as finding an axiomatization of the logic $\mathsf{LUD^+}$ of both local and global uniform dependence.

A final set of questions concerns \emph{computability}. For several logics in this paper, we do not know if they are decidable, and we know for none of our decidable logics what is its exact computational complexity.



\vspace{-1ex}

\paragraph{Richer languages for metric spaces} While the language of $\mathsf{LCD}$ is a good abstract vehicle for dependence in topological spaces in general, even the extended logic
$\mathsf{LUD}$ of uniformly continuous dependence seems poor in expressive power over the rich structure of metric spaces. One extension adds explicit accuracy variables and constants to the language, with corresponding modalities. One can then talk about $\eps$-closeness, compute with margins of error, and determine the complete modal logic for explicit versions of continuity and uniform continuity.\footnote{This setting covers Margin-of-Error Principles for knowledge  that have been widely discussed in the philosophical literature, cf. \cite{BaBeMOE} for a logical analysis.}

\vspace{-1.5ex}

\paragraph{Uniform spaces} There is also a qualitative approach to the epistemic surplus structure in  metric spaces by using so-called \emph{uniform spaces}, \cite{Isbell}. Here a family of reflexive symmetric `closeness relations' is given on the topologized state space, closed under  sequential compositions representing combined refinement. We can now define continuity and uniform continuity in terms of closeness intuitions (close arguments should yield close values), and study dynamic modalities for the effects of closeness refinement. We have a logical formalism for this setting and a plausible sound axiomatization, but proving completeness has so far eluded us.

\vspace{-1.5ex}

\paragraph{Computable dependence as Scott-continuity} While uniform continuity strengthens knowable dependence in terms of  state-independent approximations, another strengthenings would make knowability a form of \emph{computability}. In this interpretation, the exact value of a variable might be the limit of an infinite process of computation. The computable approximations can then be  understood as approximate values, living in the same domain as the exact value (rather than as open neighborhoods of the exact value), corresponding to  partial results obtainable at  finite stages of the computation. A natural semantic setting for this interpretation are \emph{Scott domains}, and  Domain Theory \cite{AJ95} fits  within our  topological framework, when we use the \emph{Scott topology}. Again this  allows us to use relational models in which the given topology does not match the modal semantics but the computational convergence given by Scott topology. Computability (to any desired degree of approximation) of a variable $Y$ from a variable $X$  then amounts to existence of a \emph{Scott-continuous dependence function} between them. We have a first completeness result for a logic of Scott-continuous dependence which extends the logic $\mathsf{LCD}$ with an axiom reflecting the special structure of domains..\footnote{Domain Theory is an elegant abstract approach to computability, but it does not contain actual manipulation with code. For that, we might have to provide a logical analysis in our style for Recursive Analysis, \cite{Goodstein}, getting closer to the syntax of dependence functions as laws that we can compute with.}

\medskip

We plan to explore these issues in a forthcoming paper \cite{BaBe2}. The next topics we list here concern more drastic extensions of our framework.

\vspace{-1.5ex}
\paragraph{Point-free topology} In the approach of this paper, dependence in empirical settings assumes the original set-based functional dependence and adds extra conditions of continuity. However, the very spirit of approximation might seem to make the use of underlying points debatable. One could also work in a \emph{point-free topology}, \cite{Vick96}, where there are only opens with some suitable (lattice) structure, and points only arise through representation theorems. Then, the requirement of continuity has to be replaced with conditions on maps from open approximations for values to open approximations for arguments. Again, our current topo-dependence models might then arise only in the limit as the result of some representation construction. For some first thoughts on this inn a category-theoretic perspective, see \cite{YeLi}.

\vspace{-1.5ex}
\paragraph{Learning dependence functions} Learning the actual state of the world, the paradigmatic issue in this paper, is a task when  relevant dependencies exist, and are known in an informative manner. But in science, we also want to learn \emph{the regularities themselves}. One way to go here is `lifting' our setting of single dependence models to families of these, as in the \emph{dependence universes} introduced in \cite{BaBe21}, where we may not know which precise dependencies govern our actual world. In this lifted setting, there is now a richer repertoire of relevant update actions: one can perform measurements, but one might also learn about regularities through mathematical reasoning or in other ways.

\vspace{-1.5ex}

\paragraph{Dynamical systems} Much of science is about dynamical systems transitioning to new states over time. While our logics can include temporal variables denoting points in time on a par with other variables, there is also a case for enriching our dependence logics with temporal operators describing the future and past of system evolution where the universe of states has a global stat transition function. What this would require is a topological extension of the dependence  logic $\mathsf{DFD}$ for dynamical systems in \cite{BaBeLi23}, which would link up also with  the more abstract temporal logic of topological dynamical systems of \cite{KM07}

\vspace{-1.5ex}
\paragraph{Probabilistic dependence and statistical learning} The topological view of opens as results of measurements pursued in this paper says little about how measurements are combined in scientific practice. A major challenge emerges here: how should one interface our approach  with the use of statistical techniques in Measurement Theory?

\section{Conclusion}

This paper has presented an extensive conceptual analysis of dependence in topological models for empirical inquiry. We  gave new formal languages incorporating our main notions, and proved representation, completeness and decidability results for the resulting logics of continuous and uniformly continuous dependence.

\medskip

These offerings can be viewed in several ways. First, we extended the area of Epistemic Topology with notions and systems that we believe to be of philosophical  interest.  We believe that our logics sit at the intersection of several epistemological perspectives. Our introduction of dependence merges two traditions to information and knowledge that have often been considered separate or even incompatible: Epistemic Logic and Situation Theory. Also, our  two-step topologization of dependence logic is reminiscent of the move from classical to constructive systems like Intuitionistic Logic, where open sets replace arbitrary sets, and where, equally crucially, continuous functions replace arbitrary functions. In addition, one can  see our systems as a further step toward a general logic underneath Formal Learning Theory. Finally, while these are all interfaces with formal epistemology, we believe that our questions-based perspective also links up with General Epistemology, and it might be very interesting, for instance, to see what the moves to continuous functions and even to uniformly continuous functions mean in a more informal common sense setting.

\medskip

At the same time,  our technical results also  enrich the logical study of mathematical theories, such as topological semantics for modal logics, \cite{BeBe07}. The system $\mathsf{LCD}$ is a multi-index topological modal logic which can also talk about continuous dependencies between variables. This suggests extended dependence versions of classical results in the area like the Tarski-McKinsey Theorem for completeness over standard metric spaces, \cite{MT44}. The models constructed in our completeness proof do have logic-induced metric structure, but it remains to be seen if completeness holds wrt the standard Euclidean metric.   But also new types of modal correspondence results arise. For instance, the epistemic desideratum that knowledge-that of propositions should imply uniform know-how for inquiry corresponds, in a precise modal sense, with a restriction to (locally) compact classes of metric spaces. Other connections of this  technical sort have been identified at earlier points in the paper, such as interfaces with Domain Theory and with Topological Dynamic Logic.

\medskip

Thus, Epistemic Topology, as pursued in the style of this paper, is both Epistemology endowed with topological tools and Topology/Analysis enriched with epistemological perspectives.

\bigskip
\noindent {\bf Acknowledgments} We thank the constructive audiences at many venues in Amsterdam, Pittsburgh, Paris, Beijing and Stanford where these ideas or some of their precursors have been presented since 2013. In particular, we thank Adam Bjorndahl, Konstantin Genin, Valentin Goranko, Kevin Kelly, Dazhu Li, J.~V{\"a}{\" a}n{\"a}nen, Fan Yang and Lingyuan Ye for their helpful feedback.

\bibliographystyle{apalike}
\bibliography{draft}



\bigskip

\appendix

\section{Appendix: Proofs of Completeness and Decidability for the Logical Systems $\mathsf{LCD}$ and $\mathsf{LUD}$}

It is easy to verify the soundness of the axioms and rules of $\mathsf{LCD}$ (and respectively $\mathsf{LUD}$) on topo-dependence models (and respectively metric dependence models). As for completeness, it is enough to prove (for both logics) completeness with respect to a \emph{smaller} class of models, namely metric \emph{dependence models}. However, the usual Modal Logic techniques \cite{BdRV01} (in particular, canonical model construction and filtration) will only give us (finite) \emph{non-standard} relational models (more precisely, what we will call \emph{preorder models}). To unify the two types of models (metric/topological and relational), we will introduce a further generalization, that subsumes both under a more abstract notion: \emph{topo-models}.

\subsection{Abstract topo-models}

Essentially, topo-models are abstract models, that combine a standard relational semantics for the $LFD$ modalities $D_X\varphi$ (in terms of \emph{equivalence relations}) together with a standard topological semantics for the $\mathsf{LCD}$ modalities $K_X\varphi$ in terms of topological \emph{interior operators}, while treating all atoms as just standard atomic sentences whose semantics is given by a \emph{valuation}, as usual in Modal Logic. We just need to impose a number of constraints on the relations, topology and valuation, to ensure the soundness of our axioms.

\begin{defi}[Topo-models for $\mathsf{LCD}$ and $\mathsf{LUD}$]
A \emph{topo-model for $\mathsf{LUD}$} is a relational-topological structure
$$\bM=(W, =_X, \tau^W_X, \|\bullet\|)_{X\subseteq V},$$
where: $W$ is a set of abstract states (`possible worlds'); for each set $X\subseteq V$ of variables, $=_X$ is an \emph{equivalence relation} on $W$ and $\tau^W_X\subseteq \mathcal{P}(W)$ is a \emph{topology} on $W$; and $\|\bullet\|$ is a valuation, that associates to each atom $\alpha$ of the form $Px_1\ldots, x_n$, $D_XY$, $K_XY$ or $U(X;Y)$ some subset $\|\alpha\|\subseteq W$ of the set of states, in the usual way. These ingredients are required to satisfy the following conditions:
\begin{description}
\item[(1)] if $s=_X w$ and $x_1, \ldots, x_n\in X$, then we have $s\in \|Px_1\ldots x_n\|$ iff $w\in \|Px_1\ldots x_n\|$;
\item[(2)] $=_\emptyset$ is the total relation $W\times W$, and $\tau^W_\emptyset=\{\emptyset,W\}$ is the discrete topology on $W$;\footnote{Given condition (7), the second part of this clause follows in fact from the first part.}
\item[(3)] if $s=_X w$ and $s\in \|D_XY\|$, then $s=_Y w$ and $w\in \|D_XY\|$;
\item[(4)] \emph{Inclusion}, \emph{Additivity} and \emph{Transitivity} for all atoms of the form $D_XY$, $K_XY$ or $U(X;Y)$:
\\ $\|D_XY\|=\|K_XY\|=\|U(X;Y)\|=W$, whenever $X\supseteq Y$;
\\ $\|D_X(Y\cup Z)\|= \|D_XY\|\cap \|D_XZ\|$, and similarly for $\|K_XY\|$ and $\|U(X;Y)\|$;
\\ $\|D_XY\|\cap \|D_YZ\|\subseteq \|D_XZ\|$, and similarly for $\|K_XY\|$ and $\|U(X;Y)\|$);
\item[(5)] $\|K_XY\|\in\tau^W_X$;
\item[(6)] if $s\in \|K_XY\|$, then every $X$-open neighborhood $U\in \tau^W_Y(s)$ includes some $Y$-open neighborhood
$O\in \tau^W_X(s)$ (i.e., $O\subseteq U$);
\item[(7)] if $s\in O\in \tau^W_X$ and $s=_X t$, then $t\in O$;
\item[(8)]  $\|K_XY\|\subseteq \|D_XY\|$;
\item[(9)] $\|D_\emptyset Y\|\subseteq \|K_\emptyset Y\|$;
\item[(10)] $\|U(X;Y)\|\in \{\emptyset, W\}$;
\item[(11)] $\|U(X;Y)\|\subseteq \|K_XY\|$ (or equivalently: if $\|U(X;Y)\|=W$, then $\|K_XY\|=W$).
\end{description}

\emph{Topo-models for $\mathsf{LCD}$} are simply obtained by skipping from the above definition all references to the uniform dependence atoms $U(X;Y)$ and all clauses involving them.
\end{defi}

\begin{defi}[Abstract Topo-Semantics]
As already announced, the topo-model semantics of $\mathsf{LUD}$ is specified (when given a topo-model $\bM$ for $\mathsf{LUD}$)  by using the standard relational semantics for the operators $D_X\varphi$ (seen as  modalities for $=_X$), as well as the interior semantics for the topological modalities $K_X\varphi$, while the semantics of the all atoms $\alpha$ (be they propositional $Px_1\ldots$ or dependence atoms $D_XY, K_XY$ or $U(X;Y)$) is given by the valuation $\|\alpha\|$:
$$s\models \alpha  \,\, \, \mbox{ iff } \, \,\, s\in \|\alpha\| \,\, \mbox{ (for all atoms $\alpha$) }$$
$$s\models D_X \varphi  \,\, \, \mbox{ iff } \, \,\, \forall w\in W \, (s=_X w \Rightarrow w\models w\models \varphi)$$
$$s\models K_X \varphi  \,\, \, \mbox{ iff } \, \,\, \exists U\in \tau^W_X(s)\, \forall w\in U \, w\models\varphi$$
Of course, the semantics of Boolean connectives $\neg\varphi$ and $\varphi\wedge\psi$ is the usual one.
The topo-model semantics of $\mathsf{LCD}$ is obtained (when given a topo-model $\bM$ for $\mathsf{LCD}$) by simply skipping the clause for $U(X;Y)$ in the above definition.
\end{defi}

\medskip

\par\noindent\textbf{Topo-dependence models are topo-models}
 Every topo-dependence model of the form $\bM=(M,S, \mathbf{x})_{x\in Var}$, based on a typed topological model $M=(\mathbb{D}_x,\tau_x, I)_{x\in Var}$, gives rise to an associated topo-model $\bM^\dag=(S, =_X, \tau_X^S, \|\bullet\|)$ for $\mathsf{LCD}$, obtained by putting:
$$\tau_X^S \, :=\, \mathbf{X}^{-1} (\tau_X)$$
$$s=_X w \,  \, \mbox{ iff } \,\, \mathbf{X}(s)=\mathbf{X}(w)$$
$$\|Px_1\ldots x_n\|\, \, = \,\, \{s\in W: (\mathbf{x_1}(s), \ldots, \mathbf{x_n}(s)) \in I(P)\}$$
$$\|D_XY\|\,\,  =\,\, \{s\in W: \forall w\in S \, (s=_X w \Rightarrow s=_Y w) \}$$
$$\|K_XY\|  \,\,  = \,\, \{s\in W:  \exists O\in \tau_X^S(s) \, \forall t\in O \, \forall U \in \tau_Y^S(t)\, \exists V\in \tau_X^S(t) \, V\subseteq U \},$$
where $\tau_X$ is the restriction to $\mathbb{D}_X$ of the product topology on $\Pi_{x\in X} \mathbb{D}_x$ and $\mathbf{X}$ is the joint empirical variable associated to the set of variables $X$.

It should be obvious that the topo-model semantics for $\mathsf{LUD}$ on the associated topo-model  $\bM^\dag$ is the same as the semantics on the original topo-dependence model $\bM$. Therefore, from now on we will identify a topo-dependence model with the associated topo-model, and thus include topo-dependence models in the class of topo-models. We will later be able to precisely characterize this class among all topo-models: topo-dependence models correspond to \emph{standard} topo-models for $\mathsf{LCD}$.

\begin{defi}[Topo-morphisms]\label{TopoMorphism} Given topo-models $\bM$ and $\bM'$, a map $\pi: W\to W'$ is a \emph{topo-morphism} if it satisfies the following properties:
\begin{enumerate}
\item for each set $X\subseteq Var$ of variables, $\pi$ is an interior map (i.e., a map that is both open and closed in the standard topological sense for functions) between the topologies $\tau^W_X$ and $\tau^{W'}_X$;
   \item $\pi$ is a relational p-morphism (in the usual Modal Logic sense)\footnote{A $p$-morphism (also known as a bounded morphism) is a functional bisimulation between relational structures. See \cite{BdRV01} for details.} wrt to the relations $=_X$ and the atoms $Px_1\ldots x_n$, $D_XY$, $K_XY$ and (in the case of $\mathsf{LUD}$ topo-models) $U(X;Y)$.
 \end{enumerate}

A topo-model $\bM'$ is said to be a \emph{topo-morphic image} of another topo-model $\bM$ if there exists a surjective topo-morphism from $\bM$ to $\bM'$.
\end{defi}

\begin{prop}\label{Preservation}
If $\pi: W\to W'$ is a topo-morphism between  topo-models $\bM$ and $\bM'$, then the states $s\in W$ and $\pi(s)\in W'$ satisfy the same formulas $\varphi$ of $\mathsf{LCD}$:
$$s\models_\bM \varphi \mbox{ iff } \pi(s)\models_{\bM'} \varphi$$
\end{prop}

The proof is an easy induction on the complexity of $\varphi$, and in fact it simply combines the standard proofs of the similar result for relational and topological structures for modal logic. An immediate consequence is:

\begin{cor}\label{Preservation2}
If $\bM'$ is a topo-morphic image of $\bM$, then the same formulas of $\mathsf{LCD}$, and respectively $\mathsf{LUD}$, are satisfiable/globally true in $\bM'$ and in $\bM$.
\end{cor}

\subsection{Special case 1: preorder models}

We say that a topo-model $\bM$ is \emph{Alexandroff} if all topologies $\tau_X^W$ are Alexandroff (i.e., the open sets are closed under arbitrary intersections).
Alexandroff topo-models can be given a \emph{purely relational description}, as a special case of modal-relational  models for $\mathsf{LCD}$ and $\mathsf{LUD}$, seen as purely \emph{modal languages}:

\begin{defi}[Preorder models]
A \emph{preorder model for $\mathsf{LUD}$} is a relational model $\bM=(W, =_X, \leq_X, \|\bullet\|)$, where: $W$ is a set of abstract states (`possible worlds'); $=_X$ are equivalence relations on $W$, and $\leq_X$ are preorders\footnote{A preorder is a binary relation that is reflexive and transitive.} on $W$ (one for each set $X\subseteq V$ of variables); and $\|\bullet\|$ is a valuation that associates to each atom $\alpha$ of the form $Px_1\ldots, x_n$, $D_XY$, $K_XY$ or $U(X;Y)$ some subset $\|\alpha\|\subseteq W$ of the set of states, in the usual way.  These ingredients are required to satisfy the following conditions:
\begin{description}
\item[(1)] if $s=_X w$ and $x_1, \ldots, x_n\in X$, then we have $s\in \|Px_1\ldots x_n\|$ iff $w\in \|Px_1\ldots x_n\|$;
\item[(2)] both $=_\emptyset$ and $\leq_\emptyset$ are equal to the total relation $W\times W$;\footnote{Once again, the second part of this clause follows in fact from the first part together with condition (7).}
\item[(3)] if $s=_X w$ and $s\in \|D_XY\|$, then $s=_Y w$ and $w\in \|D_XY\|$;
\item[(4)] atoms $D_XY$, $K_XY$ and $U(X;Y)$ satisfy Inclusion, Additivity and Transitivity;
\item[(5)] if $s\leq_X w$ and $s\in \|K_X Y\|$, then $w\in \|K_X Y\|$;
\item[(6)] if $s\leq_X w$ and $s\in \|K_X Y\|$, then $s\leq_Y w$;
\item[(7)] if $s=_X w$, then $s\leq_X w$;
\item[(8)]  $\|K_XY\|\subseteq \|D_XY\|$;
\item[(9)] $\|D_\emptyset Y\|\subseteq \|K_\emptyset Y\|$;
\item[(10)] $\|U(X;Y)\|\in \{\emptyset, W\}$;
\item[(11)] $\|U(X;Y)\|\subseteq \|K_XY\|$
\end{description}

\emph{Preorder models for } are again simply obtained by skipping from the above definition all references to the atoms $U(X;Y)$ and all clauses involving them.
\end{defi}

\begin{defi}[Preorder Semantics]
The preorder semantics for $\mathsf{LCD}$ uses the standard modal semantics for $D_X$ (as relational modalities for $=_X$) and $K_XY$ (as relational modalities for the preorders $\leq_X$), while the semantics of  all atoms $\alpha$ (be they propositional $Px_1\ldots$ or dependence atoms $D_XY, K_XY$ or $U(X;Y)$) is again given by the valuation $\|\alpha\|$:
$$s\models \alpha  \,\, \, \mbox{ iff } \, \,\, s\in \|\alpha\| \,\, \mbox{ (for all atoms $\alpha$) }$$
$$s\models D_X \varphi  \,\, \, \mbox{ iff } \, \,\, \forall w\in W (s=_X w \Rightarrow w\models w\models \varphi)$$
$$s\models K_X \varphi  \,\, \, \mbox{ iff } \, \,\, \forall w\in W (s\leq_X w \Rightarrow w\models w\models \varphi)$$
and again the semantics of the Boolean connectives is given by the usual clauses.
The preorder semantics of $\mathsf{LCD}$ is  obtained by  skipping the clause for $U(X;Y)$ in the above definition.
\end{defi}

Note that preorder models are just a special class of standard relational models for the extended language with atoms of the form $P\mb{x}$ and $K_X Y$ (as well as $U(X;Y)$ in the case of $\mathsf{LUD}$). This means that we can apply to this semantics all the known concepts and results of Modal Logic \cite{BdRV01}, concerning $p$-morphisms, bisimulations, filtrations, etc.

\begin{prop}\label{equiv}
Preorder models for $\mathsf{LCD}$ (and respectively $\mathsf{LUD}$) are the same as Alexandroff topo-models for $\mathsf{LCD}$ (and respectively $\mathsf{LUD}$). A map between two preorder models is a topo-morphism between the corresponding topo-models iff it is a relational p-morphism
wrt to both relations $=_X$ and $\leq_X$ and all the atoms ($Px_1\ldots x_n$, $D_XY$, $K_XY$, as well as $U(X;Y)$ in the case of $\mathsf{LUD}$).
\end{prop}

\medskip

\par\noindent\textbf{Equivalence between preorder models and Alexandroff topo-models}
As a consequence of Proposition \ref{equiv}, we will \emph{identify Alexandroff topo-models with the corresponding preorder models} (and later we will show that they  correspond precisely to \emph{standard} preorder models).

\subsection{Special case 2: pseudo-metric models}

We now introduce a special class of abstract topo-models that generalizes our intended metric dependence models for $\mathsf{LUD}$.

\begin{defi}[Pseudo-metric models for $\mathsf{LCD}$ and $\mathsf{LUD}$]
Let $\bM=(S, =_X, \tau^W_X, \|\bullet\|)_{X\subseteq V}$ be a topo-model for $\mathsf{LUD}$. We say that $\bM$ is a \emph{pseudo-metric model for $\mathsf{LUD}$} if the topologies $\tau_X^S$ are induced by pseudo-metrics\footnote{Recall that a \emph{pseudo-metric} on a set $S$ is a map $d: S\times S\to R$ satisfying three conditions, for all $s,w,t\in S$: $d(s,s)=0$; $d(s,w)=d(w,s)$; $d(s,w)\leq d(s,t)+ d(t,w)$.}
$d_X^S$ on $S$, i.e. we have that $$\tau_X= \mbox{ the topology generated by the basis } \{\mathcal{B}_X(s, \eps): s\in S, \eps >0\},$$
where $\mathcal{B}_X(s, \eps):=\{t\in S; d_X^S(s,t) <\eps\}$. \emph{Pseudo-metric models for $\mathsf{LCD}$} are again obtained by skipping the clauses involving $U$ (i.e. they are just topo-models for $\mathsf{LCD}$ whose topologies are induced by pseudo-metrics).
When we want to fix some specific pseudo-metrics, we exhibit them in a explicit manner, denoting the given pseudo-metric model by $\bM=(S, =_X, d^S_X, \|\bullet\|)_{X\subseteq V}$.
\end{defi}

We now recover our metric dependence models as a special case of pseudo-metric topo-models.

\medskip

\par\noindent\textbf{Metric dependence models are pseudo-metric topo-models} Every metric dependence model of the form $\bM=(M,S, \mathbf{x})_{x\in Var}$, based on a typed metric model $M=(\mathbb{D}_x,d_x, I)_{x\in Var}$, gives rise to a pseudo-metric topo-model $\bM^\dag$ for $\mathsf{LUD}$, obtained by putting:
$$d_X^S(s,w)\, \,\,\, =\,\,\, d_X(\mathbf{x}(s), \mathbf{x}(w)), \, \, \mbox{ where $d_X(\vec{u}, \vec{v})=sup\{d_x(u_x, v_x): x\in X\}$}
$$
$$s=_X w \,  \, \mbox{ iff } \,\, \mathbf{x}(s)=\mathbf{x}(w) \,\, (\mbox{ iff } \,\, d_X^S(s,w)=0)$$
$$\|Px_1\ldots x_n\|= \{s\in W: (\mathbf{x_1}(s), \ldots) \in I(P)\}$$
$$\|D_XY\| \,\,\, =\,\,\, \{s\in W: \forall t\in S \, (d_X^S(s,t)=0 \Rightarrow d_Y^S(s,t)=0)\}$$
$$\|K_XY\|  \,\,\, =\,\,\, \{s\in W: \exists \delta_0>0 \, \forall t\in \mathcal{B}_X(s, \delta_0) \, \forall \varepsilon>0 \, \exists \delta>0 \,
 \mathcal{B}_X(s, \delta_0)\cap \mathcal{B}_X(t,\delta)\subseteq \mathcal{B}_Y(t, \varepsilon)\}$$
$$\|U(X;Y)\|=W \,\,\, \mbox{ iff } \,\,\, \forall \varepsilon>0 \, \exists \delta>0 \, \forall t,w\in S \, (d_X^S(t,w) < \delta \Rightarrow d_Y^S(t,w) < \varepsilon)$$
$$\|U(X;Y)\|=\emptyset, \mbox{ otherwise.}$$
It should be obvious that the topological structure of this associated pseudo-metric model matches that of the associated topo-model (as introduced in the previous section), which explains our use of the same notation $\bM^\dag$ for both constructions. It should be equally obvious that the topo-model semantics for $\mathsf{LUD}$ on the associated pseudo-metric model $\bM^\dag$ is the same as the semantics on the original metric dependence model $\bM$. So, as in the case of general topo-models, from now on we will identify a metric dependence model with the associated pseudo-metric model, and thus include metric dependence models in the class of pseudo-metric models (and later show that they precisely correspond to \emph{standard} pseudo-metric models.)

\subsection{Standard models and their representation}

\begin{defi}[Standard topo-models]
An $\mathsf{LCD}$ topo-model  $\bM=(W, =_X, \tau^W_X, \|\bullet\|)$ is \emph{standard} if we have, for all $s\in W$ and $X,Y\subseteq Var$:
\begin{description}
\item[(St$_1$)] $\tau^W_X$ is the topology generated by $\bigcup\{\tau_x^W: x\in X\}$ (where we skipped set brackets for singletons as usual, writing $\tau_x^W$ for $\tau_{\{x\}}^W$);
\item[(St$_2$)] $=_X$ is the topological indistinguishability relation, given by: $s=_X w$ iff $\tau_X^W(s)=\tau_X^W(w)$;
    \item[(St$_3$)] $s\in \|D_XY\|$ iff $\forall w\in W\, ( s=_X w \Rightarrow s=_Y w)$;
\item[(St$_4$)] $s\in \|K_XY\|$ iff $\exists O_0\in \tau^W_X(s) \, \forall t\in O_0 \, \forall U\in \tau^W_Y(t)\, \exists O\in \tau^W_X(t) \, O\subseteq U$.
\end{description}

\smallskip

\par\noindent\textbf{Simplified presentation of standard topo-models} It should be clear that a standard topo-model is uniquely determined by its set of points $W$, its basic topologies $\tau_x^W$ and the valuation of the atoms of the form $Px_1\ldots x_n$, so it can be identified with the structure
$$(W,\tau_x^W, \|\bullet\|_{x\in V}),$$
where the valuation ranges \emph{only} over atoms $Px_1\ldots x_n$), while the topologies $\tau^W_X$ are defined using clause (St$_1$) above,
the relations $=_X$ are defined by the clause (St$_2$) above (i.e., as topological indistinguishability relations for $\tau_X^W$), and the semantics of the atoms $D_XY$ and $K_XY$ is given by the clauses (St$_3$) and (St$_4$) above.
\end{defi}

\begin{prop}\label{Standard1}
For every topo-dependence model $\bM$, the associated topo-model $\bM^\dag$ is standard.
\end{prop}

The proof of this is immediate. The converse also holds:

\begin{prop}\label{Standard2}
Each standard topo-model $\bM=(W, =_X, \tau^W_X, \|\bullet\|)$ induces a topo-dependence model $\bM^\flat$, based on the same set $S:=W$ of states, where each empirical variable $\mb{x}$ has as its range $\mathbb{D}_x$ is given by
$$\mathbb{D}_x:= \{(x,[w]_x): w\in W\}, \,\, \mbox{ where $[w]_x:=\{s\in W: w=_x s\}$};$$
the map $\mb{x}: S=W\to  \mathbb{D}_x$ is given by
$$\mb{x}(w):= (x, [w]_x);$$
the topology $\tau_x$ on $\mathbb{D}_x$  is given by
$$\tau_x :=\{ \mb{x}(O): O\in \tau^x\}, \,\, \mbox{ where $\mb{x}(O):=\{\mb{x}(w):w\in O\}$};$$
and finally, the interpretation $I(P)$ of each predicate symbol is given by
$$I(P):=\{(\mb{x_1} (w), \ldots, \mb{x_n}(w): w\in W, x_1, \ldots, x_n\in V \mbox{ s.t. } w\in\| Px_1\ldots x_n\| \}.$$
Moreover, the abstract topo-semantics on $\bM$ agrees with the topo-dependence semantics on the associated abstract topo-dependence model.
\end{prop}

\medskip

\par\noindent\textbf{Equivalence between standard topo-models and topo-dependence models}
Based on Facts \ref{Standard1} and \ref{Standard2}, from now on we will \emph{identify topo-dependence models with the corresponding standard topo-models}.

\medskip

\par\noindent\textbf{Special case: standard preorder models} A \emph{standard preorder model} is just a preorder model that is standard when considered as a topo-model (i.e. it satisfies conditions (St$_1$)-(St$_5$) in the definition of standard topo-models). It is useful to give a more direct characterization of the standard preorder models in purely relational terms:

\begin{prop}\label{Standard3}
A structure $\bM=(W, =_X, \leq_X,\|\bullet\|)_{X\subseteq V}$, consisting of a set $W$ of states, equivalence relations $=_X$ and preorder relations $\leq_X$  on $W$, as well as a valuation for the language of $\mathsf{LCD}$, is a standard preorder model iff it satisfies the following conditions, for all $s\in W$ and $X,Y\subseteq Var$:
\begin{description}
\item[($0$)]  if $s=_X w$ and $x_1, \ldots, x_n\in X$, then we have $s\in \|Px_1\ldots x_n\|$ iff $w\in \|Px_1\ldots x_n\|$;
\item[($1$)] $\leq_X$ is the intersection of all relations $\leq_x$ with $x\in X$ (where for singletons we again wrote $\leq_x$ for $\leq_{\{x\}}$);
\item[($2$)] $s=_X w$ iff both $s\leq_X w$ and $w\leq_X s$;
\item[($3$)] $s\in \|D_XY\|$ iff $\forall w\in W\, ( s=_X w \Rightarrow s=_Y w)$;
\item[($4$)] $s\in \|K_XY\|$ iff $\forall t,w\in W \, (s\leq_X t\leq_X w \Rightarrow t\leq_Y w)$.
\end{description}
\end{prop}

\smallskip

\par\noindent\emph{\textbf{Simplified presentation of standard pre-models} It should again be clear that a standard pre-model is uniquely determined by its set of points $W$, its basic preorders $\leq_x$ and the valuation of the atoms of the form $Px_1\ldots x_n$, so it can be identified with the structure
$$(W,\leq_x, \|\bullet\|)_{x\in V}$$
(where the valuation ranges \emph{only} over atoms $Px_1\ldots x_n$), with the \emph{only requirement} on these structures being condition ($0$) of Proposition \ref{Standard3}. In this simplified presentation, all the other components are definable by standardness:
the relations $\leq_X$ are defined by clause ($1$) in Proposition \ref{Standard3}, the
relations $=_X$ are defined by the clause ($2$) in the same Proposition, and the valuation of the atoms $D_XY$ and $K_XY$ is given by the clauses ($3$) and ($4$) in this Proposition.}

\medskip

\par\noindent\textbf{Equivalence between standard preorder models and Alexandroff topo-dependence models}
As a special case of the above-mentioned equivalence between standard topo-models and topo-dependence models, we can now \emph{identify standard preorder models with the corresponding  Alexandroff topo-dependence models}.

\medskip

\begin{defi}[Standard pseudo-metric models]
A pseudo-metric model  of the form $\bM=(W, =_X, d^W_X, \|\bullet\|)$ for $\mathsf{LUD}$ is \emph{standard} if it is standard as a topo-model (i.e. it satisfies all conditions (St$_1$)-(St$_4$) in the definition of standard topo-models), and in addition we have, for all $s\in W$ and $X,Y\subseteq Var$:
\begin{description}
\item[(St$_5$)]  $s\in \|U(X;Y)\|$ iff $\forall \varepsilon>0 \, \exists \delta>0 \, \forall s,t\in S \, (d_X^S(s,t) < \delta \Rightarrow d_Y^S(s,t) < \varepsilon)$.
\end{description}
\end{defi}

Note that, unlike in the case of preorder models, standard pseudo-metric models are \emph{not} the same as pseudo-metric models whose underlying topo-model is standard: one needs the additional condition (Stg$_5$). It is useful to have a direct characterization of standard pseudo-metric models in terms of the pseudo-metrics $d_X^W$:

\begin{prop}\label{Standard4}
A structure  $\bM=(W, =_X, d^W_X, \|\bullet\|)$ consisting of a set $W$ of states, equivalence relations $=_X$ and pseudo-distances $d_X^W$ on $W$, as well as a valuation for the language of $\mathsf{LUD}$, is a standard $\mathsf{LUD}$ pseudo-metric model iff it satisfies the following conditions, for all $s\in W$ and $X,Y\subseteq Var$:
\begin{description}

\item[($0$)]  if $s=_X w$ and $x_1, \ldots, x_n\in X$, then we have $s\in \|Px_1\ldots x_n\|$ iff $w\in \|Px_1\ldots x_n\|$;

\item[($1$)] $d^W_{X}$ is the Chebyshev pseudo-metric induced by the basic pseudo-metrics $d_x^W$:
$$d_X^W(s, w) \, :=\, max\{d_x^W(s,w): x\in X\},$$
with the convention that $max \,\emptyset=0$;

\item[($2$)] pseudo-metric indistinguishability conincides with $X$-equality: $s=_X w$ iff $d_X^W(s,w)=0$;
\item[($3$)] $s\in \|D_XY\|$ iff $\forall w\in W \, ( d(_X^W(s, w)=0 \Rightarrow d_Y^W(s,w)=0)$;
\item[($4$)] $s\in \|K_XY\|$ iff  $\exists \delta_0>0 \, \forall t\in \mathcal{B}_X(s, \delta_0) \,
\forall \varepsilon>0 \, \exists \delta>0 \, \forall w\in \mathcal{B}_X(s, \delta_0)\,
(d_X^W (t,w) < \delta \Rightarrow d_y^W(t,w) < \varepsilon)$;
\item[($5$)]  $s\in \|U(X;Y)\|$ iff $\forall \varepsilon>0 \, \exists \delta>0 \, \forall t,w\in W \, (d_X^W(t,w) < \delta \Rightarrow d_Y^W(t,w) < \varepsilon)$.
\end{description}
Standard $\mathsf{LCD}$ pseudo-metric models can be similarly characterized by restricting the valuation to the language of $\mathsf{LCD}$, and dropping the last clause.
\end{prop}

\smallskip

\par\noindent\emph{\textbf{Simplified presentation of standard pseudo-metric models} Again, a standard pseudo-metric model is uniquely determined by its set of points $W$, its basic pseudo-distances $d_x^W$ and the valuation of the atoms of the form $Px_1\ldots x_n$, so it can be identified with the structure
$$(W, d_x^W, \|\bullet\|)_{x\in V},$$
where the valuation ranges \emph{only} over atoms $Px_1\ldots x_n$), while
the pseudo-distances $d_X^W$ are defined by clause ($1$) in Proposition \ref{Standard4}, the
relations $=_X$ are defined by the clause ($2$) in the same Proposition, and the valuation of the atoms $D_XY$, $K_XY$ and $U(X;Y)$ is given by the clauses ($3$), ($4$) and ($5$) in this Proposition.}

\medskip

We now extend Propositions \ref{Standard1} and \ref{Standard2} to pseudo-metric models:

\begin{prop}\label{Standard5}
For every metric dependence model $\bM$, the associated pseudo-metric model $\bM^\dag$ is standard.
\end{prop}

As in the case of topo-dependence models, the proof of this is immediate, and the converse also holds:

\begin{prop}\label{Standard6}
Every standard pseudo-metric model (given in its simplified presentation) $\bM=(W, d^W_x, \|\bullet\|)$ induces a metric dependence model $\bM^\flat$, based on the same set $S:=W$ of states, together with: each empirical variable $\mb{x}$ has as its range $\mathbb{D}_x$ is given by
$$\mathbb{D}_x:= \{(x,[w]_x): w\in W\}, \,\, \mbox{ where $[w]_x:=\{s\in W: w=_x s\}=
\{s\in W: d_x^W(w, s)=0\}$};$$
the map $\mb{x}: S=W\to  \mathbb{D}_x$ is given by
$$\mb{x}(w):= (x, [w]_x)$$
(this map is obviously surjective);
the metrics $d_x$ on $\mathbb{D}_x$  are given by
$$d_x (\mb{x}(s), \mb{x}(w)):= d_x^W(s,w)\,\,\footnote{This can be checked to be well-defined, based on the definition of $\mb{x}$, its surjectivity and the following easy consequence of the triangle inequality: if $d_x^W(s,s')=d_x^W(w,w')=0$, then $d_x^W(s,w)=d_x^W(s', w')$).};$$

and finally, the interpretation $I(P)$ of each predicate symbol is given by
$$I(P):=\{(\mb{x_1} (w), \ldots, \mb{x_n}(w): w\in W, x_1, \ldots, x_n\in V \mbox{ s.t. } w\in\| Px_1\ldots x_n\| \}.$$
Moreover, the pseudo-metric semantics on $\bM$ agrees with the metric dependence semantics $\bM^\flat$.
\end{prop}

\medskip

\par\noindent\textbf{Equivalence between standard pseudo-metric models and metric dependence models}
Based on Facts \ref{Standard5} and \ref{Standard6}, we will \emph{identify metric dependence models with the corresponding standard pseudo-metric models}.

\subsection{Completeness wrt finite preorder models}

\begin{prop}\label{PreFMP}
The system $\textbf{LCD}$ is complete for finite preorder models for $\mathsf{LCD}$, and similarly the system $\textbf{LUD}$ is complete for finite preorder models for $\mathsf{LUD}$. \end{prop}

\begin{proof}  We only give the proof for the harder case of $\mathsf{LUD}$, and specify which steps should be skipped to obtain the proof for $\mathsf{LCD}$.
The proof uses a variation of the standard Modal Logic method of \emph{finite canonical models}, based on combining the techniques of canonical models and  filtration.

Let $\varphi_0$ be our fixed consistent formula  (where `consistent' refers to the appropriate proof system, $\textbf{LCD}$ and respectively $\textbf{LUD}$), and let $V$ be the \emph{finite} set of variables that are actually occurring in $\varphi$. Take now $\Phi$ to be the \emph{smallest set of formulas in $Fml_\Delta$, with the following closure properties}:
\begin{itemize}
\item $\varphi_0\in \Phi$;
\item $\Phi$ is closed under subformulas and single negations $\sim\varphi$;
\item $D_X Px_1 \ldots x_n$, $D_X D_X Y, K_X K_XY, D_X K_XY, U(X;Y)\in \Phi$, for all sets $X,Y\subseteq V$ and all $x_1, \ldots, x_n\in V$ (where the closure clause for $U(X;Y)$ can be skipped
in the case of $\mathsf{LCD}$);
\item if $K_X\psi\in \Phi$, then $D_X K_X\psi\in \Phi$.
\end{itemize}
We can easily show that \emph{$\Phi$ is finite}.

We will build a finite preorder model for $\varphi_0$ consisting of (some) maximally consistent subsets of $\Phi$ (where again  `consistent' refers to the appropriate proof system, $\textbf{LCD}$ and respectively $\textbf{LUD}$). First, we fix  a maximally consistent subset $s_0$ of $\Phi$ s.t. $\varphi_0\in s_0$. (The existence of $s_0$ follows by the standard enumeration argument for the finite case of the Lindenbaum Lemma.) Second, for each set $X\subseteq V$, we define two relations $=_X$ and $\leq_X$ between maximally consistent subsets of $s, w\subseteq \Phi$, by putting:
$$s=_X w  \,\,\, \mbox{ iff } \,\, \, \forall (D_Y\psi)\in \Phi \left(\, \{D_Y\psi, D_X Y\}\subseteq s \,\, \Leftrightarrow \, \,
 \{D_Y\psi, D_X Y\}\subseteq w \, \right);$$
$$s\leq_X w  \,\,\, \mbox{ iff } \,\, \, \forall (K_Y\psi)\in \Phi \left(\, K_Y\psi, K_X Y\}\subseteq s \,\, \Rightarrow \, \,
 \{K_Y\psi, K_X Y\}\subseteq w \, \right).$$

We define now a finite preorder model
$$\bM=(W, =_X^W, \leq_X^W, \|\bullet\|), \,\, \mbox{ where: }$$
\begin{itemize}
\item $W \,\, :=\,\, \{w\subseteq \Phi: w \mbox{ maximally consistent subset with } s_0=_\emptyset w\}$;
\item the relations $=_X^W$ and $\leq_X^W$ are just the
\emph{restrictions of the above-defined relations} $=_X$ and $\leq_X$ to the set $W$;
\item $\|\bullet\|$ is the \emph{canonical valuation} for all atoms $P\mb{x}$, $D_XY$, $K_XY$ and $U(X;Y)$, according to which by definition $w$ satisfies an atom iff the atom belongs to $w$.\footnote{More precisely, we put:
$\|P \mb{x}\| =\{s\in W: P\mb{x}\in s\}$; $\|D_XY\|=\{s\in W: D_XY\in s\}$;
$\|K_XY\|=\{s\in W: K_XY\in s\}$; and $\|U(X;Y)\|=\{s\in W: U(X;Y)\in s\}$, i.e. $\|U(X;Y)\|=W$ iff $U(X;Y)\in s_0$, and $\|U(X;Y)\|=\emptyset$ otherwise). We skip the clause for $U(X;Y)$ in the case of $\mathsf{LCD}$.}
\end{itemize}

Just by looking at the above definitions of $=_X$ and $\leq_X$, it should be obvious that the relations $=_X$ (and so also their restrictions) $=_X^W$ are equivalence relations, and that the relations $\leq_X$ (and so also their restrictions $\leq_X^W$) are preorders. This can be further strengthened to:

\begin{lem}\label{Preorder Lemma}(\emph{Preorder Model Lemma}) $\bM$ is a preorder model.
\end{lem}

\begin{proof}
We need to check that $\bM$ satisfies conditions (1)-(11) in the definition of preorder models.

For (1), let $s,w\in W$ s.t. $s=_X^W w$ and $x_1, \ldots, x_n\in X$, and assume that $s\in \|Px_1\ldots x_n\|$. By the definition of the canonical valuation, this means that $Px_1\ldots x_n\in s$. By the Determined Atoms axiom of $\mathsf{LFD}$ and the closure conditions on $\Phi$ (as well as the closure of maximally consistent subsets of $\Phi$ under modus ponens in $\Phi$), we have $D_{\{x_1, \ldots, x_n\}} Px_1\ldots x_n\in s$. Using the fact that $x_1, \ldots, x_n\in X$ and the Inclusion axiom (and again the closure conditions on $\Phi$), we see that $D_X\{x_1, \ldots x_n\}\in s$, and hence $\{D_{\{x_1, \ldots, x_n\}} Px_1\ldots x_n, D_X\{x_1, \ldots x_n\}\}\subseteq s$. By the definition of $=_X$ applied to $s=_Xw$, we obtain that $\{D_{\{x_1, \ldots, x_n\}} Px_1\ldots x_n, D_X\{x_1, \ldots x_n\}\}\subseteq w$, and using the Factivity of $D$ (axiom $T$) and the closure of $\Phi$ under subformulas, we conclude that $Px_1\ldots x_n\in w$, which by the definition of canonical valuation gives us that $w\in \|Px_1\ldots x_n\|$, as desired.

For (2): By the definition of $W$, $=_\emptyset^W$ is indeed the total relation on $W$.\footnote{As already remarked, the second part of condition (2) on preorder models is redundant, following from the first part together with condition (7).}

For (3): Assume $s=_X^W w$ and $s\in \|D_XY\|$, and so $D_XY\in s$ (by the canonical valuation), hence $D_XY\in w$ (by the definition of $=_X$), and thus $w\in \|D_XY\|$ (again by the canonical valuation). All that is left is to check that $s=_Y^W w$. For this, let $(D_Z\psi)\in \Phi$, and we need to show that $\{D_Z\psi, D_Y Z\}\subseteq s \Leftrightarrow \{D_Z\psi, D_Y Z\}\subseteq w$. We check the left-to-right implication (since the converse then follows by the symmetry of $=_X^W$): assume that $\{D_Z\psi, D_Y Z\}\subseteq s$. Since $D_XY, D_YZ\in s$, we can use the $D$-Transitivity axiom to obtain $D_XZ\in s$. So we have $\{D_Z\psi, D_XZ\} \subseteq s$, which together with $s=_X^W w$ gives us that  $\{D_Z\psi, D_XZ\} \subseteq w$, hence $D_Z\psi\in w$. On the other hand, from $D_YZ\in s$ we derive $D_Y D_YZ\in s$ by the axiom of Determined Dependence, and so we have $\{D_Y D_YZ, D_XY\}\subseteq s$, which together with $s=_X^W w$ gives us that $\{D_Y D_YZ, D_XY\}\subseteq w$, hence $D_Y D_Y Z\in w$, and thus by Factivity (axiom $T$) we have $D_Y Z\in w$. Putting these all together, we conclude that $\{D_Z\psi, D_YZ\}\subseteq w$, as desired.

Condition (4) follows from the definition of canonical valuation, together with the axioms of Inclusion, Additivity and Transitivity for $D_XY$, $K_XY$ and $U(X;Y)$.

For (5) and (6): The proof is completely similar to the above proof of the two conclusions of condition (3), using the definition of $\leq_X$ (instead of $=_X$) and the corresponding axioms for $K_X$ ($K$-Veracity, $K$-Transitivity and Knowable Dependence) instead of the axioms for $D_X$ ($D$-Factivity, $D$-Transitivity and Determined Dependence).

For (7): Assume $s=_X^W w$. To show that $s\leq_X^W w$, let $K_Y\psi\in \Phi$ s.t. $\{K_Y\psi, K_XY\}\subseteq s$, and we need to prove that $\{K_Y\psi, K_XY\}\subseteq w$. To show that $K_Y\psi \in w$, note that $K_Y\psi \in s$ implies $D_Y K_Y\psi\in s$ (by the $\textbf{LCD}$ axioms of $K_Y$-Positive Introspection and Knowable Determination, as well as our closure conditions on $\Phi$), and similarly $K_XY\in s$ implies $D_XY\in s$ (by the $\textbf{LCD}$ axiom of Knowable Dependence), and so we obtain $\{D_Y K_Y\psi, D_XY\}\subseteq s$. From this and $s=_X w$, we infer that $\{D_Y K_Y\psi, D_XY\}\subseteq w$, and thus in particular $D_Y K_Y\psi\in w$, which by the Factivity axiom gives us $K_Y\psi\in w$, as desired. Next, we also need to show that $K_XY\in w$: for this, note that $K_XY\in s$ implies $K_XK_XY\in s$ (by Knowability of Epistemic Dependence), hence we also have $D_X K_XY\in s$ (by Knowable Dependence and the closure conditions on $\Phi$), and thus $\{D_X K_XY, D_XX\}\subseteq s$ (by $D$-Inclusion). From this together with $s=_X w$, we obtain that  $\{D_X K_XY, D_XX\}\subseteq w$, so in particular $D_X K_XY\in w$, from which we infer $K_XY\in w$ (by Factivity). Putting all these together, we conclude that $\{K_Y\psi, K_XY\}\subseteq w$.

Conditions (8), (9) and (10) follow directly from the definition of canonical valuation, together with the Connecting Axioms of $\textbf{LCD}$.
\end{proof}

\begin{lem}\label{Diamond Lemma}(\emph{Diamond Lemma})
For all states $s\in W$ and formulas $\sim\varphi$ (where we recall that $\sim$ is single negation), we have:
\begin{enumerate}
\item If $D_X\varphi\in \Phi$, then
 $(\sim D_X\varphi)\in s\in W$ iff there exists $w\in W$ with $s=_X w$ and $(\sim \varphi)\in w$;
\item If $K_X\varphi\in \Phi$, then $(\sim K_X\varphi)\in s\in W$ iff there exists $w\in W$ with $s\leq_X w$ and $(\sim \varphi)\in w$.
\end{enumerate}
\end{lem}

\begin{proof} The proof of this result goes along standard lines, cf. \cite{BdRV01}: the non-trivial (left-to-right) implications are shown by constructing a consistent set of formulas in $\Phi$, that includes $\sim\varphi$ and that forces every $\Phi$-maximally consistent superset to be in the appropriate relation with $s$.\footnote{In fact, the proof of the first item was given in detail in \cite{BaBe21}, while the proof of the second item is completely analogous (using the appropriate $K_X$-axioms instead of the corresponding $D_X$-axioms).}
\end{proof}

As an immediate consequence, we obtain:

\begin{lem}\label{Truth Lemma}(\emph{Truth Lemma}) For every $s\in W$, we have that, for all formulas $\varphi\in \Phi$:
$$s\models_\bM \varphi \mbox{ iff } \varphi\in s.$$
\end{lem}

\begin{proof} The proof proceeds by induction on $\varphi$, as usual. The atomic case is simply enforced by our choice of the canonical valuation, and the Boolean connective cases are trivial, while the modal cases follow from the above Diamond Lemma.
\end{proof}

Putting these together, we obtain the desired conclusion, thus showing the completeness of our logics for finite preorder models.
\end{proof}

As a consequence, the set of theorems of $\textbf{LCD}$ is recursive, and the same holds for the set of theorems of $\textbf{LUD}$
\begin{prop}\label{PreFMP2}
The logics axiomatized by the proof systems $\textbf{LCD}$ and $\textbf{LUD}$ are decidable.
\end{prop}

Of course, this does not yet prove the decidability of our intended logics $\mathsf{LCD}$ and $\mathsf{LUD}$ (but only of the given proof systems), but this will follow when we will prove the completeness of these systems with respect to our intended models
(i.e., topo-dependence models  for $\textbf{LCD}$, and metric dependence models for $\textbf{LUD}$).

\subsection{Completeness wrt topo-dependence models: unravelling}\label{preord-unravel}

In this section we show the decidability and completeness of the logic $\mathsf{LCD}$ wrt topo-dependence models. By results in previous sections, it is enough to prove completeness wrt standard topo-models. For this, we will prove the following

\begin{lem}\label{rep-preord} (\emph{Representation theorem for finite preorder models})
Every finite preorder model for $\mathsf{LCD}$ over a finite set of variables $V$ is a topo-morphic image of a standard topo-model for $\mathsf{LCD}$ (more precisely a standard preorder model), and thus it is modally equivalent to a topo-dependence model.
\end{lem}

The rest of this section is dedicated to the proof of this result.

\medskip

\par\noindent\textbf{The idea of the proof} Throughout the section, we fix a finite preorder model of the form $\bM=(S, =_X, \leq_X, \|\bullet\|)_{X\subseteq V}$ (for $\mathsf{LCD}$), and a designated $s_0\in S$ in this model.\footnote{The specific choice of $s_0$ is irrelevant, as all the other states in $S$ will be accessible from $s_0$ via the total relation $=_\emptyset$.} To find a \emph{standard} preorder model that is modally equivalent to $\bM$, we use the well-known modal logic technique of \emph{unravelling}. The proof goes along familiar lines, by unravelling $\bM$ into a tree of all possible ``histories", redefining the relations on this tree to ensure standardness, and defining a $p$-morphism from the tree into our original model $\bM$, by mapping each history to its last element. The proof is very similar to the corresponding proof for $LFD$ in \cite{BaBe21}, which is itself a modification of the standard modal completeness proof for distributed knowledge (meant to deal with the new dependence atoms $D_XY$).

Admittedly, our construction is \emph{slightly more complicated than necessary for the purposes of this section}, since it will include some redundant numerical components that will play no role in this section, but will be useful later in dealing with the metric case. More precisely,
instead of labeling the transitions in our unravelled tree by using simple \emph{sets} $X$ of variables (`agents') as our labels (as in e.g. the standard proof for $\mathsf{LFD}$), we label them by expressions of the form $X^\beta$, for any set of variables $X\subseteq V$ and any real number $\beta\in [0,1)$. The expression $X^\beta$ is here just a more compact way of writing the ordered pair $(X,\beta)$.
The actual real numbers $\beta$ will play no role in this section, except for the distinction between $\beta=0$ (corresponding to $=_X$-transitions) and $\beta\neq 0$ (corresponding to $\leq_X$-transitions). So in principle we could use only two values for $\beta$ in this section, and the reader can think of all the other real values in our construction as irrelevant `syntactic sugar'; but we do keep all these values for future use in the next section, where they will help us structure our tree of histories as a (standard) pseudo-metric model.\footnote{Even for this purpose, not all real values in $[0,1)$ will actually be needed. If one prefers to keep the resulting model countable, we can restrict to rational values $\beta\in Q\cap [0,1)$. In fact, we can restrict the numbers $\beta$ used in our completeness proofs to any set values $B\subseteq [0,1)$, with the property that $B$ includes $0$, as well as some countable sequence of numbers $(\beta_n)_{n\in N}$ with $0<\beta_n<1$ and $lim_{n\to \infty} \beta_n=0$.}

\medskip

\par\noindent\textbf{Histories} Let
$\bM=(S, =_X, \leq_X, \|\bullet\|)_{X\subseteq V}$ be a finite preorder model for $\mathsf{LCD}$, and let $s_0\in S$ be a designated state in this model.
A \emph{history} is any finite sequence (of any length $n\geq 0$) of the form $h=(s_0, X_1^{\beta_1}, s_1, \ldots, X_n^{\beta_n}, s_n)$, where $s_0$ is the designated world, and for each $k\geq 1$ we have: $s_k\in S$; $X_k\subseteq V$; $\beta_k\in [0,1)$ is a non-negative real number less than $1$; $s_{k-1} \leq_{X_k} s_k$; and finally $s_{k-1} =_{X_k} s_k$ whenever $\beta_k=0$. We denote by $W$ the \emph{set of possible histories}. The map $last:W\to S$, which sends every history $h=(s_0, X_1^{\beta_1}, s_1, \ldots, X_n^{\beta_n}, s_n)$ to its last element  $last(h):=s_n$, will play a key role in our construction.

\medskip

\par\noindent\textbf{One-step transition relations on histories}
There are three relevant one-step relations on histories: the \emph{immediate predecessor relation} $h\to h'$, as well as two \emph{labelled versions} $h\to_Y^= h'$ and $h\to_Y^\leq h'$ of this relation, by putting (for all histories $h,h'$ and sets $Y\subseteq V$):
$$h\to h' \,\, \mbox{  iff } \,\, h'=(h, X^\beta, s') \mbox{ for some } X, \beta,s';$$
$$h\to_Y^{=} h'\,\, \mbox{  iff } \,\, h'=(h, X^0, s') \mbox{ for some } X, s' \mbox{  with } last(h)\models D_X Y;$$
$$h\to_Y^{\leq} h'\,\, \mbox{  iff } \,\, h'=(h, X^\beta, s') \mbox{ for some } X, \beta,s' \mbox{  with } last(h)\models K_X Y.$$
We denote by $\ot$, $\ot_Y^=$ and $\ot_Y^\leq$ the \emph{converses} of these one-step relations.

\medskip

\par\noindent\textbf{Order, $X$-equality and $X$-order on histories} We can now define the \emph{order} relation $h\preceq h'$ on histories, as well as analogues $=_X^W$ and $\leq_x^W$ of the relations $=_X$ and $\leq_X$ on histories in $W$:
$$\preceq \, \,\, \, :=\, \, \,\, \to^*$$
$$=_X^W  \, \,\, :=\, \,\,  (\to_X^=\cup \ot_X^=)^*$$
$$\leq_X^W \, \,\, :=\, \,\, (\to_X^=\cup \ot_X^= \cup \to_X^\leq)^*$$
where we denote by $R^*$ the reflexive-transitive closure of relation $R$.

\medskip

\par\noindent\textbf{Tree of histories, neighbors and (non-redundant) paths}
Two histories $h,h'$ are \emph{neighbors} if either one is the immediate predecessor of the other.
A \emph{path} from a history $h$ to another history $h'$ is a chain of histories $(h_0,\ldots, h_n)$, having $h_0=h$ and $h_n=h'$ as its endpoints, and s.t. for every $k$, histories $h_k$ and $h_{k+1}$ are neighbors. A path is \emph{non-redundant} if no history appears twice in the chain.
The set of histories $(W, \preceq)$ endowed with the partial order $\preceq$ has an easily visualizable \emph{tree structure}:
\begin{itemize}
\item there is a unique history with no predecessor, namely $(s_0)$;
\item every other history $h\not=(s_0)$ has a unique immediate predecessor;
\item for every history $h$, there are only finitely many histories $e\preceq h$;
\item every pair of histories $h,h'\in S$ has a meet $inf(h,h')$ w.r.t. to the order $\preceq$; moreover, by repeating this operation, we can see that triplets of histories also have meets $inf(h,h', h''):=inf(h, inf(h', h''))$,\footnote{More
generally, every non-empty set of histories has a meet. However, the tree structure $(S,\preceq)$ is not a meet semi-lattice, since the empty set of histories $\emptyset$ has no meet: that would be a top element, but the tree has no top.}
\item there exists a unique non-redundant path between any two histories;
\item for all histories $h, h', h''\in W$, we have $inf(h,h', h'')= min \{inf(h,h'), inf(h,h'')\}$, and thus we either have $inf(h,h', h'')=inf(h,h')$ or $inf(h,h',h'')=inf(h,h'')$.
\end{itemize}


\par\noindent\textbf{Interval notation} The \emph{interval} $[h,h']$ between two histories $h, h'$ is the set of all histories on the non-redundant path from $h$ to $h'$, i.e.:
$$[h,h']\,\, :=\,\, \{h'' \in S: h'' \mbox{ lies on some non-redundant path from } h \mbox{ to } h'\}.$$
It is easy to see that in general we have for all histories $h, h', h''\in W$:
$$[h,h']\subseteq [h,h'']\cup [h'', h'];$$
and this inclusion becomes equality whenever $h''$ lies between $h$ and $h'$ (i.e., in the interval between them):
$$h'' \in [h, h'] \mbox{ implies } [h,h'] = [h,h'']\cup [h'', h'].$$

\medskip

\par\noindent\textbf{Tree-distance and history length}
The \emph{tree-distance} $dist(h,h')$ between histories $h$ and $h'$ is defined
by putting $dist(h,h'):=n-1$, where $n$ is the cardinality of the set $[h,h']$ (i.e., the length of the non-redundant path $(h_0,\ldots, h_n)$ from $h$ to $h'$). It is easy to see that $dist$ is a metric on $S$, and that moreover $dist(h,h')=1$ iff $h, h'$ are neighbors. The \emph{length} $l(h)$ of a history $h=(s_0, X_1^{\beta_1}, s_1, \ldots, X_n^{\beta_n}, s_n)$ is the tree-distance $l(h):=n$ between $h$ and the root $(s_0)$.

\medskip

\par\noindent\textbf{Path-characterizations of $X$-equality and $X$-order} It is easy to see that, given histories $h, h'$, if
 $(h_0, \ldots,h_i, \ldots h_n)$ is the \emph{non-redundant path} from $h$ to $h'$ where $h_i=inf(h,h')$ is the\emph{ meet} of the two histories, then the $X$-equality and $X$-order relations can be characterized as follows:

\begin{itemize}
\item $h=_X^W h'$ holds iff $h=h_0\ot_X^=h_1 \ldots \ot_X^= h_i \to_X^= \ldots h_{n-1} \to_X^= h_n= h'$;
\item $h\leq_X^W h'$ holds iff $h=h_0\ot_X^=h_1  \ot_X^= h_i (\to_X^=\cup \to_X^\leq)  \ldots h_{n-1} (\to_X^= \cup \to_X^\leq) h_n = h'$.
\end{itemize}

To check this, it is enough to note that the definitions of $=X^W$ and $\leq_X^W$ imply the existence of appropriate paths from $h$ to $h'$ while that the tree structure implies that any such path must include the non-redundant path $(h_0, \ldots,h'', \ldots h_n)$, then finally to check that this non-redundant path is ``appropriate" for  $=X^W$ and respectively for $\leq_X^W$ only if it satisfies the corresponding clause in the above characterization.\footnote{Technically speaking, proving this requires induction on the tree-distance $d(h,h')$.}

\medskip

\begin{quote}
\textbf{Claim 1}. The relations
$=_Y^W$ and $\leq_Y^W$ have the following properties:
\begin{enumerate}
\item $=_Y^W$ are equivalence relations, and $\leq_Y^W$ are preorders;


\item if $h=_Y^W h'$, then $h\leq_Y^W h'$;
\item the intersection of $\leq_Y^W$ with its converse $\geq_Y^W$ is the relation $=_Y^W$;


\item the relations $=_Y^W$ and $\leq_Y^W$ are additive:
$$h=_Y^W h' \mbox{ iff } h =_y h' \mbox{ for all $y\in Y$};$$
$$h\leq_Y^W h' \mbox{ iff } h \leq_y h' \mbox{ for all $y\in Y$};$$


\item if $h=_Y^W h'$ then $last(h)=_Y last(h')$;


\item if $h\leq_Y^W h'$ then $last(h)\leq_Y^W last(h')$;


\item if $h =_X^W h'$ and $last(h)\models D_XY$, then $h =_Y^W h'$ and $last(h')\models D_XY$;


\item if $h \leq_X^W h'$ and $last(h)\models K_XY$, then $h \leq_Y^W h'$ and $last(h')\models K_XY$;

\item if $h=_Y^W h'$ and $h''\in [h,h']$, then $h=_Y^W h'=_Y^W h''$;
\item if $h\leq_Y^W h'$ and $h''\in [h,h']$, then $h\leq_Y^W h'' \leq_Y^W h'$.
\end{enumerate}
\end{quote}

\begin{proof}
Items 1 and 2 follow immediately from the definitions of these relations on histories.

Item 3: The fact that $h\leq_Y^W h'\leq_Y^W h$ implies $h=_Y^W h'$ is easily proved by induction on the tree-distance $d(h,h')$ between the two histories, using the above characterizations of $=_Y^W$ and $\leq_Y^W$ in terms of the non-redundant path, as well as the fact that $\to_Y^\leq$ and $\ot_Y^{\leq}$ are disjoint by definition.

For item 4, we first check the additivity of $\to_Y^=$ and $\to_Y^\leq$, i.e. that: $h\to^=_Y h'$ iff $h \to^=_y h'$ for all $y\in Y$; and $h\to^\leq_Y h'$ iff $h \to^\leq_y h'$ for all $y\in Y$. (These follow directly from the definitions of $\to^=_Y$ and $\to^\leq_Y$ and the Additivity axioms for $D_XY$ and $K_XY$.) The desired additivity of $=_Y^W$ and $\leq_Y^W$ follows from these observations by induction of the tree distance $d(h,h')$ (using the above characterizations of  $=_Y^W$ and $\leq_Y^W$ in terms of the non-redundant path).

Items 5 and 6 similarly follow by easy induction (on the tree distance, using the above path characterization) from the corresponding statements for the one-step transitions, i.e.:  $h\to^=_Y h'$ implies $last(h)=_Y last(h')$; and similarly $h\to^\leq_Y h'$ implies $last(h)\leq_Y last(h')$. To prove the first of these one-step statements: from $h \to^=_Y h'$, we infer that $h'=(h, X^0, s')$ for some $X, s'$ with $last(h)\models D_XY$; but in order for $h'$ to be a well-formed history, we must have $last(h)=_X s'$, which together with $last(h)\models D_XY$ gives us $last(h)=_Y s'=last(h')$, as desired. The second one-step statement is similar: from $h \leq^=_Y h'$, we infer that $h'=(h, X^\beta, s')$ for some $X, s'$ with $last(h)\models K_XY$; but in order for $h'$ to be a well-formed history, we must have $last(h)\leq_X s'$, which together with $last(h)\models K_XY$ gives us $last(h)\leq_Y s'=last(h')$, as desired.

Similarly, items 7 and 8 follow by easy induction (on the tree distance) from the corresponding statements for the one-step transitions. To check the first of these, assume $h \to^=_X h'$ with $last(h)\models D_XY$. By definition, the first of these means that we have  $h'=(h, Z^0, s')$ for some $Z, s'$ with $last(h)\models D_ZX$. Together with $last(h)\models D_XY$, this gives us $last(h)\models D_ZY$ (by the Transitivity axiom for $D$), and hence we have $h\to^=_Y$, as desired. The case of the relatiion $\to^\leq_X$ (for item 8) is completely analogous.

Finally, items 9 and 10 follow by applying the above characterizations of $h =_Y^W h'$ and respectively $h\leq_Y^h h'$ in terms of the non-redundant path from $h$ to $h'$, and then noticing that $h''\in [h, h']$ implies that the same characterizing conditions apply to the subpaths from $h$ to $h''$, and respectively from $h''$ to $h'$.
\end{proof}

\medskip

\par\noindent\textbf{The standard preorder model on histories} To structure the set $W$ of histories as a standard preorder model (in simplified presentation)
$$\bM^\infty_\leq = (W, \leq_x^W, \|\bullet\|),$$
we take as our preorders $\leq_x^W$ to be just the relations $\leq_{\{x\}}^W$ defined above, while the valuation for atoms of the form $P\vec{x}$ is inherited from our finite preorder model $\bM$:
$$\| Px_1\ldots x_n\| = \{h\in W: last(h)\models Px_1\ldots x_n\}.$$
This is meant to be a \emph{standard} preorder model, so technically speaking the relations $\leq_X^W$ and $=_X^W$ will be defined in this model by standardness (as intersections, in terms of the basic relations $\leq_x^W$, using clauses (1) and (2) of Proposition \ref{Standard3}), rather than the way we defined them above on histories. But Claim 1 ensures that they match the above definitions, as made explicit in the following result:

\begin{quote}
\textbf{Claim 2}. The structure $\bM^\infty_\leq$ is indeed a standard preorder model, with its induced $X$-equality and $X$-preorder relations matching the relations $\leq_X^W$ and $=_X^W$, as defined above on histories.
\end{quote}

\begin{proof} The fact that the relations $\leq_X^W$ (as defined above on histories) coincide with the intersections $\bigcap_{x\in X}\leq_x^W$  follows from the additivity of $\leq_X^W$ (item 4 of Claim 1), while the fact that the relations $=_X^W$ coincide with the intersections $\leq_X^W\cap \geq_X^W$ follows from item 3 of Claim 1.
So, to prove standardness (given the simplified presentation), it is enough to check condition ($0$) of Proposition \ref{Standard3}. For this, assume that we have $h=_X^W h'$ and $h\in \|Px_1\ldots x_n\|$. By our choice of valuation on $\bM^\infty_\leq$, this means that $last(h)\models Px_1\ldots x_n$. By item 5 in Claim 1,  $h=_X^W h'$ implies that $last(h)=_X last(h')$ in $\bM$. This, together with $last(h)\models Px_1\ldots x_n$ (and the fact that $\bM$ is a preorder model), give us $last(h')\models Px_1\ldots x_n$ (by clause 1 in the definition of preorder models), from which we can conclude that $h'\in \|Px_1\ldots x_n\|$ (again by our choice of valuation on $\bM^\infty_\leq$), as desired.
\end{proof}

We now finish  the proof of our representation result (Lemma \ref{rep-preord} above), with the following:

\begin{quote}
\textbf{Claim 3}.
The function $last:S \to W$, that maps every history $h\in S$ to its last state $last(h)\in W$, is a surjective $p$-morphism from the standard preorder model $\bM^\infty_\leq$ to the given finite preorder model $\bM$.
\end{quote}

\begin{proof} The \emph{surjectivity} of the map $last$ is immediate: given $s\in W$, if we take the history $h=(s_0, \emptyset^0, s)$ of length $1$ (which is a correct history, since $s_0=_\emptyset s$ by the definition of preorder models), then $last(h)=s$. So we just need to check the clauses in the definition of $p$-morphisms on preorder models.

\emph{Atomic clause $Px_1\ldots x_n$}: holds by our choice of valuation in $\bM^\infty_\leq$.

\emph{Atomic clause $D_XY$}: We need to show that we have $h\models D_XY$ iff $last(h)\models D_XY$.

For the \emph{left-to-right implication}, suppose that $h\models D_XY$, and let $s:=last(h)$;
we must prove that $s\models D_XY$. For this, take the history
$h':=(h, X^0, s)$. Then we have $h=_X^W h'$, which together with $h\models K_XY$, gives $h=_Y^W h'$ (since $\bM^\infty_\leq$ is a preorder model). Since $=_Y^W$ is the reflexive-transitive closure of $\to_Y^=\cup \ot_Y^=$ while $h$ is an immediate predecessor of $h'=(h, X^0, s)$,  $h=_Y^W h'$ can only hold in this case iff we have $h\to_Y^= h'$. By the definition of $\to_Y^=$ and the structure of $h'=(h, X^0, s)$, this means that we must have $s=last(h')\models D_XY$, as desired.

For the \emph{converse implication}, suppose that $last(h)\models D_X Y$; we need to show that $h\models D_XY$ as well. For this, let $h'\in W$ be any  history in $W$ with $h=_X^W h'$, and we need to prove that $h=_Y^W h'$. Indeed, from $h=_X^W h'$ and $last(h)\models K_X Y$, we get $h=_Y^W h'$ (by item 7 of Claim 1), as desired.

\emph{Atomic clause $K_XY$}: We need to show that we have $h\models K_XY$ iff $last(h)\models K_XY$.

The \emph{left-to-right implication} is similar to the above proof for $D_XY$: suppose that $h\models K_XY$, and let $s:=last(h)$;
we need to prove that $s\models K_XY$. For this, take again the history
$h':=(h, X^0, s)$. Then we have $h\leq_X^W h'$ (in fact $h=_X^W h'$), which together with $h\models K_XY$, gives us $h\leq_Y^W h'$ (since $\bM^\infty_\leq$ is a preorder model). Once again, since $\leq_Y^W$ is  the reflexive-transitive closure of $\to_Y^\leq$ and $h$ is an immediate predecessor of $h'=(h, X^0, s)$,  $h\leq_Y^W h'$ can only hold in this case iff we have $h\to_Y^\leq h'$. By the definition of $\to_Y^\leq$ and the structure of $h'=(h, X^0, s)$, this means that we must have $s=last(h')\models K_XY$, as desired.

For the \emph{converse implication}, let $last(h)\models K_X Y$; we need to show that $h\models K_XY$ as well. For this, let $h',h''\in S$ be any  histories with $h\leq_X^W h'\leq_X^W h''$, and we must prove that $h'\leq_Y^W h''$. From $h\leq_X^W h'$ and $last(h)\models K_X Y$, we get $last(h')\models K_XY$ (by item 8 of Claim 1 above).
Using this and $h'\leq_X^W h''$, we obtain $h'\leq_Y^W h''$ (by applying again item 8 of Claim 1), as desired.

The \emph{forth clauses} for $=_X^W$ and $\leq_X^W$ are given by parts 5 and 6 of Claim 1 above.

The \emph{back clause} for $=_X^W$:
Assuming $last(h)=_X s$, we take $h':=(h,X^0,s)$. This is a well-formed history in $W$, satisfying $h=_X^W h'$ and $last(h')=s$, as desired.

The \emph{back clause} for $\leq_X^W$:
Assuming $last(h)\leq_X s$, we take $h':=(h, X^\beta, s)$ for any  number $\beta\in (0,1)$ chosen in our construction. Once again, this is a well-defined history in $W$,
obviously satisfying the conditions $h\leq_X^W h'$ and $last(h')=s$, as required.
\end{proof}

\medskip

Putting together Proposition \ref{PreFMP} and Lemma \ref{rep-preord} with the well-known facts about the preservation of modal formulas under $p$-morphisms (or more generally, under bisimulations), as well as the modal equivalence between standard preorder models and Alexandroff topo-dependence models, we immediately obtain the completeness of $\textbf{LCD}$ wrt topo-dependence models (and as a consequence the decidability of the logic $\mathsf{LCD}$ of topo-dependence models).

\subsection{Completeness wrt metric dependence models: refinement}

We now wrap up our completeness (and decidability) proofs for $\mathsf{LCD}$ and $\mathsf{LUD}$ wrt (pseudo-locally Lipschitz) metric dependence models. By previous results, it is enough to prove completeness wrt standard pseudo-metric models (whose pseudo-metrics satisfy the pseudo-locally Lipschitz condition). For this, we show the following:

\begin{lem}\label{rep-metric} (\emph{Pseudo-metric representation theorem})
Every finite preorder model for $\mathsf{LCD}$ and respectively $\mathsf{LUD}$ (over a finite set of variables $V$) is a topo-morphic image of a standard pseudo-metric model for $\mathsf{LCD}$, and respectively $\mathsf{LUD}$. Moreover, the metric dependence model associated to this standard pseudo-metric models is pseudo-local Lipschitz.
\end{lem}

\medskip\par\noindent\textbf{Plan of the proof}.
Essentially, the proof proceeds by (a) defining pseudo-distances $d_X^W: W\times W\to [0,1]$ between histories in the unravelled tree model $\bM^\infty_\leq=(W, =_X^W, \leq_X^W, \|\bullet\|)$ introduced in the previous section; then (b) showing that, importantly,  this pseudo-metric topology is a \emph{refinement of the Alexandroff topology} associated to the preorders $\leq_X^W$ from the previous section, and (c) proving that the resulting structure $\bM^\infty_d=(W, =_X^W, d_X^W, \|\bullet\|)$ is a \emph{standard pseudo-metric model}; and finally (d) showing that the map $last:W \to S$ (sending each history to its last element) is a (surjective) topo-morphism from the $\mathsf{LUD}$ topo-model $\bM^\infty_d$ to (the $\mathsf{LUD}$ topo-model associated) to our original preorder model for $\mathsf{LUD}$.

\medskip\par\noindent\textbf{The proof in detail} In the rest of this section, we present the details of our proof of Lemma \ref{rep-metric}. We give the proof of this representation result for $\mathsf{LUD}$ only; the proof for $\mathsf{LCD}$ can be obtained by simply skipping all the steps and clauses involving the uniform dependence atoms $U(X;Y)$. We start as in the previous section: let $\bM=(S, =_X, \leq_X, \|\bullet\|)_{X\subseteq V}$ be a finite preorder model for $\mathsf{LUD}$. Since this is also in particular a finite preorder model for $\mathsf{LCD}$, we can apply the unravelling technique from the previous section to obtain the structure $\bM^\infty_\leq=(W, =_X^W, \leq_X^W, \|\bullet\|)$, as introduced in the previous section. All the notations and results in the previous section still apply, and we will make use of them, in particular of the fact that $\bM^\infty_\leq$ is a standard preorder model for $\mathsf{LCD}$, and that the map $last:W \to S$ (sending each history to its last element) is a surjective $p$-morphism between preorder models for $\mathsf{LCD}$, and hence a (surjective) interior map between the corresponding topo-models for $\mathsf{LCD}$. For $h\in W, s\in S$, we will use the notations
$$h\uparrow_X^W :=\{h'\in W: h\leq_X^W h'\} \,\,\,\,\,\,\, \mbox{ and }
\,\,\,\,\,\,\,  s\uparrow_X:=\{s'\in S:s\leq_X s'\}$$
to denote the corresponding $X$-upsets in $W$ and respectively $S$. Note that the fact that the surjective map $last: \bM^\infty_\leq\to \bM$ is a $p$-morphism is equivalent to saying that, for all histories $h\in W$, we have:
$$last(h)\uparrow_X= last(h\uparrow_X^W).$$

As announced in the Proof Plan, we will endow the set $W$ of histories with a pseudo-metric structure, obtained by defining \emph{ultra-pseudo-metric distances} $d_X^W: W\times W\to [0,1]$, and at the same time we will \emph{extend  the valuation} to cover atoms of the form $U(X;Y)$, to obtain a pseudo-metric refinement $\bM^\infty_d=(W, =_X^W, d_X^W, \|\bullet\|)$ of the topology of our unravelled preorder model  $\bM^\infty_\leq$. To do this, we need to introduce some preliminary notions.

\medskip

\par\noindent\textbf{The $X$-root of a history} When studying the Alexandroff topologies induced by the preorders $\leq_X$ on the set $W$ of all histories, the following notion is of special interest: the \emph{$X$-root} of a history $h$ is the shortest history $h_X$ s.t. $h_X\leq_X^W h$.\footnote{Here, ``shortest" means that there is no proper sub-history $e\prec h_X$ with $e\leq_X^W h$.} It is easy to see that \emph{every history has a unique $X$-root}, which we will henceforth denote by $h_X$; moreover, we obviously have $h_X\preceq h$.
In fact, if we denote by $h(X)$ the $\leq_X$-connected component\footnote{A set $A\subseteq W$ of histories is $\leq_X$-connected if for every two histories $h,h'\in A$ there exists some history $h''\in A$ with $h''\leq_X h, h'$. The largest $\leq_X$-connected subset of $W$ that contains a history $h$ is called the $\leq_X$-connected component of $h$, and denoted by $h(X)$. In fact, $h(X)$ is the equivalence class of $h$ with respect to the equivalence relation $\asymp_X$ given by putting $h\asymp_X h '$ iff there exists $h''\leq_X h, h'$.} that contains $h$, it is easy to see that $h(X)$ is a subtree of $W$ (wrt $\preceq$), having $h_X$ as its root. 
Obviously, states that belong to the same $\leq_X$-connected component have the same $X$-root, i.e. we have that:
$$h\leq_X^W h' \mbox{ implies } h_X=h'_X.$$

\medskip

\par\noindent\textbf{Density and $X$-density of a history}
In addition to history length, we need another measure of the complexity of a history $h$.
The \emph{density} $\delta^h$ of a history $h$ is defined as the \emph{minimum non-zero index $\beta$ occurring in the history $h$} (where we put $\delta^h=1$ when there are no such non-zero indices in $h$).
More precisely, we put
$$\delta^{h}\, \, :=\,\, inf  \{\beta>0: (e, Y^\beta, s)\preceq h \mbox{ for some } Y\subseteq V, s\in S\},$$
with the convention that $inf \,\emptyset =1$. Note that, in our setting, this convention is the natural analogue of the more standard convention that $inf \,\emptyset=\, \infty$: since $d_X^W$ takes values in $[0,1]$, the maximum possible distance is $1$, playing the role of `infinity'. A distance $d_X^W(h,h')=1$ will indicate that the two histories are so ``far" from each other, that one cannot reach $h'$ from $h$ by any number of ``small" changes of $X$-value. Also note that the infimum in the above definition is applied to a finite set of real values, so it can be replaced by minimum whenever the set is non-empty.

The density of a history $h$ is a measure of how ``far" apart are its preceding histories.
However, for defining our $X$-topology on the tree $W$, not all the past of a history $h$ is relevant, but only the past of its $X$-root.
The \emph{$X$-density} of history $h$ is simply the \emph{density of its $X$-root}, i.e., the quantity
$$\delta^{h_X}.$$


\par\noindent\textbf{Closeness} To define our $X$-distances, we need a notion of closeness between histories that differs from the tree distance: indeed, even a one-step transition may involve a large jump as far as the intended $X$-distance is concerned.
Two \emph{neighboring} histories $h=(s_0, \ldots, s)$ and $h'=(h, Y^\beta, s')$ are said to be \emph{$X$-close} if any of the following conditions holds:
$$\mbox{ either $\bM\models U(Y;X)$; or $s\models_\bM K_YX$ and $\beta<\delta^{h_X}$; or $s\models_\bM D_YX$ and $\beta=0$.}$$
For \emph{arbitrary} histories, $h$ and $h'$ are said to be \emph{$X$-close} if all the neighboring histories on the non-redundant path from $h$ to $h'$ are $X$-close. Two histories are said to be \emph{$X$-far} if they are \emph{not} $X$-close. Note that, since $U(X;\emptyset)$ is a valid formula, \emph{every two histories are $\emptyset$-close}; in other words: \emph{if two histories are $X$-far then $X\neq\emptyset$}. Finally, in the case that $X=\{x\}$ is a singleton, we can drop the set brackets, saying that two histories are $x$-close ($x$-far) if they are $\{x\}$-close ($\{x\}$-far).

\medskip

\begin{quote}
\textbf{Claim 1}. For all histories $h=(s_0, \ldots, s)$ and all sets $X\subseteq V$, we have:
\begin{enumerate}
\item $0<\delta^h \leq 1$;
\item if $h\preceq h'$, then $\delta^h\geq \delta^{h'}$;
\item $\delta^h\leq \delta^{h_X}$;
\item if $h\leq_X^W h'$, then $\delta^{h_X}=\delta^{h'_X}$;
\item if $\bM\models K(Y;X)$ (and so also, in particular, if $\bM\models U(Y;X)$), then $\delta^{h_Y}\leq \delta^{h_X}$;
\item $\delta^{h_{X\cup Y}} = min \{\delta^{h_X}, \delta^{h_Y}\}$;
\item $\delta^{h_{X}}=  max \{\delta^e: e\leq_X^W h  \}$;
\item $X$-closeness is additive:  two histories are $X\cup Y$-close iff they are both $X$-close and $Y$-close;
\item $X$-closeness is an equivalence relation.
\end{enumerate}
\end{quote}

\begin{proof}: Parts 1 and 2 follow directly from the definitions of $\delta^h$. Part 3 follows from part 1 and the fact that $h_X\preceq h$. Part 4 follows from the fact that $h\leq_X^W h'$ implies $h_X=h'_X$.

Part 5 follows from the observation that $\bM\models K(Y;X)$ implies that every history $e\leq_Y^W h$ has also the property that $e\leq_X^W h$, hence in particular we have that $h_Y\leq_X h$. By the definition of $h_X$, this implies that $h_X\preceq h_Y$.

For part 6: since $h_X, h_Y\preceq h$, by the tree property we must have either $h_X\preceq h_Y$ or $h_Y\preceq h_X$. Our claim is symmetric in $X$ and $Y$, so without loss of generality we can assume that $h_X\preceq h_Y$, hence by part 2 we have $\delta^{h_X} \geq \delta^{h_Y}$, i.e.
$min\{\delta^{h_X}, \delta^{h_Y}\}= \delta^{h_Y}$. On the other hand, since $h_X\leq_X^W h$ and $h_X\preceq h_Y \preceq h$, by Claim 1(10) in Section \ref{preord-unravel} we must also have $h_Y\leq_X^W h$; putting this together with $h_Y\leq_Y^W h$ and using the additivity of the relation $\leq^W$, we obtain that $h_Y\leq_{X\cup Y} h$; moreover, by definition of $h_Y$, no strict predecessor of $h_Y$ is $\leq_Y^W h$, hence no such predecessor can be $\leq_{X\cup Y}^W h$; this, together with the fact that $h_Y\leq_{X\cup Y} h$, shows that $h_Y$ satisfies the defining property of $h_{X\cup Y}$, so we have $h_Y =h_{X\cup Y}$. In conclusion, we obtain $min\{\delta^{h_X}, \delta^{h_Y}\}= \delta^{h_Y}= \delta^{h_{X\cup Y}}$, as desired.

For part 7, we prove two inequalities separately. The fact that $\delta^{h_X}\leq max \{\delta^e: e\leq_X^W h\}$ follows immediately from the fact that $h_X\leq_X^W h$. To check the converse inequality $\delta^{h_{X}}\geq  max \{\delta^e: e\leq_X^W h  \}$, let $e$ be any history s.t. $ e\leq_X^W h$. Then by part 4 we have $\delta^{e_X}=\delta^{h_X}$, and so by part 3 we have $\delta^e\leq \delta^{e_X}=\delta^{h_X}$, as desired.

Part 8 follows directly from the definition of $X$-closeness, using the Additivity Axioms for $U(X;Y)$, $K_XY$ and $D_XY$, as well as part 6.

For part 9, it is obvious that $X$-closeness is reflexive and symmetric. For transitivity, let $h, h', h''\in W$ be s.t. $h$ is $X$-close to $h'$, and $h'$ is $X$-close to $h''$. By definition, this means that all neighboring histories in $[h,h']$ are $X$-close, and that the same holds for all neighboring histories in $[h', h'']$. Using the fact that $[h, h'']\subseteq [h,h']\cup [h',h'']$, it follows that all neighboring histories in $[h, h'']$ are also $X$-close, and thus by definition $h$ and $h''$ are $X$-close, as desired.
\end{proof}

\medskip

\par\noindent\textbf{The pseudo-metric model on histories}
We proceed now to define our ultra-pseudo-metric model (in simplified presentation) $\bM^\infty_d=(W, d_x^W, \|\bullet\|)$. As the \emph{set of states}, we take again the set $W$ of all histories, endowed with the same valuation on atoms of the form $P\vec{x}$ as on $\bM^\infty_\leq$:
$$\|Px_1\ldots x_n\|=\{h\in W: last(h)\models Px_1\ldots x_n\}.$$
Since $\bM^\infty_d$ is meant to be a \emph{standard} pseudo-metric model, the only thing left is to \emph{define the basic pseudo-metrics} $d_x^W:W\times W\to [0,\infty)$. Essentially, the pseudo-distance between two `close' histories will be given by the maximum of the real numbers $\beta\in [0,1)$ encountered on the non-redundant path between them; while the pseudo-distance between `far' histories will be $=1$ by definition. More precisely, for every $x\in V$ we put:

\vspace{-3ex}

\begin{eqnarray*}
d_x^W(h,h') &:=& max\{\beta\in [0,1): \exists e, (e, Y^\beta, s)\in [h, h'] \mbox{ with } s\in S,Y\subseteq V\}, \mbox{ if $h$ and $h'$ are $x$-close}; \\
d_x^W(h,h') & := & 1, \mbox{ if $h$ and $h'$ are $x$-far},
\end{eqnarray*}

\vspace{-1ex}

where we assume the convention that $max \,\emptyset=0$, thus ensuring that $d_x^W(h,h)=0$.

For \emph{sets} of variables $X\subseteq V$, the pseudo-distances $d_X^W$ are defined like in any standard pseudo-metric model, using the Chebyshev distance
$$d_X^W(h,h') \, :=\, max\{d_x^W(h,h'): x\in X\}.\footnote{Once again, the convention that $max\, \emptyset=0$ ensures that $d_\emptyset^W(h,h')=0$ for all $h,h'\in W$, thus matching the fact that $=_\emptyset$ is the universal relation.}$$
As usual, for every history $h\in W$, set $X\subseteq V$ and number $\varepsilon\in (0,1]$, we use the notation
$$\mathcal{B}_X(h, \varepsilon):=\{h'\in W: d_X(h,h')<\varepsilon\}$$
to denote the ``open ball of radius $\varepsilon$ centered at $h$".

\smallskip

We now proceed to prove a number of useful results. First, using the additivity of $X$-closeness (part 8 of Claim 1), we can immediately extend the above defining equations for $d_x^W$ to \emph{sets} $X$ of variables (given the definition of $d_X^W$ in terms of Chebyshev distance):


\begin{quote}
\textbf{Claim 2}. For every set of variables $X\subseteq V$, the functions $d_X^W$ satisfy the conditions:
\begin{eqnarray*}
d_X^W(h,h') & = & 0, \,\, \,\, \mbox{ if $X=\emptyset$};\\
d_X^W(h,h') &=& max\{\beta\in [0,1):  \exists e, (e, Y^\beta, s)\in [h, h'] \mbox{ with } s\in S,Y\subseteq V\},\\
& &   \,\,\,\,\,\,\,\,\,\, \mbox{ if $h,h'$ are $X$-close with $X\neq \emptyset$};\\
 d_X^W(h,h') \, & = &  1, \,\,\,\,\,\, \mbox{if $h$ and $h'$ are $X$-far (hence $X\neq\emptyset$)}.
\end{eqnarray*}
\end{quote}

As an immediate consequence of Claim 2, we have the following

\medskip

\par\noindent\textbf{Important Observations} For all histories $h, h'\in W$, we have:
\begin{enumerate}
\item $d_X^W(h, h') <1$ iff $h$ and $h'$ are $X$-close;
\item if $d_X^W(h, h'), d_Y^W(h,h') <1$ with $X, Y\neq \emptyset$, then $d_X^W(h,h')=d_Y^W(h,h')$.
\end{enumerate}

\medskip

Next, we look at distances between neighboring histories, spelling out a few immediate consequences of Claim 2:

\smallskip

\begin{quote}
\textbf{Claim 3}.  For all neighboring histories $h=(s_0, \ldots, s)$, and $h'=(h,Y^\beta, s')$, and all non-empty sets $X,Z\subseteq V$ (with $X,Z\neq \emptyset$), we have the following:
\begin{enumerate}
\item if $d_X^W(h,h') =0$, then $\beta=0$, $d_Z^W(h,h')\in \{0,1\}$, $s\models D_YX$ and $s=_Ys'$ (hence also $s=_Xs'$);
\item if $X\subseteq  Y$, then $d_X^W(h, h')= \beta$;
\item if $d_X^W(h,h')<1$, then $d_X^W(h,h')=\beta<1$, $s\models D_YX$ and $s\leq_X s'$;
\end{enumerate}
\end{quote}

\begin{proof} All parts follow from Claim 2 by simple case inspection. For part 1, we also use the constraints on histories and the fact that $S^0_X$ is the same as $X$, as well as the fact that, by Claim 2,  $\beta=0$ is  consistent only with $d_Z(h,h')\in \{0,1\}$). For part 2, we use the fact that $X\subseteq Y$ implies $\bM \models U(Y;X)$, which implies $s\models D_YX$. For part 3, we use the fact that both $U(Y;X)$ and $K_YX$ imply $D_YX$, and the constraints on histories that imply that $s\leq_Y^W s'$,
as well as the fact that $d_X^W(h,h')<1$ implies that $d_X^W(h,h')$ and that we have either $s\models K_YX$ (which together with $s\leq_Y^W s'$ gives us $s\leq_X s'$) or $\beta=0$ (in which case the constraints on histories give us that $s=_X s'$, hence also $s\leq_X s'$).
\end{proof}

\medskip

We now look at the key properties of our distances $d_X^W$ between arbitrary histories:

\smallskip

\begin{quote}
\textbf{Claim 4}.
Each $d_x^W$
satisfies the following properties:
\begin{enumerate}
\item $d_X^W(h,h)=0$;
\item $d_X^W(h,h')=d_X^W(h',h)$;
\item $d_X^W(h,h')=0$ holds iff $h=_X^W h'$;
\item if $h''\in [h,h']$, then $d_X^W(h,h') = max\{d_X^W(h,h'') , d_X^W(h'', h')\}$;
\item if $e,e'\in [h,h']$, then $d_X^W(e,e')\leq d_X^W(h,h')$;
\item $d_X^W(h, h')\leq max\{d_X^W(h,h'') , d_X^W(h'', h')\}$ (for \emph{arbitrary} $h''$);
\item if $h\preceq h'$ and $d_X^W(h,h')< 1$, then $h\leq_X^W h'$ (and hence $last(h)\leq_X last(h')$);
\item if $d_X^W(h,h')< 1$, then $inf(h,h')\leq_X^W h, h'$ (and hence $\delta_X^h=\delta_X^{h'}=\delta_X^{inf(h,h')}$);
\item $\mathcal{B}_X(h,\delta^h)\subseteq h\uparrow_X$;
\item $last(h)\uparrow_X \subseteq last (\mathcal{B}_X(h,\varepsilon))$ for all $\varepsilon>0$;
\item if $\bM\models U(Y;X)$ and $d_Y^W(h,h')<1$, then $h$ and $h'$ are $X$-close;
\item if $last(h)\models_\bM K_YX$ and $d_Y^W(h,h')< \delta^h$, then $h$ and $h'$ are $X$-close.
\end{enumerate}
\end{quote}

\begin{proof} Part 1 follows from Claim 2
and the convention that $max\, \emptyset =0$. Part 2 follows immediately from the symmetrical shape of the inductive definition of $d_X(h,h')$, and the in-built symmetry of $d_X(h,h')$ in the case that $dist(h,h')=1$.

For part 3: The special case in which $h$ and $h'$ are neighboring histories, say $h'=(h,Y^\beta, s')$ should be obvious, using Claim 3(1), via the following sequence of equivalences:
$d_X^W(h, h')=0$ iff $(\beta=0 \wedge last(h)\models D_YX)$
iff $(h'=(h, Y^0, \beta)\wedge last(h)\models D_YX)$
iff $h\to_X^= h'$ iff $h=_X^W h'$.
The general case follows from this special case by induction on the length of the tree distance  between $h$ and $h'$ (using the general definition of $d_X^h$ for non-neighboring histories, as well as the transitivity of $=_X^W$).

Part 4 follows immediately from the definition of $d_X^W$, together with the above-mentioned fact that:
$$h'' \in [h, h'] \mbox{ implies } [h,h'] = [h,h'']\cup [h'', h'].$$
Part 5 is another immediate consequence of the definition of $d_X^W$, while part 6 follows from the same definition together with another fact mentioned above, namely that the inclusion
$$[h,h'] \subseteq [h,h'']\cup [h'', h']$$
holds for arbitrary histories $h''$.

For part 7, assume $h\preceq h'$ and $d_X^W(h,h')<1$. Then, for all pairs $e,e'\in [h,h']$ satisfying $e\to e'$, we have $d_X^W(e,e')\leq d_X^W(h,h')<1$ (by  part 6), hence by part 3 of Claim 3 we have $last(e)\leq_X last(e')$. Putting all these together along the non-redundant path from $h$ and $h'$ and applying the transitivity of $\leq_X^W$, we obtain the desired conclusion.

To prove part 8, assume $d_X^W(h,h')<1$ and put $h'':=inf(h,h')$. Since $h''\in [h,h']$, by part 5 we have  $d_X^W(h'', h), d_X^W(h'', h')\leq d_X^W(h,h')<1$. This together with the fact that $h''=inf(h,h')\preceq h, h'$, gives by part 7 that $h''\leq_X^W h,h'$, as desired.

For part 9, assume that $h\in \mathcal{B}_X(h,\delta^h)$. From this and the fact that $\delta^h\leq 1$, we infer that $d_X^W(h, h')<1$, so if we put  $h'':=inf(h,h')$, then by part 8 we have $h''\leq_X^W h'$.
It is clearly enough to prove that $h''=_X^W h$ (hence $h\leq_X^W h''\leq_X^W h'$, and thus $h'\in h\uparrow_X^W$, as desired), which by  part 3 is equivalent to showing that $d_X^W(h'', h)=0$. For this, we need to prove that $d_X^W(e,e')=0$ holds for all pairs $e,e'\in [h'', h]$ with $e \to e'$. To prove this, assume towards a contradiction that we have $d_X^W(e,e')\neq 0$ for some such pair. Then, by the definition of $\delta^h$ and the fact $e,e'\preceq h$ (since $h''\preceq h$ and $e,e'\in [h'', h]$), we have that $\delta^h\leq d_X^W(e,e')$. But on the other hand, from $e,e'\in [h'',h]$ and $h''\preceq h, h'$ we obtain by  part 5 that $d_X^W(e,e')\leq d_X^W(h,h'')\leq d_X^W(h,h')<\delta^h$, which contradicts $\delta^h\leq d_X^W(e,e')$.

For part 10, assume that $\varepsilon>0$ and $s\in last(h)\uparrow_X$, i.e. $last(h)\leq_X s$. Choose some $\beta\in (0,\varepsilon)$, and take the history $h':=(h,X^\beta, s)$. The fact that $last(h)\leq_X s$ ensures that $h'$ is a well-defined history with $last(h')=s$. By part 2 of Claim 3, we also have that $d_X^W(h,h')=\beta<\varepsilon$, i.e. $h'\in \mathcal{B}_X(h,\varepsilon)$, and thus $s=last(h')\in last(\mathcal{B}_X(h,\varepsilon))$, as desired.

For part 11, suppose that $\bM\models U(Y;X)$, and let $h,h'\in W$ be histories with $d_Y(h, h') <1$, hence by Claim 2 $h$ and $h'$ are $Y$-close. To show that $h$ and $h'$ are also $X$-close, we need to check that all neighboring histories $e,e'\in [h,h']$  are $X$-close. For this, let $e\in [h,h']$ with $last(e)=s$ and  $e'=(e, Z^\beta, s')\in [h,h']$ be such neighboring histories. By part 5, we have
$d_Y^W(e,e')\leq d_Y^W(h,h')<1$. By Claim 2, this implies that we are in one of the following three cases: either $\bM\models U(Z;Y)$ and $\beta\in (0,1)$, or $s\models K_ZY$ and $\beta\in (0,\delta^e_Y)$, or else $s\models D_ZY$ and $\beta=0$. Using the assumption that  $\bM\models U(Y;X)$ (which implies that $s\models K_YX$ and $s\models D_YX$) and the transitivity axioms for $U,K$ and $D$, as well as the fact that $\delta^{e_Y}\leq \delta^{e_X}$ (due to Claim 1.(5) and $\bM\models U(Y;X)$) we can easily see that each of these cases yields the analogous case for $X$: in other words, we have either $\bM\models U(Z;X)$ and $\beta\in (0,1)$, or $s\models K_ZX$ and $\beta\in (0,\delta^e_X)$, or else $s\models D_ZX$ and $\beta=0$. In all cases, $h$ and $h'$ are $X$-close, as desired.

For part 12, let $h$ and $h'$ be histories with $last(h)\models_\bM K_YX$ and $d_Y^W(h,h')< \delta^h\leq 1$, hence $d_Y^W(h,h')=\beta<\delta^h$, and we can also assume that $\bM\not\models U(Y;X)$ (otherwise the desired conclusion follows immediately by part 11). Applying again Claim 2, $h$ and $h'$ are $Y$-close. To show that $h$ and $h'$ are also $X$-close, we need to check that all neighboring histories $e,e'\in [h,h']$ histories are $X$-close. For this, let $e\in [h,h']$ and  $e'=(e, Z^\beta, s')\in [h,h']$ be such neighboring histories. Since $d_Y^W(h,h')<\delta^h$, by part 9 we have $h'\in \mathcal{B}_Y(h, \delta^h)\subseteq h\uparrow_X$, i.e.
$h\leq_Y^W h'$. From this and the fact that $e\in [h,h']$, we obtain that $h\leq_Y^W e$, which together with $last(h)\models K_YX$, gives us that $last(e)\models K_YX$ and $h\leq_X^W e$ (by Claim 1.(15) in Section \ref{preord-unravel}), and so $\delta^{h_X}=\delta^{e_X}$ by Claim 1.(4) above. Since $h,h'$ are $Y$-close and $[e,e'\in [h,h']$, we have that $e,e'$ are $Y$-close. By the definition of $Y$-closeness and the assumption that $\bM\not\models U(Y;X)$, it follows that we are in one of the following two cases: either $last(e)\models K_ZY$ and $\beta< \delta^{e_Y}$, or else $last(e)\models D_ZY$ and $\beta=0$. In the first case, from $last(e)\models K_ZY$ and $last(e)\models K_YX$ we obtain $last(e)\models K_ZX$, which combined with $\beta<\delta^h\leq \delta^{h_X}=\delta^{e_X}$ (due to parts 3 and 4 of Claim 1 and the fact that $h\leq_X^W e$) gives
us that $e,e'$ are $X$-close, as desired. In the second case, in which we have $last(e)\models D_ZY$ and $\beta=0$, we combine these with the
fact that $last(e)\models K_YX$ implies $last(e)\models D_YX$ and use additivity to conclude that we have $last(e)\models D_ZX$ and $\beta=0$, thus $e,e'$ are again $X$-close, as desired.
\end{proof}

\medskip

The facts in Claim 4 are the main ingredients ensuring that our construction will work as intended. We now proceed to show this, in Claims 5-8, that wrap up the proof of our representation result.

\begin{quote}
\textbf{Claim 5}. Each $d_X^W$ is an ultra-pseudo-metric whose topology is a refinement of the Alexandroff topology given by the corresponding preorder $\leq_x^W$ and whose metric-topological indistinguishability relation (relating states $s,w$ s.t. $d_X^W(s,w)=0$) coincides with the $X$-equality relation $=_X^W$. Thus, $\bM^\infty_d$ is a refinement of the (standard) topo-model $\bM^\infty_\leq$ to a standard pseudo-metric model. Moreover, the metric dependence model associated to $\bM^\infty_d$ is pseudo-locally Lipschitz (with Lipschitz constant $1$), and so $\bM^\infty_d$ validates the `paradisiacal' implication $k_XY\Rightarrow U_XY$.
\end{quote}

\begin{proof} The fact that $d_X^W$ is an ultra-pseudo-metric follows immediately by putting together items 1, 2 and 6 of Claim 4. By part 9 of Claim 4, we have that every upset $h\uparrow_X=\{h'\in W: h\leq_X^W h'\}$ includes an $X$-open ball centered at $h$ (namely the ball $\mathcal{B}_X(h, \delta^h)=\{h'\in W: d_X^W(h,h')< \delta^h\}$), which implies that the pseudo-metric topology induced by $d_X^W$ is a refinement of the preorder topology induced by the preorder $\leq_X^W$. Part 3 of Claim 4 shows that the metric indistinguishability relation coincides with $X$-equality $=_X^W$. Standardness is by the definition of $d_X^W$ as the Chebyshev distance based on the $d_x^W$'s. Finally, the fact that the associated metric dependence model is pseudo-locally Lipschitz (with constant $K=1$) follows immediately from the (second part of the) Important Observations following Claim 2. \end{proof}

\medskip
\medskip

\begin{quote}
\textbf{Claim 6}. $h \models_{\bM^\infty_d} K_YX$ iff $last(h)\models_{\bM} K_YX$.
\end{quote}

\begin{proof} We can assume that $X\neq\emptyset$, otherwise both sides hold by definition.

For the left-to-right direction, suppose that $h\in W$ is a history with $h\models_{\bM^\infty_d} K_YX$. We will only use a \emph{weaker consequence} of this, namely the fact that $h\models_{\bM^\infty_d} k_YX$ (\emph{pointwise continuous} dependence at $h$).\footnote{But note that by the last part of Claim 5, this property is actually equivalent to $K_XY$ on the `paradisiacal' model $\bM^\infty_d$.}  Choose any $\eps\in (0,1)$. Then, by the continuous dependence $k_YX$ at point $h$, there exists some $\delta>0$ s.t.
$$(**) \,\,\,\,\,\,\,\,\,\, \mathcal{B}_Y (h,\delta)\subseteq \mathcal{B}_X(h,\eps).$$
Take now the history $h':=(h, Y^\beta, s)$, where $s:=last(h)$ and $\beta$ is any chosen number with $0<\beta < \delta$. Since
$d_Y^W(h,h')=\beta< \delta$, we obtain by $(**)$ that $d_X^W(h,h')<\eps<1$. Given the definition of $d_X$ on neighboring histories $h$ and $h'$ and the assumption that $X\neq\emptyset$, we can infer from $d_X(h,h')<1$ and $\beta\neq 0$ that we have $s=last(h)\models_\bM K_YX$, as desired.

For the converse, let $h\in W$ be a history with $s:=last(h)\models_\bM K_YX$. We need to show that $h\models_{\bM^\infty_d} K_YX$.
For this, we choose $\delta_0:=\delta^h$, and have to prove that $X$ is continuously dependent on $Y$ on the open neighborhood $\mathcal{B}_Y(h,\delta_0)$. In fact, we will prove a \emph{stronger} statement, namely that  $X$ is \emph{uniformly} continuously dependent on $Y$ on this neighborhood.\footnote{Once again, note that by Claim 5 the apparently stronger condition $U_XY$ is actually equivalent to $K_XY$ on $\bM^\infty_d$.} Let $\eps>0$ be any arbitrary positive number. To show uniformly continuous dependence on our given neighborhood, we need to show that
$$(?) \,\,\,\,\,\,\,\,\,\, \exists \delta>0 \mbox{ s.t. } \forall h', h'' \in \mathcal{B}_Y(h,\delta_0) \,  \left(\, d_Y(h',h'') < \delta\Rightarrow d_X^W(h', h'') < \eps\, \right).$$
For this, we choose $\delta:=  min( \delta_0, \eps)=min (\delta^h, \eps)$, so we have $\delta\leq \delta^h, \eps$. Let  $h', h''\in \mathcal{B}_Y(h,\delta_0)$ be s.t. $d_Y(h',h'') < \delta$; so we have  $d_Y^W(h,h'), d_Y^W(h,h'')< \delta^h$ (since $\delta_0=\delta^h$ and $\delta\leq \delta^h$). From these, together with our assumption that $last(h)\models_\bM K_YX$, we infer (by part 12 of Claim 4) that $h$ is $X$-close to both $h'$ and $h''$. Since closeness is an equivalence relation, it follows that $h'$ and $h''$ are $X$-close. By the Important Observation following Claim 2, we obtain that $d_X^W(h', h'')=d_Y^W(h', h'')<\delta\leq \eps$, as desired. \end{proof}

\bigskip

Next, we prove the analogous assertion for uniform dependence:

\begin{quote}
\textbf{Claim 7}. $\bM^\infty_d \models U(Y;X)$ iff $\bM\models U(Y;X)$.
\end{quote}

\begin{proof}
For the left-to-right direction, suppose that $\bM^\infty_d \models U(Y;X)$, and assume towards a contradiction that $\bM\not\models U(Y;X)$.
Choose some number $\eps\in (0,1)$. Applying the definition of $\bM^\infty_d \models U(Y;X)$ to $\eps$, it follows that there exists some number $\beta\in (0,1)$ with the property that
$$(*) \,\,\,\,\,\,\,\,\,\, \forall h,h'\in W(d_Y(h,h') < \beta \Rightarrow d_X(h,h') < \eps).$$
Fix now some $\gamma\in (0,\beta)$, and consider the histories  $h:=(s_0, \emptyset^\gamma, s_0)$ and $h':=(h, Y^\gamma, s_0)=(s_0, \emptyset^\gamma, s_0, Y^\gamma, s_0)$.  By Claim 3.(2), we have $d_Y(h, h')=\gamma<\beta$, so by $(*)$ we have $d_X(h,h')< \eps<1$, hence $d_X(h, h')\neq 1$, and thus (by Claim 2) $h$ and $h'$ are $X$-close. From this, together with the assumption that $\bM\not \models U(Y;X)$ and the fact that the history $h'=(h, Y^{\gamma}, s_0)$ has $\gamma\neq 0$, we obtain (using the definition of $X$-closeness) that $s_0=last(h)\models K_YX$ and $\beta< \delta^{h_X}$. We consider now two cases: $(s_0)\leq_X^W h$ and $(s_0)\not\leq_X^W h$. If $(s_0)\leq_X^W h$, then in fact we must have $(s_0)\leq_X^W h$ (since $(s_0)$ and $h$ are neighboring histories), so by the definition of $\to_X^\leq$
we infer that $s_0\models K_\emptyset X$, and hence $s_0\models U(Y;X)$ (by the axiom $K(X)\Rightarrow U(Y;X)$ of the system $\mathbf{LUD}$), thus contradicting our assumption that $\bM\not\models U(Y;X)$.
In the alternative case when $(s_0)\not\leq_X^W h$, it follows that $h=(s_0, \emptyset^\gamma, s_0)$ is its own $X$-root (since it contains no proper sub-history $h''\prec h$ with $h''\leq_X^W h$), i.e., $h=h_X$, and thus we infer $\delta^{h_X}=\delta^h=\gamma<\beta< \delta^{h_X}$,
 obtaining the contradiction $\delta^{h_X}< \delta^{h_X}$.

\smallskip

For the converse direction, suppose that $\bM\models U(Y;X)$. To prove that $\bM^\infty_d\models U(Y;X)$, let $\eps$ be any number with $0<\eps<1$, and we need to show that there exists $\delta>0$ s.t.
$$\forall h, h'\in W (d_Y^W(h, h') < \delta \Rightarrow d_X^W(h,h') < \eps).$$
For this, we take $\delta:=\eps$. To prove the desired conclusion, let $h,h'\in W$ be histories with $d_Y^W(h, h') < \delta=\eps<1$. By part 11 of Claim 4 and the fact that $\bM\models U(Y;X)$, this implies that $h$ and $h'$ are $X\cup Y$-close, and thus (by the Important Observation following Claim 2) we have $d_X^W(h,h')=d_Y^W(h,h') <\delta=\eps$, as desired. \end{proof}

\medskip

\begin{quote}
\textbf{Claim 8}. The map $last: W \to S$, mapping each history $h\in W$ to its last element $last(h)\in S$, is a uniform topo-morphism from $\bM^\infty_d$ to $\bM$.
\end{quote}

\begin{proof} In the proof of Lemma \ref{rep-preord} (our representation result for preorder models), we showed that the map $last$ is a surjective topo-morphism from $\bM^\infty_\leq$ to $\bM$. In particular, using the atomic clauses for $P\vec{x}$ and $D_XY$ in the definition of that topo-morphism, together with the fact that the models $\bM^\infty_d$ and $\bM^\infty_\leq$ agree on the semantics of these atoms (which, in the case of $D_XY$ follows from the fact that both $\bM^\infty_d$ and $\bM^\infty_\leq$ are standard and that by Claim 5 the pseudo-metric indistinguishability coincides with $X$-equality $=_X^W$), we obtain the corresponding atomic clauses for $P\vec{x}$ and $D_XY$ in the definition of topo-morphism from $\bM^\infty_d$ to $\bM$.

The atomic clauses for $K_XY$ and $U(X;Y)$ follow from Claims 6 and 7.

The fact that the map $last$ is a $p$-morphism wrt $=_X$ was already checked in the proof of Lemma \ref{rep-preord}. Since Claim 5 tells us that the relation $=_X^W$ is the same as the pseudo-metric indistinguishability relation, it follows that $last$ is a topo-morphism with respect to this relation.

Finally, we need to prove that $last$ is an interior map, i.e. both open and continuous, when seen as a map from $\bM^\infty_d$ to $\bM$.
For this, we will make use of the fact that (by the proof of Lemma \ref{rep-preord}), $last$ is an interior map (hence continuous) when seen as a map from $\bM^\infty_\leq$ to $\bM$.

To show that $last$ is open, let $U$ be a open set in $\tau_X^W$. We need to prove that $last(U)$ is upward-closed in $S$ wrt $\leq_X$. For this, let $last(h)\in last(U)$ with $h\in U$, and we have to show that $last(h)\uparrow_X$ is included in $last(U)$. Since $h\in U$ and $U$ is open in the pseudo-metric topology induced by $d_X^W$, there must exist some $\eps>0$ s.t. $\mathcal{B}_X(h,\eps)\subseteq U$. By part 10 of Claim 4, we obtain $last(h)\uparrow_X\subseteq last(\mathcal{B}_X(h,\eps))\subseteq last(U)$, as desired.

To show that $last$ is continuous, let $V\subseteq S$ be upward-closed in $S$ wrt $\leq_X$. We need to prove that $last^{-1}(V)$ is open in the pseudo-metric topology $\tau_X^W$. Since the map $last$ was shown in the previous section to be a $p$-morphism from $\bM^\infty_\leq$ to $\bM$ (or equivalently, an interior map between the two preorder topologies), it follows that $last^{-1}(V)$ is upward-closed in $W$ wrt $\leq_X^W$ (i.e., open in the Alexandroff topology induced by $\leq_X^W$). But our pseudo-metric topology $\tau_X^W$ is a refinement of the Alexandroff topology induced by $\leq_X^W$, and therefore $last^{-1}(V)$ is also open in the pseudo-metric topology. \end{proof}

\bigskip

This finishes the proof of our representation theorem, hence also our completeness and decidability proof for $\mathbf{LUD}$ with respect to metric dependence models. Since the model we constructed is pseudo-locally Lipschitz, this also proves completeness for pseudo-locally Lipschitz models.

\end{document}